\documentclass[leqno,a4paper,12pt]{article}
\usepackage{amsmath,amsthm,bbm}
\usepackage[sans]{dsfont}

\newtheorem{Thm}{\indent Theorem}[section]
\newtheorem{MThm}[Thm]{\indent Main Theorem}
\newtheorem{Prop}[Thm]{\indent Proposition}
\newtheorem{Lem}[Thm]{\indent Lemma}
\newtheorem{Cor}[Thm]{\indent Corollary}

\newtheorem{Var}[Thm]{\indent Variant}

\theoremstyle{definition}
\newtheorem{Def}[Thm]{\indent Definition}
\newtheorem{Rem}[Thm]{\indent Remark}
\newtheorem{Ex}[Thm]{\indent Example}

\newtheorem{Aux}[Thm]{\indent Auxiliary Construction}
\newtheorem{Cons}[Thm]{\indent Construction}

\newtheorem{Scho}[Thm]{\indent Scholie}

\def\qed{{\hskip0pt\unskip\unskip\nobreak\hfil\penalty50
          \hskip1em\hbox{}\nobreak\hfil
          {\bf q.e.d.}%
          \parfillskip=0pt\finalhyphendemerits=0
          \par}\medskip}

\newenvironment{Proof}
               {{\it Proof.}\quad}
               {\qed}

\newenvironment{Proofof}[1]
               {{\it Proof of #1.}\quad}
               {\qed}

\newcommand{\Prime}{\kern3\fontdimen1\font$'$\kern-7\fontdimen1\font}

\long\def\forget#1{}

\long\def\beginSIDEREMARK#1\endSIDEREMARK
    {{\par\bigskip\advance\leftskip by 2cm
                  \advance\rightskip by -2cm\noindent
      {\bf Our own side remark:} #1
      \par\bigskip\noindent}}

\long\def\beginFORGET#1\endFORGET{#1}
\long\def\beginFORGET#1\endFORGET{}

\def\?{\ ???\ \immediate\write16{}%
\immediate\write16{Warning: There was still a question mark . . . }%
\immediate\write16{}}

\usepackage{amsmath}

\usepackage{exscale}
\usepackage{amssymb}

\input cyracc.def

\newcommand{\BA}{{\mathbb{A}}}

\newcommand{\BC}{{\mathbb{C}}}

\newcommand{\BG}{{\mathbb{G}}}

\newcommand{\BI}{{\mathbb{I}}}

\newcommand{\BQ}{{\mathbb{Q}}}
\newcommand{\BR}{{\mathbb{R}}}

\newcommand{\BV}{{\mathbb{V}}}
\newcommand{\BW}{{\mathbb{W}}}
\newcommand{\BX}{{\mathbb{X}}}

\newcommand{\BZ}{{\mathbb{Z}}}

\newcommand{\Fc}{{\mathfrak{c}}}

\newcommand{\Fn}{{\mathfrak{n}}}

\newcommand{\FH}{{\mathfrak{H}}}

\newcommand{\FX}{{\mathfrak{X}}}

\newcommand{\CA}{{\cal A}}

\newcommand{\CC}{{\cal C}}
\newcommand{\CD}{{\cal D}}

\newcommand{\CF}{{\cal F}}
\newcommand{\CG}{{\cal G}}
\newcommand{\CH}{{\cal H}}

\newcommand{\CK}{{\cal K}}
\newcommand{\CL}{{\cal L}}

\newcommand{\CR}{{\cal R}}
\newcommand{\CS}{{\cal S}}

\newcommand{\CV}{{\cal V}}
\newcommand{\CW}{{\cal W}}

\newcommand{\ad}{\mathop{\rm ad}\nolimits}
\newcommand{\adm}{\mathop{\rm adm}\nolimits}

\newcommand{\App}{\mathop{\rm App}\nolimits}

\newcommand{\coh}{\mathop{\rm coh}\nolimits}

\newcommand{\Loc}{\mathop{\rm Loc}\nolimits}

\newcommand{\Cent}{\mathop{\rm Cent}\nolimits}

\newcommand{\GL}{\mathop{\rm GL}\nolimits}
\newcommand{\Gm}{\mathop{\BG_m}\nolimits}

\newcommand{\Rad}{\mathop{\rm Rad}\nolimits}
\newcommand{\Rep}{\mathop{\rm Rep}\nolimits}
\newcommand{\Res}{\mathop{\rm Res}\nolimits}
\newcommand{\res}{\mathop{\rm res}\nolimits}

\newcommand{\Stab}{\mathop{\rm Stab}\nolimits}

\newcommand{\sgn}{\mathop{{\rm sgn}}\nolimits}

\newcommand{\loccit}{[loc.$\;$cit.]}

\def\tei{\, | \,}
\def\halb{\frac{1}{2}}

\def\id{{\rm id}}

\newbox\mybox
\def\arrover#1{\mathrel{
       \setbox\mybox=\hbox spread 1.4em{\hfil$\scriptstyle#1$\hfil}
       \vbox{\offinterlineskip\copy\mybox
             \hbox to\wd\mybox{\rightarrowfill}}}}
\def\larrover#1{\mathrel{
       \setbox\mybox=\hbox spread 1.4em{\hfil$\scriptstyle#1$\hfil}
       \vbox{\offinterlineskip\copy\mybox
             \hbox to\wd\mybox{\leftarrowfill}}}}

\def\ontoover#1{\mathrel{
       \setbox\mybox=\hbox spread 1.4em{\hfil$\scriptstyle#1$\hfil}
       \vbox{\offinterlineskip\copy\mybox
             \hbox to\wd\mybox{\rightarrowfill\hskip-2.8mm
                               $\rightarrow$}}}}
\def\leftontoover#1{\mathrel{
       \setbox\mybox=\hbox spread 1.4em{\hfil$\scriptstyle#1$\hfil}
       \vbox{\offinterlineskip\copy\mybox
             \hbox to\wd\mybox{$\leftarrow$\hskip-2.8mm
                               \leftarrowfill}}}}
\def\longto{\longrightarrow}
\def\into{\hookrightarrow}
\def\onto{\ontoover{\ }}
\def\longonto{\ontoover{\ }}
\def\isoto{\arrover{\sim}}

\def\longinto{\lhook\joinrel\longrightarrow}

\usepackage[curve,matrix,arrow,cmtip]{xy}
\NoComputerModernTips

\def\myxymessage{\def\messagetext
   {Here an xy-pic diagram was omitted to speed up compilation . . . }
   \immediate\write16{\messagetext}
   \hbox{\bf \messagetext}}
\def\filxymatrix#1{\myxymessage}
\def\filxyarray#1{\myxymessage}

\newdir^{ (}{{}*!/-3pt/\dir^{(}}
\newdir_{ (}{{}*!/-3pt/\dir_{(}}
\newdir^{ )}{{}*!/+3pt/\dir^{)}}
\newdir_{ )}{{}*!/+3pt/\dir_{)}}

\def\rscript#1{\hbox to 0pt{$\scriptstyle#1$\hss}}

\let\oldbullet\bullet
\def\bullet{{\mathchoice{\oldbullet}%
                        {\oldbullet}%
                        {\scriptscriptstyle\oldbullet}%
                        {\oldbullet}}}

\newcommand{\argdot}{{\;\bullet\;}}
\newcommand{\argast}{{\;\ast\;}}

\newcommand{\Ab}{\mathop{\CA b}\nolimits}

\newcommand{\tQ}{\mathop{\tilde{Q}}\nolimits}

\newcommand{\uQ}{\mathop{\underline{Q}}\nolimits}
\newcommand{\utQ}{\mathop{\underline{\tQ}}\nolimits}

\begin{document}

%

\hfuzz=3pt
\overfullrule=10pt                   


\setlength{\abovedisplayskip}{6.0pt plus 3.0pt}
\setlength{\belowdisplayskip}{6.0pt plus 3.0pt}
\setlength{\abovedisplayshortskip}{6.0pt plus 3.0pt}
\setlength{\belowdisplayshortskip}{6.0pt plus 3.0pt}

\setlength{\baselineskip}{13.0pt}
\setlength{\lineskip}{0.0pt}
\setlength{\lineskiplimit}{0.0pt}

%
%

\title{Shimura data and corners: cohomology
\forget{
\footnotemark
\footnotetext{To appear in ....}
}
}
\author{\footnotesize by\\ \\
\mbox{\hskip-2cm
\begin{minipage}{6cm} \begin{center} \begin{tabular}{c}
J\"org Wildeshaus \\[0.2cm]
\footnotesize Universit\'e Sorbonne Paris Nord \\[-3pt]
\footnotesize LAGA, CNRS (UMR~7539)\\[-3pt]
\footnotesize F-93430 Villetaneuse\\[-3pt]
\footnotesize France\\
{\footnotesize \tt wildesh@math.univ-paris13.fr}
\end{tabular} \end{center} \end{minipage}
\hskip-2cm}
\\[2.5cm]
}
\maketitle
\begin{abstract}
\noindent
The purpose of this article is to determine the gluing data associated to degeneration 
of local systems to the boundary of the Baily--Borel compactification of a Shimura variety.  \\


\end{abstract}


\bigskip
\bigskip
\bigskip

\noindent {\footnotesize Math.\ Subj.\ Class.\ (2020) numbers: 
14G35 (11F75, 14F25, 54B40, 54C10).
}

\eject

\tableofcontents

\bigskip

%
%

\setcounter{section}{-1}
\section{Introduction}
\label{Intro}



Let $M$ be a \emph{Shimura variety}, and $M^*$ its \emph{Baily--Borel compactification}
\cite{AMRT,BB}. The aim of this article is to provide tools allowing to control \emph{cohomo\-logy}
of the \emph{boundary} $\partial M^* := M^* - M$ of $M^*$. These tools are based on our previous analysis
of the proper, surjective map from the \emph{Borel--Serre
compactification} $M^{BS}$ \cite{BS} to (the space of complex points of) $M^*$ \cite{W1}. \\

Our tools are of sheaf theoretic origin.
In order to discuss their precise nature, it is convenient to consider a finite filtration
$\emptyset = U_{r+1} \subset U_r \subset \ldots \subset U_2 \subset U_1 = X$ of a 
compact topological space
$X$ by open sub-sets $U_p \,$, $p=1, \ldots, r+1$. Denote by $j_p$ the open immersion of $U_p \,$, and by $i_p$ the immersion
of the stratum $M_p := U_p - U_{p+1}$ into $X$, $p = 1, \ldots r$. Thus, $i_1$ is closed,
and $i_r$ is open. For any complex $\CF$ of sheaves 
of Abelian groups on $X$, that is bounded from below, we have a convergent spectral sequence
\[
(E_X) \quad\quad\quad\quad\quad\quad E^{p,q}_1 = H_c^{p+q} (M_p, i_p^* \CF) \Longrightarrow H^{p+q}(X,\CF) \; ,
\]
relating cohomology of $X$ to cohomology with compact supports of the strata $M_p \,$, $p = 1, \ldots r$. \\

A first sheaf theoretic approximation of ``control'' of cohomology of $X$ would be to 
describe the restrictions $i_p^* \CF$ of $\CF$ to the individual strata $M_p \,$. 
\forget{Ingredient~(A) appears necessary for a description of the $E_1$-terms of
our spectral sequence.} \\

Depending on which aspect of cohomology of $X$ one wishes to control, 
it might be sufficient to control the $E_1$-terms of the spectral sequence $(E_X)$. 
A typical such aspect is the validity of a property~(P) on Abelian groups, stable under passage
to sub-quotients and extensions. In order to establish (P) for $H^{\! \argdot \!}(X,\CF)$,
it suffices to do so for all $E_1$-terms $H_c^{p+q} (M_p, i_p^* \CF)$. 
For an example where this principle is exploited in the context of Shimura varieties,
we refer to \cite[Sect.~8.7]{MT} (see \cite[Prop.~6]{MT} in particular). \\

In general, property~(P) is satisfied for $H^{\! \argdot \!}(X,\CF)$
if and only if it is sa\-tisfied for all $E_\infty$-terms of the spectral sequence $(E_X)$. But if $(E_X)$
does not dege\-ne\-rate at $E_1$, property~(P)
might be true for $H^{\! \argdot \!}(X,\CF)$, without holding  
for all $E_1$-terms. Hence the need to control the $d^{p,q}_1$,
and \emph{a priori}, all higher differentials $d^{p,q}_n: E^{p,q}_n \to E^{p+n,q-n+1}_n$ as well. \\

The differential $d^{p,q}_1: E^{p,q}_1 \to E^{p+1,q}_1$ is the composition of the boundary homomorphism
\[
E^{p,q}_1 = H_c^{p+q} (M_p , i_p^* \CF) \longto H_c^{p+q+1} (U_{p+1}, j_{p+1}^* \CF)
\]
and the restriction
\[
H_c^{(p+1)+q} (U_{p+1}, j_{p+1}^* \CF) \longto H_c^{(p+1)+q} (M_{p+1} , i_{p+1}^* \CF) = E^{p+1,q}_1 \; .
\]
It is induced, by applying the cohomological functor $H_c^{\! \argdot}$, from a morphism
\[
(i_p)_! i_p^* \CF \longto (i_{p+1})_! i_{p+1}^* \CF [1]
\]
in the derived category of complexes of sheaves of Abelian groups on $X$. The latter corresponds to
the morphism
\[
i_p^* \CF \longto Ri_p^! (i_{p+1})_! i_{p+1}^* \CF [1] = i_p^* R(i_{p+1})_* i_{p+1}^* \CF \; ,
\]
that is equal to $i_p^*(\ad)$, where $\ad: \CF \to R(i_{p+1})_* i_{p+1}^* \CF$ equals the adjunction. \\

A more precise sheaf theoretic approximation of ``control'' of cohomology of $X$ would be to describe the 
family of morphisms 
\[
i_p^*(\ad): i_p^* \CF \longto i_p^* R(i_{p+1})_* i_{p+1}^* \CF \; .
\]
\forget{Ingredient~(B) appears necessary for a description of the differentials 
relating the $E_1$-terms of our spectral sequence.} 

In a somewhat more elaborate way than what has just been recalled, 
the differentials $d^{p,q}_n: E^{p,q}_n \to E^{p+n,q-n+1}_n$, for $1 \le n \le r$, are related to
the morphisms 
$i_p^*(\ad_{\alpha,\beta})$, where $\alpha \subset \beta$ are sub-sets of $\{ p, p+1, \ldots, p+n \}$
containing $p$,
$I_\alpha \CF$ is the image of $\CF$ under the composition of the functors $R(i_m)_* i_m^*$, for
$m \in \alpha$, and $\ad_{\alpha,\beta}$ is the composition of adjunctions 
\[
I_\alpha \CF \longto I_\beta \CF \; .
\]
For us, sheaf theoretic ``control'' of cohomology of $X$ means to describe the family of morphisms
\[
i_p^*(\ad_{\alpha,\beta}): i_p^* I_\alpha \CF \longto i_p^* I_\beta \CF \; .
\]
Let us refer to this family of morphisms as the \emph{gluing data} for the restrictions $i_p^* \CF$. \\

In the situation of interest for us, the space $X$ equals (the space of complex points of)
the boundary $\partial M^*$ of $M^*$, the ``level'' 
of $M$, \emph{i.e.}, the underlying open compact group $K$, is supposed \emph{neat} \cite[Sect.~0.6]{P},
the strata $M_p$ are those of the \emph{canonical stratification} of $\partial M^*$ (see \emph{e.g.}
\cite[Def.~5.8]{W1}), and $\CF$ comes about as the restriction 
$(R j_* \CV)_{\tei \partial M^*}$ to $\partial M^*$ of the direct image 
$R j_* \CV$ of a complex of sheaves $\CV$ on $M$, whose cohomology objects are local systems
($j:=$ the open immersion of $M$ into $M^*$). 
Our Main Theorem~\ref{7MT} gives a description of their gluing data $(i_p^*(\ad_{\alpha,\beta}))_{p,\alpha \subset \beta} \, $.
This involves in particular a description of the $i^*_p I_\alpha (R j_* \CV)_{\tei \partial M^*} \,$: 
for the singleton $\alpha = \{ p \}$, we recover the well-known result (see \emph{e.g.} \cite[(6.2)]{LR}) allowing to compute 
$i_p^* (R j_* \CV)_{\tei \partial M^*}$
as \emph{cohomology} of the inertia group $H_p$ of the stratum $M_p \, $. If 
more generally $\alpha$ contains $p$, then
$i^*_p I_\alpha (R j_* \CV)_{\tei \partial M^*}$ turns out to be 
computed as \emph{par-cohomology} of $H_p \,$, by which we mean ---
for the purpose of this introduction only ---
a certain direct sum of cohomology of
intersections $H_{\alpha}$ of $H_p$ with parabolic sub-groups (which depend on the elements of $\alpha$, that are distinct from $p$)
of the ambiant algebraic group. While this observation can be deduced, at least in specific cases and using
isomorphisms, whose canonicity is not obvious, from successive application of the case of a singleton (whose
elements run through those of $\alpha$, in decreasing order), our identification of the $i_p^*(\ad_{\alpha,\beta})$,
although easy to state, appears new: the morphisms 
\[
i_p^*(\ad_{\alpha,\beta}): i_p^* I_\alpha \CF \longto i_p^* I_\beta \CF \; ,
\]
for $\alpha \subset \beta$, are given by restriction $\Res^{H_{\alpha}}_{H_{\beta}}$ from $H_{\alpha}$ to $H_{\beta} \,$! \\

The main application we have in mind concerns complexes $\CV$ of sheaves underlying \emph{variations
of Hodge structure} on $M$. Then 
\[
(E_{\partial M^*}) \quad\quad
E^{p,q}_1 = H_c^{p+q}  \bigl( M_p, i_p^* (R j_* \CV)_{\tei \partial M^*} \bigr) 
\Longrightarrow H^{p+q} \bigl( \partial M^*,(R j_* \CV)_{\tei \partial M^*} \bigr) 
\]
is a spectral sequence of \emph{mixed Hodge structures}. In view of our theory of the \emph{interior}, or 
\emph{intersection motive} (see in particular \cite[Sect.~4]{W4},
\cite[Def.~2.4]{W5}), we are interested,
for specific choices of $\CV$, in controlling the \emph{weights} occurring in
$H^{\! \argdot}( \partial M^*,(R j_* \CV)_{\tei \partial M^*})$.
In this context, properties~(P) of interest are on the one hand, the absence of a given weight
in a Hodge structure, and on the other hand, its presence. According to which pro\-perty~(P) one treats,
different issues concer\-ning
the spectral sequence $(E_{\partial M^*})$ occur: for example, weight $w$ might occur in a certain term
$E^{p,q}_1$. But in order to ensure its presence in $H^{p+q}( \partial M^*,(R j_* \CV)_{\tei \partial M^*})$,
one needs to prove that $w$ ``survives'' until $E^{p,q}_\infty$. 
This precise problem occurs \emph{e.g.} in the context of Siegel threefolds (see \cite[proof of Prop.~2.9]{W2}).
It will be reinterpreted using the description of the gluing data
from Main Theorem~\ref{7MT}, at the end of Section~\ref{7} (see Example~\ref{7R} and Remark~\ref{7S}). \\

The dual approach to control $H^{\! \argdot}( \partial M^*,(R j_* \CV)_{\tei \partial M^*})$ necessitates information
on the exceptional inverse images $i_p^! (R j_* \CV)_{\tei \partial M^*}$. 
Our Main Theorem~\ref{8MT} provides that information, together with
a description of the cano\-ni\-cal morphisms 
\[
i_p^! (R j_* \CV)_{\tei \partial M^*} \longto i_p^* (R j_* \CV)_{\tei \partial M^*} \; ,
\]
and their cones. In addition to cohomology, that description uses 
\emph{cohomo\-lo\-gy with compact supports} and \emph{boundary cohomology} of the inertia group $H_p \,$. \\

The proofs of both Main Theorem~\ref{7MT} and Main Theorem~\ref{8MT} dictate the structure
of this paper. They use proper base change to the
Borel--Serre compactification $M^{BS}$, meaning that both correspond, \emph{via}
application of the higher direct image, to results (Main Theorem~\ref{4MT} and Main Theorem~\ref{5MT}, respectively)
concerning the boundary of $M^{BS}$. The first three sections
prepare these results. Sections~\ref{4} and \ref{5} contain their formulation and their proofs.
Section~\ref{6} makes the application of the higher direct image to Main Theorems~\ref{4MT} and \ref{5MT} explicit,
yielding Theorems~\ref{6MT1} and \ref{6MT2}. Section~\ref{7} translates Theorem~\ref{6MT1} into
the language of group cohomology. Section~\ref{8} translates Theorem~\ref{6MT2} into
group cohomology and its variants: cohomology with compact supports and boundary cohomo\-lo\-gy. \\ 

Let us now give a more detailed description of the individual sections. \\

Section~\ref{1} starts with a key notion, namely, 
that of a \emph{contractible map} of topological spaces (Definition~\ref{1A}). 
Its interest stems from the fact that the open immersion of the Shimura variety $M$
into its Borel--Serre compactification $M^{BS}$ is contractible. More generally, this holds
for the inclusion of any of the strata of $M^{BS}$ into its closure. 
Given a contractible map $f: A \to B$, and a local system $\CF$ on $A$, 
the complex $R f_* \CF$ is concentrated in degree zero (Proposition~\ref{1B}). For any local 
system $\CG$ on $B$, the adjunction $\CG \longto Rf_* f^* \CG$ is an isomorphism
(Corollary~\ref{1C}). This latter result will be the essential ingredient of the proof of Main Theorem~\ref{4MT}.
If $f$ is a contractible open immersion, with closed complement $k: Z \into B$, then the co-localization triangle
\[
k_*Rk^! \longto  \id_B \longto R f_* f^* \longto k_*Rk^! [1] \; ,
\]
together with Corollary~\ref{1C}, shows that $Rk^!$ vanishes on complexes $\CF$ of sheaves on $B$, that are bounded  
from below, and whose cohomology objects are local systems. Proposition~\ref{1D} provides a refinement
of this observation in a situation where $f$ is only contractible ``up to an error term $Z^0$ contained in $Z$'':
denoting by $\j$ the immersion of $Z^0$ into $Z$, there is a canonical isomorphism between 
$Rk^! \CF$ and $\j_! (\CF_{\tei Z^0})$.
Proposition~\ref{1D} will be the essential ingredient of the proof of Main Theorem~\ref{5MT}. 
The identification of objects in the respective images of $Rk^!$ and $\j_!$ provides a first hint 
why the description of $i_p^! (R j_* \CV)_{\tei \partial M^*}$ should involve    
cohomology with compact supports (Main Theorem~\ref{8MT}). \\

In Section~\ref{2}, we establish the combinatorial set-up for 
our description of the gluing data $(i_p^*(\ad_{\alpha,\beta}))_{p,\alpha \subset \beta} \,$.
In the geometrical context of interest for us, \emph{i.e.}, that of the boundary $\partial M^*$, 
the strata $M_p$ are indexed by proper \emph{admissible} parabolic
sub-groups. We therefore start (Definition~\ref{2A}) by recalling the relation $\preceq \,$, which on the
level of such parabolics corresponds to the relation ``$M_p$ lies in the closure of $M_q$''. 
The indices $\alpha$ of the gluing data then correspond to chains $Q_1 \prec Q_2 \prec \ldots \prec Q_r \,$,
hence the nature of our index set denoted $\CC_{(P,\FX)}$ (Definition~\ref{2D}). 
To finish the section, we give a simple criterion allowing to check
the inclusion ``$\alpha \subset \beta$'' between such indices 
in certain cases (Theorem~\ref{2F}). \\
  
Section~\ref{3} provides the sheaf theoretical enrichment of the set-up created in Section~\ref{2}.
Our main results are comparison statements on the gluing data $(i_p^*(\ad_{\alpha,\beta}))_{p,\alpha \subset \beta} \,$,
hence the need for target categories, in which these compa\-ri\-sons take place. For the result concerning
the Borel--Serre compacti\-fi\-cation $M^{BS}$ (Main Theorem~\ref{4MT}), the target category 
is denoted $\CC_{(G,\FX)}^{K,BS}$ (Definition~\ref{3C}). For the results concerning
the Baily--Borel compactification $M^*$ (Theorem~\ref{6MT1} and Main Theorem~\ref{7MT}), the target category 
is denoted $\CC_{(G,\FX)}^{K,*}$ (Definition~\ref{3E}). The link between $\CC_{(G,\FX)}^{K,BS}$ and $\CC_{(G,\FX)}^{K,*}$ 
is provided by the higher direct image of the map from $M^{BS}$ to $M^*$ (Proposition~\ref{3F}). \\       

In Section~\ref{4}, we give two constructions
of objects in $\CC_{(G,\FX)}^{K,BS}$, starting from complexes $\CF$ of sheaves on the boundary 
$\partial M^{BS} := M^{BS} - M$ of $M^{BS}$, that are bounded  
from below. The first (Construction~\ref{4Con1})
uses successive dege\-ne\-ration along the pre-images $M_p'$ of the strata $M_p$ of $\partial M^*$. 
It anticipates the analoguous construction on the level of the
$M_p \,$, and the use of Proposition~\ref{3F}.
The second (Construction~\ref{4Con2}) is much more straightforward to set up: here,
we use simple pull-back to the $M_p'$. Main Theorem~\ref{4MT} states that
the two constructions yield identical results, provided the cohomology objects of the input data $\CF$ are
local systems. \\

Section~\ref{5} contains the statement and proof of Main Theorem~\ref{5MT}. We are thus concerned with 
of the effect of the exceptional inverse images $R k^!$  on certain complexes $\CF$ of sheaves on $\partial M^{BS}$,
where $k$ runs through the immersions of the pre-images $M_p'$ of the strata $M_p \,$. It turns out
that the hypothesis of Proposition~\ref{1D} is verified only if the index $p$ corresponds to a parabolic, which
is not only admissible, but maximal proper. Under this condition, the identification $Rk^! \CF \cong \j_! (\CF_{\tei Z^0})$
holds, for $\j$ equal to the immersion of the interior $Z^0$ of $M_p'$ (Main Theorem~\ref{5MT}~(a)).
If $p$ corresponds to an $r$-fold intersection of distinct maximal proper parabolics, then we obtain
an isomorphism between $Rk^! \CF$ and the shift $\j_! (\CF_{\tei Z^0})[-(r-1)]$ (Main Theorem~\ref{5MT}~(b)). \\
  
Section~\ref{6} contains Theorems~\ref{6MT1} and \ref{6MT2}, which are the results of application
to Main Theorems~\ref{4MT} and \ref{5MT}, respectively, of the higher direct image associated to the map $M^{BS} \to M^*$. \\

Section~\ref{7} starts by providing the correct framework for par-cohomology of the inertia groups $H_p$
of the strata $M_p$ (Definition~\ref{7A} -- Proposition~\ref{7Da}). 
Our framework takes into account equivariance under the quotient by $H_p$ of a certain larger group $H$
containing $H_p$ as a normal sub-group. This is of importance, since that same quotient $H/H_p$ appears as the fundamental
group of $M_p$ itself; representations thereof thus correspond to local systems on $M_p \,$. 
Definition~\ref{7E} -- Definition~\ref{7B} prepare the statement of Main Theo\-rem~\ref{7MT}. 
In particular, we define the representation theoretical analogue $\CR_{(G,\FX)}^{K,*}$ of the category
$\CC_{(G,\FX)}^{K,*}$ from Definition~\ref{3E}. The correspondence between representations of $H/H_p$ and local systems
of $M_p$ induces a \emph{canonical construction} functor $\CR_{(G,\FX)}^{K,*} \to \CC_{(G,\FX)}^{K,*} \,$, while
the framework developed in the course of the first part of the section 
allows to define a functor from local systems on $M$ to $\CR_{(G,\FX)}^{K,*} \,$.   
Main Theo\-rem~\ref{7MT} then basically asserts that the composition of the two functors,
applied to the complex of sheaves $\CV$ on our Shimura variety $M$,
yields the gluing data associated to $(R j_* \CV)_{\tei \partial M^*} \,$.
We give two variants of Main Theorem~\ref{7MT}, concerning first, 
the case where $\CV$ is given by a complex of algebraic representations (Variant~\ref{7Var1}),
and second, the adelic setting (Variant~\ref{7Var2}). 
We conclude the section with the aforementioned Example~\ref{7R}, where we make the description of the
gluing data from Main Theorem~\ref{7MT} explicit for Siegel threefolds. \\

Section~\ref{8} starts by listing the basic properties of cohomology with compact supports and boundary cohomo\-lo\-gy
of the groups $H_p \,$. In the course of 
the preparation of this paper, it turned out that neither of the two concepts are 
documented in the literature. 
The blog \cite{Bel} helped to realize that there are indeed fundamental issues to be resolved
(not the least subtle being well-definedness). 
In order to avoid overloading an already long work, we decided to address these
issues separately \cite{W3}, and to summarize the ne\-cessary results in Proposition~\ref{8a} -- Proposition~\ref{8B}. 
We then formulate and prove Main Theorem~\ref{8MT}. As in the previous section, we give two variants,
concerning the case of algebraic representations (Variant~\ref{8Var1}),
and the adelic setting (Variant~\ref{8Var2}). We conclude by re-interpreting Example~\ref{7R} using Main Theorem~\ref{8MT}
(Remark~\ref{8V}). \\

I wish to thank F.~Lemma, A.~Mokrane, S.~Morra and J.~Tilouine for useful remarks and discussions. \\

{\bf Conventions}: For an abstract group $H$, let us denote by $\Rep H$ the category of representations of $H$ in Abelian groups. 
Thus, if $H$ is reduced to the singleton, $\Rep H$ equals the category $\Ab$ of Abelian groups. 
For an affine algebraic group $H$ (defined over some base field), and a field of coefficients $F$,  
we shall denote by $\Rep_F H$ the
category of algebraic representations of $H$ in finite dimensional $F$-vector spaces. 
For a topological space $V$, we denote by $D^+(V)$ the derived category  of complexes of sheaves 
of Abelian groups on $V$, that are bounded from below. 

\forget{For an affine algebraic or real Lie group $H$, let us denote by $H^0$ its neutral connected component (in the Zariski topology if $H$ is algebraic, in the ordinary 
topology if $H$ is a real Lie group). Note that if $H$ is affine algebraic over $\BR$, then
$H^0(\BR)$ might be strictly larger than $H(\BR)^0$.

$g$ for elements of $P$ in order not to confuse with the notation for the projection $p$.
Several meanings of $\overline{\ast}$ (overline)...
}


\bigskip

%
%

\section{Cohomology of contractible maps}
\label{1}



Recall the following notion.

\begin{Def}[{\cite[Def.~3.4]{W1}}] \label{1A}
A continuous map $f:A \to B$ of topological spaces
is called \emph{contractible} if the topology on $B$
admits a basis $(V_i)_i$, for which the pre-images
$f^{-1}(V_i) \subset A$ are contractible, for all $i$.
\end{Def}

The property of a continuous map being contractible is local on its target.
If $f$ is contractible, then the image of $f$ is dense in $B$.
Our interest in contractible maps comes from the following fact.

\begin{Prop} \label{1B}
Let $f:A \to B$ be a contractible map, and $\CF$ a locally constant sheaf of Abelian groups
on $A$. \\[0.1cm]
(a)~The sheaf $\CG := R^0f_* \CF$ on B is locally constant. \\[0.1cm]
(b)~The adjunction
\[
\ad: f^* \CG \longto \CF 
\] 
is an isomorphism. \\[0.1cm]
(c)~The object $Rf_* \CF$ of the derived category of sheaves of Abelian groups on $B$
is concentrated in degree zero. \\[0.1cm]
(d)~Up to unique isomorphism, the pair $(\CG,\ad)$ is unique with respect to the requirements
that $\CG$ be locally constant, and $\ad$ an isomorphism $f^* \CG \isoto \CF$.
\end{Prop}

\begin{Proof}
If $U \subset A$ is an open, contractible sub-set, 
then the restriction $\CF_{\tei U}$ is constant, say of value $F$, and
\[ 
H^n \bigl( U, \CF \bigr) = 0
\]
unless $n =0$, in which case the adjunction 
\[
F \longto H^0 \bigl( U',\CF \bigr)
\]
is an isomorphism, for any open, contractible sub-set $U'$ of $U$. 

This latter observation, together with our hypothesis on $f$, implies that for every
open sub-set $V$ of $B$ such that $f^{-1}(V)$ is contractible, 
the restriction $(R^0 f_* \CF)_{\tei V}$ is constant,
and the adjunction 
\[
(f^* R^0 f_* \CF)_{\tei f^{-1}(V)} \longto \CF_{\tei f^{-1}(V)}
\]
is an isomorphism. This implies parts~(a) and (b) of our claim.

As for part~(c), recall first that for any point $b \in B$, and any integer $n$, there is
a canonical isomorphism
\[
\lim_{\longto} H^n \bigl( f^{-1}(V),\CF \bigr) \isoto (R^n f_* \CF)_b \; ,
\]
the direct limit running over all open neighbourhoods $V$ of $b$. Then apply the hypothesis on $f$.

In order to prove part~(d), let $\CG'$ be locally constant on $B$, and $\varphi: f^* \CG' \isoto \CF$.
Then 
\[
R^0 f_*(\varphi) : R^0 f_*f^* \CG' \isoto R^0 f_* \CF = \CG \; .
\]
Composition with the adjunction $\CG' \to R^0 f_*f^* \CG'$ yields a morphism 
\[
\psi: \CG' \longto \CG \; 
\]
satisfying $f^* (\psi) = \ad^{-1} \circ \varphi$. In particular, $f^* (\psi)$ is an isomorphism.
Therefore, kernel and co-kernel of $\psi$ are locally constant sheaves on $B$ whose pull-backs
under $f$ are zero. But the image of $f$ is dense in $B$.
\end{Proof}

\begin{Cor} \label{1C}
Let $f:A \to B$ be a contractible map. \\[0.1cm]
(a)~For any locally constant sheaf $\CF$ of Abelian groups on $A$, the adjunction
\[ 
\ad: f^* R f_* \CF \longto \CF
\]
is an isomorphism. \\[0.1cm]
(b)~For any locally constant sheaf $\CG$ of Abelian groups on $B$, the adjunction
\[ 
\ad: \CG \longto R f_* f^* \CG 
\]
is an isomorphism. \\[0.1cm]
(c)~The functors $R f_*$ and $f^*$ induce mutually inverse
equivalences of the categories of locally constant sheaves 
of Abelian groups on $A$ and on $B$, respectively.
\end{Cor}

Parts~(a) and (b) of Corollary~\ref{1C} generalize to complexes of sheaves,
bounded from below, and
whose cohomology objects are locally constant. \\

In the applications of Corollary~\ref{1C} we have in mind, the map $f$ 
will be an open immersion. 
The contractibility condition will
frequently be satisfied only ``up to an error term'' $Z^0$. Let us be more precise. 
Fix a topological space $X$, two disjoint open sub-sets $U$ and $Z^0$ of $X$,
and consider the following set-up. 
\[
\vcenter{\xymatrix@R-10pt{
        U \ar@{^{ (}->}[r]^-\circ \ar@{=}[d] &
        X - Z^0 \ar@{_{ (}->}[d] &
        Y \ar@{_{ (}->}[l] \ar@{_{ (}->}[d]^{\i} 
        \\
        U \ar@{^{ (}->}[r]_-{j}^-\circ &
        X &
        Z \ar@{_{ (}->}[l]^-{k} 
        \\   
        &
        Z^0 \ar@{=}[r] \ar@{^{ (}->}[u]^-\circ
        &
        Z^0 \ar@{^{ (}->}[u]_-{\j}^-\circ 
\\}}
\] 
Immersions situated on the same line are complementary to each other 
(exam\-ple: $j$ and $k$), the four immersions marked by ``$\circ$''
are open (exam\-ple: $\j$), and the other four are closed (example: $k$). \\

We shall consider two exact triangles of functors from $D^+(X)$ to $D^+(Z)$,
that are associated to the above setting: first, the triangle
\[
\j_! (k \j)^* \longto k^* \longto \i_* (k \i)^* \longto \j_! (k \j)^*[1]
\]
equal to $T_1 k^*$, where $T_1$ is the localization triangle
\[
\j_! \j^* \longto \id_Z \longto \i_* \i^* \longto \j_!  \j^*[1] \; ;
\]
second, the triangle
\[
Rk^! \longto k^* \longto k^* R j_* j^* \longto Rk^![1] 
\]
equal to $k^* T_2$, where $T_2$ is the co-localization triangle
\[
k_*Rk^! \longto  \id_X \longto R j_* j^* \longto k_*Rk^![1]
\]
associated to the complementary inclusions $j$ and $k$
(see \emph{e.g.} \cite[1.4.2.1]{BBD}).

\begin{Cons} \label{1Con}
We shall set up a canonical natural transformation of exact triangles 
\[
\vcenter{\xymatrix@R-10pt{
\j_! (k \j)^* \ar[r] \ar[d]_-{\ad} & k^* \ar[r] \ar@{=}[d] & 
  \i_* (k \i)^* \ar[r] \ar[d]^-{\i_* \i^* k^*(\ad)} & \j_! (k \j)^*[1] \ar[d]^-{\ad [1]}  \\
Rk^! \ar[r] & k^* \ar[r] & k^* R j_* j^* \ar[r] & Rk^![1] 
\\}}
\quad .
\] 
The natural transformation $\j_! (k \j)^* \to Rk^! $ is obtained 
by adjunction from the identity
$(k \j)^* \to R (k \j)^! = R\j^! Rk^! $ (the composition $k \j$ 
remains open).

As far as the natural transformation $\i_* (k \i)^* \to k^* R j_* j^*$ 
is concerned, note that (again since $k \j$ is open) the pull-back \emph{via}
$\j$ of 
\[
Rk^! \longto  k^* \longto k^* R j_* j^* \longto Rk^![1]
\]
equals
\[
(k \j)^* \stackrel{=}{\longto} (k \j)^* \longto \j^* (k^* R j_* j^*) \longto (k \j)^*[1] \; .
\]
In other words, $\j^* (k^* R j_* j^*) = 0$. Equivalently, the adjunction
\[
\ad: k^* R j_* j^* \longto \i_* \i^* (k^* R j_* j^*)
\]
is an isomorphism of functors from $D^+(X)$ to $D^+(Z)$. 
Define the natural transformation 
\[
\i_* \i^* k^* = \i_* (k \i)^* \to k^* R j_* j^* = \i_* \i^* (k^* R j_* j^*)
\]
as $\i_* \i^* k^*(\ad)$, where $\ad: \id_X \to R j_* j^*$ is adjoint to the identity 
$j^* \to j^*$. 

We leave it to the reader to show that the above \emph{does} yield
a natural transformation of exact triangles.
\end{Cons}

\begin{Prop} \label{1D}
(a)~Construction~\ref{1Con} is functorial
in the following sen\-se: let $f: X' \to X$ be a continuous map.
Denote the pull-back \emph{via} $f$ of
\[
\vcenter{\xymatrix@R-10pt{
        U \ar@{^{ (}->}[r]^-\circ \ar@{=}[d] &
        X - Z^0 \ar@{_{ (}->}[d] &
        Y \ar@{_{ (}->}[l] \ar@{_{ (}->}[d]^{\i} 
        \\
        U \ar@{^{ (}->}[r]_-{j}^-\circ &
        X &
        Z \ar@{_{ (}->}[l]^-{k} 
        \\   
        &
        Z^0 \ar@{=}[r] \ar@{^{ (}->}[u]^-\circ
        &
        Z^0 \ar@{^{ (}->}[u]^-\circ_-{\j}
\\}}
\] 
by 
\[
\vcenter{\xymatrix@R-10pt{
        U' \ar@{^{ (}->}[r]^-\circ \ar@{=}[d] &
        X' - {Z^0}' \ar@{_{ (}->}[d] &
        Y' \ar@{_{ (}->}[l] \ar@{_{ (}->}[d]^{\i'} 
        \\
        U' \ar@{^{ (}->}[r]_-{j'}^-\circ &
        X' &
        Z' \ar@{_{ (}->}[l]^-{k'} 
        \\   
        &
        {Z^0}' \ar@{=}[r] \ar@{^{ (}->}[u]^-\circ
        &
        {Z^0}' \ar@{^{ (}->}[u]^-\circ_-{\j'}
\\}} \quad ,
\] 
and consider the natural transformations of exact triangles
\[
\vcenter{\xymatrix@R-10pt{
\j_! (k \j)^* \ar[r] \ar[d]_-{\ad} & k^* \ar[r] \ar@{=}[d] & 
  \i_* (k \i)^* \ar[r] \ar[d]^-{\i_* \i^* k^*(\ad)} & \j_! (k \j)^*[1] \ar[d]^-{\ad [1]}  \\
Rk^! \ar[r] & k^* \ar[r] & k^* R j_* j^* \ar[r] & Rk^![1] 
\\}}
\quad\quad\quad\quad (\ref{1Con})_X
\] 
and
\[
\vcenter{\xymatrix@R-10pt{
\j_!' (k' \j')^* \ar[r] \ar[d]_-{\ad} & (k')^* \ar[r] \ar@{=}[d] & 
  \i_*' (k' \i')^* \ar[r] \ar[d]^-{\i_*' (\i')^* (k')^*(\ad)} & 
  \j'_! (k' \j')^*[1] \ar[d]^-{\ad [1]}  \\
R(k')^! \ar[r] & (k')^* \ar[r] & (k')^* R j_*' (j')^* \ar[r] & R(k')^![1] 
\\}}
\quad\quad\quad\quad (\ref{1Con})_{X'}
\] 
associated to the two settings by Construction~\ref{1Con}. Denote by
$f^Z$ the base change of $f$ to $Z$. Then the natural transformations
$\j_! (k \j)^* \to R f^Z_* \j_!' (k' \j')^* f^*$,
$k^* \to R f^Z_* (k')^* f^*$, $\i_* (k \i)^* \to R f^Z_* \i_*' (k' \i')^* f^*$,
$Rk^! \to R f^Z_* R(k')^! f^*$, $k^* \to R f^Z_* (k')^* f^*$ (again),
and $k^* R j_* j^* \to R f^Z_* (k')^* R j_*' (j')^* f^*$
organize into a na\-tural transformation of exact triangles 
$(\ref{1Con})_X \longto R f^Z_* (\ref{1Con})_{X'} f^*$
of functors $D^+(X) \to D^+(Z)$. 
In particular, the diagram
\[
\vcenter{\xymatrix@R-10pt{
\i_* (k \i)^*  \ar[r] \ar[d]_-{\i_* \i^* k^*(\ad)} & 
  R f^Z_* \i_*' (f k' \i')^* \ar[d]^-{R f^Z_* \i_*' (\i')^* (k')^*(\ad)f^*}  \\
k^* R j_* j^* \ar[r] & R f^Z_* (k')^* R j_*' (f j')^*
\\}}
\]
commutes. \\[0.1cm]
(b)~If the open immersion $U \into X - Z^0$ is contractible, then the transformation 
\[
\vcenter{\xymatrix@R-10pt{
\j_! (k \j)^* \ar[r] \ar[d]_-{\ad} & k^* \ar[r] \ar@{=}[d] & 
  \i_* (k \i)^* \ar[r] \ar[d]^-{\i_* \i^* k^*(\ad)} & \j_! (k \j)^*[1] \ar[d]^-{\ad [1]}  \\
Rk^! \ar[r] & k^* \ar[r] & k^* R j_* j^* \ar[r] & Rk^![1] 
\\}}
\quad\quad\quad\quad (\ref{1Con})_X
\] 
from 
Construction~\ref{1Con} 
is an isomorphism on the full
sub-triangulated category of $D^+(X)$ 
of complexes, whose cohomology objects are sheaves with locally constant
restriction to $X - Z^0$. 
In particular, for any locally constant sheaf
of Abelian groups $\CF$ on $X$, 
the object $Rk^! \CF$ of $D^+(Z)$ is concentrated in degree zero, and
\[
\ad:  \j_! (k \j)^* \CF \longto Rk^! \CF
\]
is an isomorphism.            
\end{Prop}

The transformations $\j_! (k \j)^* \to R f^Z_* \j_!' (k' \j')^* f^*$,
$k^* \to R f^Z_* (k')^* f^*$, 
$\i_* (k \i)^* \to R f^Z_* \i_*' (k' \i')^* f^*$,
$Rk^! \to R f^Z_* R(k')^! f^*$,
and $k^* R j_* j^* \to R f^Z_* (k')^* R j_*' (j')^* f^*$ from
Proposition~\ref{1D}~(a) are adjoint (under $((f^Z)^*,R f^Z_*)$)
to the identities (in the first three cases),
to $R(k')^! f^* (\ad_k)$,
where $\ad_k$ is the adjunction $k_!Rk^! \to \id_X$
(in the forth case),
and to $(f \circ k')^* (\ad_{f^U})$,
where $\ad_{f^U}$ is the adjunction $R j_* j^* \to R (f \circ j')_* (f \circ j')^*$
(in the last case).

\begin{Rem}
Part~(b) of Proposition~\ref{1D} specializes to Corollary~\ref{1C}~(b)
(for open immersions) if $Z^0$ is empty.
\end{Rem}

\medskip

\begin{Proofof}{Proposition~\ref{1D}}
(a): left to the reader.

\noindent (b): given Construction~\ref{1Con}, 
we need to show that under the additional hypo\-thesis on 
the immersion $j^0: U \into X - Z^0$,
the map 
\[
\i_* \i^* k^*(\ad: \CF \to R j_* j^* \CF)
\]
is an isomorphism in $D^+(Z)$,
whenever $\CF$ is a sheaf
of Abelian groups on $X$, whose restriction to $X - Z^0$ is locally constant.

But the functor $\i_* \i^* k^*$ factors through $i^*$, where $i$ denotes the
(closed) immersion of $X - Z^0$ into $X$, and
\[
i^*(\ad): i^* \CF \to i^*R j_* j^* \CF = R j^0_* (j^0)^* i^* \CF
\]
is adjoint to the identity $(j^0)^* i^* \CF \to (j^0)^* i^* \CF$. The sheaf
$i^* \CF$ being locally constant, Corollary~\ref{1C}~(b) tells us that 
$i^*(\ad)$ is an isomorphism. Hence so is $\i_* \i^* k^*(\ad)$.
\end{Proofof}


\bigskip

%
%

\section{The combinatorial set-up}
\label{2}



Let us fix \emph{mixed Shimura data}  $(P,\FX)$ \cite[Def.~2.1]{P}.
Recall the notion of \emph{admissible} parabolic sub-group of $P$ \cite[Def.~4.5]{P}.
Any maximal proper parabolic of $P$ is admissible. If the maximal semi-simple quotient of $P$
is simple, then any admissible parabolic is maximal proper, or equal to $P$. 
To any admissible parabolic sub-group $Q_j$ of $P$ is associated
a canonical normal sub-group $P_j \subset Q_j$ \cite[Sect.~4.7]{P}, underlying 
finitely many \emph{rational boundary components} $(P_j,\FX_j)$ of $(P,\FX)$ 
\cite[Sect.~4.11]{P}.

\begin{Def}[{\cite[Def.~4.1]{W1}}] \label{2A}
Let $Q_1$ and $Q_2$ be two admissible parabo\-lic sub-groups
of $P$, with associated canonical normal sub-groups $P_1 \subset Q_1$ and $P_2 \subset Q_2$.
We define the relation
\[
Q_1 \preceq Q_2
\]
to hold if boundary components $(P_1,\FX_1)$ and $(P_2,\FX_2)$ of $(P,\FX)$
can be chosen such that $(P_1,\FX_1)$ is a boundary component of $(P_2,\FX_2)$. 
\end{Def}

Equivalently \cite[Sect.~4.11]{P}, \emph{any} boundary component $(P_2,\FX_2)$ admits a boun\-dary
component of the the form $(P_1,\FX_1)$. This implies that $\preceq$ is transitive. 
Note also that $Q_1 \preceq Q_2$ only if $P_1 \subset P_2$. \\

It follows from \cite[Prop.~4.2]{W1} that $\preceq$ is anti-symmetric. Therefore,
$\preceq$ is a partial order on the set of admissible parabolic sub-groups of $P$.
Let us write $Q_1 \prec Q_2$ if $Q_1 \preceq Q_2$ and $Q_1 \ne Q_2$. 

\begin{Prop} \label{2B}
Let $Q_1 \preceq Q_2 \preceq \ldots \preceq Q_r$ a chain of admissible parabo\-lic sub-groups of $P$.
\\[0.1cm]
(a)~The intersection $\cap_{j=1}^r Q_j$ remains parabolic. \\[0.1cm]
(b)~We have the equality $Q_r = (\cap_{j=1}^r Q_j) P_r$.
\end{Prop}

\begin{Proof}
We apply induction on $r$, both (a) and (b) being trivial for $r=1$.

Following \cite[Lemma~4.19~(b), (c)]{P}, the intersection $Q_j \cap P_r$ is admissible
parabolic in $P_r$, for all $j = 1,\ldots,r-1$, and 
we get an induced chain 
\[
Q_1 \cap P_r \preceq Q_2 \cap P_r \preceq \ldots \preceq Q_{r-1} \cap P_r
\]
of parabolics of $P_r$.
According to the induction hypo\-thesis for claim~(a) (applied to any of the boundary
components $(P_r,\FX_r)$ of $(P,\FX)$),
the intersection $\cap_{j=1}^{r-1} (Q_j \cap P_r)$
is parabolic in $P_r$, while the induction hypothesis for claim~(b) implies that 
$Q_{r-1} = (\cap_{j=1}^{r-1} Q_j) P_{r-1}$.

From \cite[Lemma~4.19~(b)]{P}, we get $Q_r = (Q_{r-1} \cap Q_r)P_r$, hence by the preceding equality
\[
Q_r = Q_{r-1}P_r \cap Q_r = (\cap_{j=1}^{r-1} Q_j) P_r \cap Q_r = (\cap_{j=1}^r Q_j) P_r \; .
\]
It follows that the inclusion of the normal sub-group $P_r$ into $Q_r$ induces
an isomorphism of quotient varieties
\[
P_r / \bigl( \cap_{j=1}^{r-1} (Q_j \cap P_r) \bigr) 
= P_r / \bigl( \bigl( \cap_{j=1}^r Q_j \bigr) \cap P_r \bigr) 
\isoto Q_r / \bigl( \cap_{j=1}^r Q_j \bigr) \; .
\]
\end{Proof}

Recall the map $\adm$, defined in \cite[Def.~4.6]{W1}, from the set of all parabolic sub-groups of $P$
to the one of admissible parabolics: for a parabolic $Q$, we have $\adm(Q) = Q_1$ if and only if
$P_1 \subset Q \subset Q_1$ ($P_1=$ the canonical normal sub-group of $Q_1$ defined in 
\cite[Sect.~4.7]{P}). 

\begin{Prop} \label{2C}
Let $Q_1 \preceq Q_2 \preceq \ldots \preceq Q_r$ a chain of admissible parabo\-lic sub-groups of $P$.
Then 
\[
\adm \bigl( \cap_{j=1}^r Q_j \bigr) = Q_1 \; .
\]
\end{Prop}

\begin{Proof}
We have $P_1 \subset P_2 \subset \ldots \subset P_r$, which implies
\[
P_1 \subset \cap_{j=1}^r Q_j \subset Q_1 \; .
\]
Hence $Q_1$ has the defining property of $\adm ( \cap_{j=1}^r Q_j )$.
\end{Proof}

\begin{Def} \label{2D}
(a)~Define a partially ordered set $(\CC_{(P,\FX)},\subset)$ as follows:
\[
\CC_{(P,\FX)} := \{ (Q_1 \prec Q_2 \prec \ldots \prec Q_r) \tei r \ge 1 \; , \; Q_j \prec P
\; \text{admissible parabolic} \; \forall \, j \} \; ; 
\]
\[
(Q_1 \prec Q_2 \prec \ldots \prec Q_r) 
\subset (\tilde{Q}_1 \prec \tilde{Q}_2 \prec \ldots \prec \tilde{Q}_s)
\]
meaning
that $\{ \tilde{Q}_i \tei i = 1,\ldots,s \}$ is contained in 
$\{ Q_j \tei j = 1,\ldots,r \}$ (in other words, each $\tilde{Q}_i$ equals one of
the $Q_j$). \\[0.1cm]
(b)~Define two maps $b$ and $\bigcap$ from $\CC_{(P,\FX)}$ to the set of proper
parabolics of $P$ as follows:
\[
b(Q_1 \prec Q_2 \prec \ldots \prec Q_r) := Q_1 \; ,
\]
\[
\bigcap (Q_1 \prec Q_2 \prec \ldots \prec Q_r) := \cap_{j=1}^r Q_j \; .
\] 
(c)~Let $Q_1$ be a proper admissible parabolic sub-group of $P$. Define
\[
\CC_{(P,\FX) \tei Q_1} := b^{-1}(Q_1) \subset \CC_{(P,\FX)} \; .
\]
\end{Def}

Note that the map $\bigcap$ is order preserving.
According to Proposition~\ref{2C}, we have
\[
\adm \circ \bigcap = b \; .
\]
In particular, $\bigcap$ induces a map, denoted
by the same symbol
\[
\bigcap: \CC_{(P,\FX) \tei Q_1} \longto \adm^{-1}(Q_1) \; ,
\]
for every proper admissible parabolic $Q_1$ of $P$. 

\begin{Lem} \label{2E}
Let $Q$ be a parabolic sub-group of $P$. Put $Q_1 := \adm(Q)$, and consider the set 
\[
\CS := \{ Q_2 
\tei Q_2 \; \text{admissible parabolic} \; , \; Q_1 \preceq Q_2 \; , \; Q \subset Q_2 \} \; .
\]
If $Q_1$ is maximal proper, then $\CS$ is totally ordered
(under $\preceq$), and every element of $\CS$ is either maximal proper, or equal to $P$.  
\end{Lem}

\begin{Proof}
Using the formalism of \emph{quotient Shimura data} \cite[Prop.~2.9]{P}, 
and \emph{functoriality} of boundary components \cite[Sect.~4.16]{P}, we may, by dividing
out first, the unipotent radical, and then, the center, assume that $P$ is equal
to a finite product of simple groups $G_i \,$.

Then $Q_1$, being maximal proper, equals the product of parabolics of $G_i \,$,
all of which except one, say associated to the index $i_0 \,$, are equal to $G_i \,$. 
Therefore, its sub-group $P_1$ is the product of sub-groups of $G_i \,$,
all of which except the one associated to the index $i_0 \,$, are normal in $G_i \,$.
It follows that any parabolic sub-group of $P$ containing $P_1$ is again
the product of parabolics of $G_i \,$,
all of which except the one associated to the index $i_0 \,$, are equal to $G_i \,$. 
By dividing out $\prod_{i \ne i_0} G_i$, we are thus reduced to the case where $P$ is simple.

Under this assumption, 
the restriction of the relation $\preceq$ to the (finitely many)
admissible parabolic sub-groups containing $Q$ is a total order \cite[Rem.~(ii) on p.~91]{P}.

For simple $P$, by definition, a parabolic of $P$ is admissible if and only if
it is either maximal proper, or equal to $P$
\cite[Def.~4.5]{P}.
\end{Proof}

\begin{Rem}
(a)~Lemma~\ref{2E} is false if $Q_1$ is (admissible, but) not maximal proper: in this case,
$Q_1$ equals the intersection of two admissible parabolic sub-groups $Q_2$ and $Q_2'$
(hence $Q_1 \preceq Q_2$ and $Q_1\preceq Q_2'$), with
$Q_2 \not \preceq Q_2'$ and $Q_2' \not \preceq Q_2 \,$. \\[0.1cm]
(b)~In the situation of (a), observe that we also have
\[
\bigcap (Q_1 \prec Q_2) = Q_1 = \bigcap (Q_1 \prec Q_2') 
\]
(cmp.~part~(a) of the Theorem~\ref{2F} below).
\end{Rem}

\begin{Thm} \label{2F}
(a)~Let $Q_1$ be a maximal proper parabolic sub-group of $P$. Then
\[
\bigcap: \CC_{(P,\FX) \tei Q_1} \longto \adm^{-1}(Q_1) 
\]
is an ordered bijection. \\[0.1cm]
(b)~If the maximal semi-simple quotient of $P$ is simple, then
\[
\bigcap: \CC_{(P,\FX)} \longto \{ \text{proper parabolics of} \; P \}
\]
is an ordered bijection.
\end{Thm}

A bijection $\alpha$ 
between ordered sets $(A,\subset)$ and $(B,\le)$ is called \emph{ordered}
if for any two elements $a_1, a_2$ of $A$, we have $\alpha(a_1) \le \alpha(a_2)$
if and only if $a_1 \subset a_2$.

\medskip 

\begin{Proofof}{Theorem~\ref{2F}}
(a):~let $Q \in \adm^{-1}(Q_1)$. In particular, we have $Q \subset Q_1$. Let
$Q_2,\ldots,Q_r$ be the other, pairwise unequal, maximal proper parabolics containing $Q$.
We thus have $Q = \cap_{j=1}^r Q_j$. According to \cite[Thm.~4.5~(b)]{W1}, 
\[
Q_1 \prec Q_j \; , \; \forall \, j = 2,\ldots,r \; .
\]
Apply Lemma~\ref{2E}: the $Q_j$, $j  = 1,\ldots,r$ are totally ordered.
Without loss of generality, we may assume
\[
Q_1 \prec Q_2 \prec \ldots \prec Q_r \; .
\]
We obtain that $Q = \bigcap (Q_1 \prec Q_2 \prec \ldots \prec Q_r)$; this shows surjectivity
of $\bigcap$ \,.

As for injectivity, let 
\[
Q_1 \prec Q_2 \prec \ldots \prec Q_r \quad \text{and} \quad 
Q_1 \prec \tilde{Q}_2 \prec \ldots \prec \tilde{Q}_s
\]
be two chains in $\CC_{(P,\FX) \tei Q_1}$ satisfying
\[
\cap_{j=1}^r Q_j = Q_1 \cap \cap_{i=2}^s \tilde{Q}_j \; .
\]
Put $Q := \cap_{j=1}^r Q_j$.
According to Lemma~\ref{2E}, all $Q_j$ and $\tilde{Q}_i$
are maximal proper. Therefore, they are \emph{the} maximal proper parabolics containing $Q$.
Thus, the sets $\{ Q_j \tei j \ge 2 \}$ and $\{ \tilde{Q}_i \tei i \ge 2 \}$ are equal.

Similarly, if  
\[
Q_1 \prec Q_2 \prec \ldots \prec Q_r \quad \text{and} \quad 
Q_1 \prec \tilde{Q}_2 \prec \ldots \prec \tilde{Q}_s
\]
in $\CC_{(P,\FX) \tei Q_1}$ satisfy
\[
\cap_{j=1}^r Q_j \subset Q_1 \cap \cap_{i=2}^s \tilde{Q}_j \; ,
\]
then $\{ \tilde{Q}_i \tei i \ge 2 \}$ is contained in $\{ Q_j \tei j \ge 2 \}$. Therefore,
the bijection $\bigcap$ is ordered.

\noindent (b)~: under the assumption on $P$, a parabolic of $P$ is admissible if and only if
it is either maximal proper, or equal to $P$ \cite[Def.~4.5]{P}.
The only parabolic $Q$ whose image under $\adm$ is not maximal proper, is $Q = P$.
Now apply (a).
\end{Proofof}
 

\bigskip
\newpage

%
%

\section{The sheaf-theoretical set-up}
\label{3}



From now on, assume that our mixed Shimura data  $(P,\FX)$ satisfy hypotheses 
$(+)$ and $(U=0)$ from \cite{W1}: denote by $G$ the maximal reductive quotient of $P$.
\begin{enumerate}
\item [$(+)$] The neutral connected component $Z (G)^0$ of the center $Z (G)$ of 
$G$ is, up to isogeny, a direct product of a $\BQ$-split torus with a torus 
$T$ of compact type (\emph{i.e.}, $T(\BR)$ is compact) defined over $\BQ \,$.
\end{enumerate}
\begin{enumerate}
\item [$(U=0)$] The weight $(-2)$-part of $P$ \cite[Def.~2.1~(v)]{P} is trivial. 
\end{enumerate}

Note that hypothesis $(U=0)$ is satisfied if $(P,\FX)$ is \emph{pure}, \emph{i.e.}, if
$P=G$ is reductive. \\

Imitating the construction from \cite{BS}, the manifold with corners $\FX^{BS}$
was defined in \cite[Def.~3.3]{W1},
together with its stratification 
\[
\FX^{BS} = \coprod_R e(R) 
\]
by locally closed strata, indexed by the parabolic sub-groups $R$ of $P$.
The unique open stratum $e(P)$ equals $\FX$. This was then used \cite[Def.~3.13]{W1}
to define the \emph{Borel--Serre
compactification of $M^K (P,\FX)$} as 
\[
M^K (P,\FX) (\BC)^{BS} := P (\BQ) \backslash \bigl( \FX^{BS} \times P (\BA_f) / K \bigr) \; ,
\]
for any open compact sub-group $K$ of $P(\BA_f)$. 

\begin{Def}[{\cite[Def.~5.1]{W1}}] \label{3A}
Let $K$ be an open compact sub-group of $P(\BA_f)$. \\[0.1cm] 
(a)~Let $Q$ be a parabolic sub-group of $P$. Define 
\[
e^K \bigl( Q,P(\BA_f) \bigr) \subset M^K  (P,\FX) (\BC)^{BS}
\]
to be the image of $e(Q) \times P (\BA_f) / K \subset \FX^{BS} \times P (\BA_f) / K$
under the projection from $\FX^{BS} \times P (\BA_f) / K$. \\[0.1cm]
(b)~The \emph{canonical stratification of
the Borel--Serre compactification}
is the stratification by the $e^K ( Q,P(\BA_f) )$, where $Q$ runs through the
$P(\BQ)$-conjugate classes of parabolic sub-groups of $P$.
\end{Def}

Denoting by $\overline{e^K ( Q,P(\BA_f) )}$
the closure of $e^K(Q,P(\BA_f))$, we have
\[
\overline{e^K \bigl( Q,P(\BA_f) \bigr)} = \bigcup_{R \subset Q} e^K \bigl( R,P(\BA_f) \bigr) \; ,
\]
where $R$ runs through all parabolics contained in $Q$.

\begin{Def}[{\cite[Def.~5.3]{W1}}] \label{3B}
Let $K$ be an open compact sub-group of $P(\BA_f)$, and
$Q$ be a parabolic sub-group of $P$. Define 
\[
e^K \bigl( Q,P(\BA_f) \bigr)' \subset \overline{e^K \bigl( Q,P(\BA_f) \bigr)}
\]
as the union of the $e^K ( R,P(\BA_f) )$, where $R$ runs through the
parabolics contained in $Q$, and satisfying $\adm(R) = \adm(Q)$.
\end{Def}

In other words, 
\[
e^K \bigl( Q,P(\BA_f) \bigr)' 
             = \bigcup_{\adm_{Sh}(Q) \subset R \subset Q} e^K \bigl( R,P(\BA_f) \bigr) \; ,
\]
where $\adm_{Sh}(Q)$ is the canonical normal sub-group of the admissible parabolic 
$\adm(Q)$ \cite[Def.~4.6]{W1}. The set $e^K ( Q,P(\BA_f))'$ contains $e^K ( Q,P(\BA_f))$,
and is therefore dense in $\overline{e^K ( Q,P(\BA_f))}$.
According to \cite[Prop.~5.4~(b)]{W1}, it is open in $\overline{e^K ( Q,P(\BA_f) )}$.
Denote by $k_Q$ the immersion of $e^K ( Q,P(\BA_f))'$ into $M^K  (P,\FX) (\BC)^{BS}$.

\begin{Def} \label{3C}
Let $K$ be an open compact sub-group of $P(\BA_f)$.
Define a category $\CC_{(P,\FX)}^{K,BS}$ as follows: objects are of the form
\[
\bigl( \CK_{\uQ} , a_{\uQ \subset \utQ} \bigr)_{\uQ, \uQ \subset \utQ} \; ,
\]
where $\CK_{\uQ} \in D^+(e^K ( \bigcap \uQ,P(\BA_f) )')$,
for each $\uQ \in \CC_{(P,\FX)}$, and 
\[
a_{\uQ \subset \utQ}: R (k_{\bigcap \utQ})_* \CK_{\utQ} \longto R (k_{\bigcap \uQ})_* \CK_{\uQ}
\]
is a morphism in $D^+(M^K  (P,\FX) (\BC)^{BS})$, for each relation $\uQ \subset \utQ$
in $\CC_{(P,\FX)}$.
Morphisms between $(\CK_{\uQ} , a_{\uQ \subset \utQ})_{\uQ, \uQ \subset \utQ}$
and $(\CL_{\uQ} , b_{\uQ \subset \utQ})_{\uQ, \uQ \subset \utQ}$ are of the form
$( \alpha_{\uQ} )_{\uQ}$,
where $\alpha_{\uQ}: \CK_{\uQ} \to \CL_{\uQ}$ is a morphism in 
$D^+(e^K ( \bigcap \uQ,P(\BA_f) )')$,
for each $\uQ \in \CC_{(P,\FX)}$, such that for each relation $\uQ \subset \utQ$
in $\CC_{(P,\FX)}$, the diagram
\[
\vcenter{\xymatrix@R-10pt{
R (k_{\bigcap \utQ})_* \CK_{\utQ} \ar[r]^-{a_{\uQ \subset \utQ}} 
                                   \ar[d]_-{R (k_{\bigcap \utQ})_* \alpha_{\utQ}} 
  & R (k_{\bigcap \uQ})_* \CK_{\uQ} \ar[d]^-{R (k_{\bigcap \uQ})_* \alpha_{\uQ}}  \\
R (k_{\bigcap \utQ})_* \CL_{\utQ} \ar[r]^-{b_{\uQ \subset \utQ}} 
  & R (k_{\bigcap \uQ})_* \CL_{\uQ}
\\}}
\] 
commutes.
\end{Def}

The reader might wish to include in Definition~\ref{3C} the cocyle condition
$a_{\uQ \subset \utQ} \circ a_{\utQ \subset \underline{\tilde{\tilde{Q}}}} 
= a_{\uQ \subset \underline{\tilde{\tilde{Q}}}}$
for double relations $\uQ \subset \utQ \subset \underline{\tilde{\tilde{Q}}}$,
and the condition $\CK_{\uQ} = \CK_{\uQ'}$ whenever $\uQ$ and $\uQ'$ are
conjugate under $P(\BQ)$. However,
in the cases of interest for us, these conditions will be automatically satisfied. \\

Let us now discuss a variant of the above construction for the 
\emph{Baily--Borel compactification} $M^K (P,\FX)^*$ of $M^K (P,\FX)$.
It is defined if the Shimura data $(P,\FX) = (G,\FX)$ are pure. The space
of complex points equals
\[
M^K (G,\FX)^* (\BC) = G (\BQ) \backslash \bigl( \FX^* \times G (\BA_f) / K \bigr) \; ,
\]
where 
\[
\FX^* := \coprod \FX_j/W_j 
\]
is the disjoint union over all rational boundary components $(P_j,\FX_j)$
of $(G,\FX)$ ($W_j :=$ the unipotent radical of $P_j$)
(\cite[Chap.~III, Sect.~6.1]{AMRT} or \cite[Sect.~6.2]{P}). 
The inclusion $\FX \into \FX^*$ is open, and so is
\[
 M^K (G,\FX) \longinto M^K (G,\FX)^* \; .
 \]

\begin{Def}[{\cite[Def.~5.8]{W1}}] \label{3D}
Assume that the Shimura data $(P,\FX) = (G,\FX)$ are pure. 
Let $K$ be an open compact sub-group of $G(\BA_f)$. \\[0.1cm]
(a)~Let $Q_j$ be an admissible parabolic sub-group of $G$. Consider the disjoint union 
$\coprod \FX_j$ of the finitely many spaces $\FX_j$ underlying
rational boundary components associated to $Q_j \,$, as well as its quotient
$\coprod \FX_j / W_j \,$. Define
\[
M^K (Q_j,\FX) \subset M^K (G,\FX)^*
\]
to be the image of $\coprod \FX_j / W_j \times G(\BA_f) / K$ under the projection
\[
\FX^* \times G (\BA_f) / K \longonto M^K (G,\FX)^* (\BC) \; . 
\]
(b)~The \emph{canonical stratification of the Baily--Borel compactification}
is the stratification by the $M^K (Q_j,\FX)$, where $Q_j$ runs through the
$G(\BQ)$-conjugate classes of admissible para\-bolic sub-groups of $G$.
\end{Def}

One of the main rationality results from \cite[Sect.~12]{P} implies that  
the image of $\coprod \FX_j / W_j \times G(\BA_f) / K$ under 
\[
\FX^* \times G (\BA_f) / K \longonto M^K (G,\FX)^* (\BC) 
\]
is indeed identified with the set of $\BC$-valued points of a locally
closed sub-scheme $M^K (Q_j,\FX)$ of $M^K (G,\FX)^*$ \cite[Rem.~5.9~(b)]{W1}. 
Denote by $i_{Q_j}$ the immersion of $M^K (Q_j,\FX)$ into $M^K (G,\FX)^*$.

\begin{Def} \label{3E}
Assume that the Shimura data $(P,\FX) = (G,\FX)$ are pure. 
Let $K$ be an open compact sub-group of $G(\BA_f)$. \\[0.1cm]
(a)~Define a category $\CC_{(G,\FX)}^{K,*}$ as follows: objects are of the form
\[
\bigl( \CK_{\uQ} , a_{\uQ \subset \utQ} \bigr)_{\uQ, \uQ \subset \utQ} \; ,
\]
where $\CK_{\uQ} \in D^+(M^K (b(\uQ),\FX)(\BC))$,
for each $\uQ \in \CC_{(G,\FX)}$, and 
\[
a_{\uQ \subset \utQ}: R (i_{b(\utQ)})_* \CK_{\utQ} \longto R (i_{b(\uQ)})_* \CK_{\uQ}
\]
is a morphism in $D^+(M^K  (G,\FX)^* (\BC))$, for each relation $\uQ \subset \utQ$
in $\CC_{(G,\FX)}$. 
Morphisms between $(\CK_{\uQ} , a_{\uQ \subset \utQ})_{\uQ, \uQ \subset \utQ}$
and $(\CL_{\uQ} , b_{\uQ \subset \utQ})_{\uQ, \uQ \subset \utQ}$ are of the form
$( \alpha_{\uQ} )_{\uQ}$,
where $\alpha_{\uQ}: \CK_{\uQ} \to \CL_{\uQ}$ is a morphism in 
$D^+(M^K (b(\uQ),\FX)(\BC))$,
for each $\uQ \in \CC_{(G,\FX)}$, such that for each relation $\uQ \subset \utQ$
in $\CC_{(G,\FX)}$, the diagram
\[
\vcenter{\xymatrix@R-10pt{
R (i_{b(\utQ)})_* \CK_{\utQ} \ar[r]^-{a_{\uQ \subset \utQ}} 
                                   \ar[d]_-{R (i_{b(\utQ)})_* \alpha_{\utQ}} 
  & R (i_{b(\uQ)})_* \CK_{\uQ} \ar[d]^-{R (i_{b(\uQ)})_* \alpha_{\uQ}}  \\
R (i_{b(\utQ)})_* \CL_{\utQ} \ar[r]^-{b_{\uQ \subset \utQ}}  & R (i_{b(\uQ)})_* \CL_{\uQ}
\\}}
\] 
commutes. \\[0.1cm]
(b)~Let $Q_j$ be a proper admissible parabolic sub-group of $G$.  
Define a category $\CC_{(Q_j,\FX)}^K$ as follows: objects are of the form
\[
\bigl( \CK_{\uQ} , a_{\uQ \subset \utQ} \bigr)_{\uQ, \utQ \in \CC_{(G,\FX) \tei Q_j}} \; ,
\]
where $\CK_{\uQ} \in D^+(M^K (Q_j,\FX)(\BC))$,
for each $\uQ \in \CC_{(G,\FX) \tei Q_j}$, and 
\[
a_{\uQ \subset \utQ}: R (i_{Q_j})_* \CK_{\utQ} \longto R (i_{Q_j})_* \CK_{\uQ}
\]
is a morphism in $D^+(M^K  (G,\FX)^* (\BC))$, for each relation $\uQ \subset \utQ$
in $\CC_{(G,\FX) \tei Q_j}$.
Morphisms are defined as in $\CC_{(G,\FX)}^{K,*}$. \\[0.1cm]
(c)~Let $Q_j$ be a proper admissible parabolic sub-group of $G$.
Define
\[
i^*_{Q_j}: \CC_{(G,\FX)}^{K,*} \longto \CC_{(Q_j,\FX)}^K
\]
as the restriction 
\[
\bigl( \CK_{\uQ} , a_{\uQ \subset \utQ} \bigr)_{\uQ, \uQ \subset \utQ} \longmapsto 
\bigl( \CK_{\uQ} , a_{\uQ \subset \utQ} \bigr)_{\uQ, \utQ \in \CC_{(G,\FX) \tei Q_j}} \; .
\]
\end{Def}

If $(P,\FX) = (G,\FX)$ are pure, then the continuous map 
\[
p^K: M^K  (G,\FX) (\BC)^{BS} \to M^K (G,\FX)^* (\BC)
\]
from
\cite[Cor.~4.18]{W1} maps $e^K \bigl( Q_j, G(\BA_f) \bigr)'$ to
$M^K (Q_j,\FX)(\BC)$, for any admissible parabolic sub-group $Q_j$ of $G$ 
\cite[Comp.~4.15]{W1}.
In particular, for any parabolic $Q$ of $G$, we have
\[
p^K \bigl( e^K \bigl( Q, G(\BA_f) \bigr)' \bigr) \subset M^K (\adm(Q),\FX)(\BC) 
          \subset M^K  (G,\FX)^* (\BC) \; .
\] 

\begin{Prop} \label{3F}
Assume that the Shimura data $(P,\FX) = (G,\FX)$ are pure. 
Let $K$ be an open compact sub-group of $G(\BA_f)$. The higher direct images
\[
Rp^K_*: D^+ \bigl( e^K \bigl( Q, G(\BA_f) \bigr)' \bigr) 
                                      \longto D^+ \bigl(M^K  (\adm(Q),\FX)^* (\BC) \bigr) \; ,
\]
for all parabolic sub-groups $Q$ of $G$,
induce a functor, denoted by the same symbol
\[
Rp^K_*: \CC_{(G,\FX)}^{K,BS} \longto \CC_{(G,\FX)}^{K,*} \; ,
\]
On objects, it is given by
\[
\bigl( \CK_{\uQ} , a_{\uQ \subset \utQ} \bigr)_{\uQ, \uQ \subset \utQ} \longmapsto
\bigl( Rp^K_* \CK_{\uQ} , Rp^K_* a_{\uQ \subset \utQ} \bigr)_{\uQ, \uQ \subset \utQ} \; .
\]
On morphisms, it is given by $( \alpha_{\uQ} )_{\uQ} \mapsto ( Rp^K_* \alpha_{\uQ} )_{\uQ}$.
\end{Prop}

\begin{Proof}
This follows from the relations $\adm \circ \bigcap = b$ and 
\[
p^K \circ k_Q = i_{\adm(Q)}(\BC) \circ p^K: 
e^K \bigl( Q, G(\BA_f) \bigr)' \longto M^K  (G,\FX)^* (\BC) \; , 
\]
for any parabolic $Q$.
\end{Proof}


\bigskip

%
%

\section{First Main Theorem: multiple degeneration in the 
Borel--Serre compactification}
\label{4}



Fix mixed Shimura data  $(P,\FX)$ satisfying hypotheses 
$(+)$ and $(U=0)$, and an open compact sub-group $K$ of $P(\BA_f)$. 

\begin{Def} \label{4A}
Define the \emph{boundary of the Borel--Serre compactification} of $M^K  (P,\FX)$ to be
the closed complement 
\[
\partial M^K  (P,\FX) (\BC)^{BS} := M^K  (P,\FX) (\BC)^{BS} - M^K  (P,\FX) (\BC)
\] 
of $M^K  (P,\FX) (\BC)$ in $M^K  (P,\FX) (\BC)^{BS}$.
\end{Def}

\begin{Cons} \label{4Con1}
Let $\CV \in D^+(\partial M^K  (P,\FX) (\BC)^{BS})$. Define an object
$\deg(\CV)$ of $\CC_{(P,\FX)}^{K,BS}$ as follows.
\[
\deg(\CV) = \bigl( \CK_{\uQ} , a_{\uQ \subset \utQ} \bigr)_{\uQ, \uQ \subset \utQ} \; ,
\]
where for $\uQ = (Q_1 \prec Q_2 \prec \ldots \prec Q_r) \in \CC_{(P,\FX)}$, 
\[
\CK_{\uQ} :=  
k_{\bigcap \uQ}^* \bigl( R(k_{Q_1})_* k_{Q_1}^* R(k_{Q_2})_* k_{Q_2}^* \ldots R(k_{Q_r})_* k_{Q_r}^* \CV \bigr) \; ,
\]
and for $\uQ \subset \utQ$,
\[
a_{\uQ \subset \utQ}: R (k_{\bigcap \utQ})_* \CK_{\utQ} \longto R (k_{\bigcap \uQ})_* \CK_{\uQ}
\]
is the adjunction.
\end{Cons}
 
Note that all immersions $k_Q$ occurring in Construction~\ref{4Con1} are associated
to proper parabolic sub-groups $Q$ of $P$. Therefore, their images are contained in
$\partial M^K  (P,\FX) (\BC)^{BS}$ --- which justifies the following abuse of notation:
use the same symbol $k_Q$ to denote
the immersion of $e^K ( Q,P(\BA_f) )'$ into 
the boundary $\partial M^K  (P,\FX) (\BC)^{BS}$. \\

We leave it to the reader to complete Construction~\ref{4Con1} to a functor
\[
\deg: D^+ \bigl( \partial M^K  (P,\FX) (\BC)^{BS} \bigr) \longto \CC_{(P,\FX)}^{K,BS} \; .
\]
In order to analyse $\deg$ more closely, the following will be useful. 

\begin{Prop} \label{4Aa}
Let $\uQ = (Q_1 \prec \ldots \prec Q_r) \in \CC_{(P,\FX)}$. 
For any $s \in \{1, \ldots, r \}$, we have
\[
e^K \bigl( Q_s,P(\BA_f) \bigr)' \cap 
              \overline{e^K \bigl( Q_{s+1} \cap \ldots \cap Q_r,P(\BA_f) \bigr)} 
              = e^K \bigl( Q_s \cap \ldots \cap Q_r,P(\BA_f) \bigr)' \; .
\]
\end{Prop}

\begin{Proof}
This is a special case of \cite[Prop.~5.6]{W1}.
\end{Proof}

\begin{Cor} \label{4Ab}
Let $\CV \in D^+(\partial M^K  (P,\FX) (\BC)^{BS})$, and write
\[
\deg(\CV) = \bigl( \CK_{\uQ} , a_{\uQ \subset \utQ} \bigr)_{\uQ, \uQ \subset \utQ} \; .
\]
Let $\uQ = (Q_1 \prec Q_2 \prec \ldots \prec Q_r) \in \CC_{(P,\FX)}$. \\[0.1cm]
(a)~Let $s \in \{1, \ldots, r \}$. The object
\[
k_{Q_s}^* \ldots R(k_{Q_r})_* k_{Q_r}^* \CV 
\]
of $D^+(e^K ( Q_s,P(\BA_f) )')$ has support in 
$e^K ( Q_s \cap \ldots \cap Q_r,P(\BA_f) )'$, and
the object
\[
R(k_{Q_s})_* k_{Q_s}^* \ldots R(k_{Q_r})_* k_{Q_r}^* \CV 
\]
of $D^+(\overline{e^K ( Q_s,P(\BA_f) )})$
has support in $\overline{e^K ( Q_s \cap \ldots \cap Q_r,P(\BA_f) )}$. \\[0.1cm]
(b)~Writing $k_{\uQ,Q_1}$ for the immersion of $e^K ( \bigcap \uQ,P(\BA_f) )'$ into 
$e^K ( Q_1,P(\BA_f) )'$, the adjunction
\[
\ad : k_{Q_1}^* \ldots R(k_{Q_r})_* k_{Q_r}^* \CV \longto (k_{\uQ,Q_1})_* \CK_{\uQ}
\]
is an isomorphism. 
\end{Cor}

Note that $(k_{\uQ,Q_1})_* = R(k_{\uQ,Q_1})_*$ as $k_{\uQ,Q_1}$ is closed.

\medskip

\begin{Proofof}{Corollary~\ref{4Ab}}
In order to prove the first statement of (a) (which implies the second),
use descending induction on $s$, and Proposition~\ref{4Aa}.

As for (b), setting $\CW:= k_{Q_1}^* \ldots R(k_{Q_r})_* k_{Q_r}^* \CV$,
part~(a) means that
\[
\ad : \CW \longto (k_{\uQ,Q_1})_* k_{\uQ,Q_1}^* \CW
\]
is an isomorphism. But the equality of functors 
\[
k_{\uQ,Q_1}^* = k_{\bigcap \uQ}^* R(k_{Q_1})_* \; 
\]
shows that $k_{\uQ,Q_1}^* \CW = \CK_{\uQ} \,$.
\end{Proofof}

The components $\CK_{\uQ}$ of our second construction are easier to define
than the ones of $\deg(\CV)$. But they
organise into an object of $\CC_{(P,\FX)}^{K,BS}$ only under additional assumptions.

\begin{Prop}[{\cite[Prop.~5.4~(c)]{W1}}] \label{4B}
Let $Q$ be a parabolic sub-group of $P$. 
If $K$ is \emph{neat} in the sense of \cite[Sect.~0.6]{P}, then the open immersion
\[
e^K \bigl( Q,P(\BA_f) \bigr)' \longinto \overline{e^K \bigl( Q,P(\BA_f) \bigr)}  
\]
is contractible.
\end{Prop}
 
\begin{Cons} \label{4Con2}
Let $\CV \in D^+(\partial M^K  (P,\FX) (\BC)^{BS})$, and assume that its
cohomology objects are locally constant.
Assume also that $K$ is neat. Define an object
$\res(\CV)$ of $\CC_{(P,\FX)}^{K,BS}$ as follows.
\[
\res(\CV) = \bigl( \CK_{\uQ} , a_{\uQ \subset \utQ} \bigr)_{\uQ, \uQ \subset \utQ} \; ,
\]
where for $\uQ \in \CC_{(P,\FX)}$, 
\[
\CK_{\uQ} := k_{\bigcap \uQ}^* \CV = \CV_{\tei e^K ( \bigcap \uQ,P(\BA_f) )'}\; , 
\]
and for $\uQ \subset \utQ$,
\[
a_{\uQ \subset \utQ}: R (k_{\bigcap \utQ})_* \CK_{\utQ} \longto R (k_{\bigcap \uQ})_* \CK_{\uQ}
\]
is defined as follows: the object $R (k_{\bigcap \utQ})_* \CK_{\utQ}$ 
has support in
$\overline{e^K ( \bigcap \utQ,P(\BA_f) )}$.
According to Proposition~\ref{4B} (for $Q = \bigcap \utQ$) and Corollary~\ref{1C}~(b)
(applied to the cohomology objects of $\CV_{\tei \overline{e^K ( \bigcap \utQ,P(\BA_f) )}}$),
its restriction to $\overline{e^K ( \bigcap \utQ,P(\BA_f) )}$ is equal to 
$\CV_{\tei \overline{e^K ( \bigcap \utQ,P(\BA_f) )}}$. But $\uQ$ is contained in
$\utQ$, and hence, 
\[
e^K \bigl( \bigcap \uQ,P(\BA_f) \bigr)'
\subset \overline{e^K \bigl( \bigcap \uQ,P(\BA_f) \bigr)} 
\subset \overline{e^K \bigl( \bigcap \utQ,P(\BA_f) \bigr)} \; .
\]
Thus, there is a canonical isomorphism
\[
k_{\bigcap \uQ}^* R (k_{\bigcap \utQ})_* \CK_{\utQ} 
\isoto k_{\bigcap \uQ}^* \CV = \CK_{\uQ} \; .
\]
The map $a_{\uQ \subset \utQ}$ is defined to be adjoint to this isomorphism. \\
Componentwise adjunction defines a canonical morphism
\[
\res(\CV) \longto \deg(\CV) 
\] 
in $\CC_{(P,\FX)}^{K,BS}$.
\end{Cons}

Construction~\ref{4Con2} is functorial in $\CV$, and 
$\res(\CV) \to \deg(\CV)$ underlies a natural transformation $\res \to \deg$.

\begin{MThm} \label{4MT}
Let $\CV \in D^+(\partial M^K  (P,\FX) (\BC)^{BS})$, and assume that its
cohomology objects are locally constant.
Assume also that $K$ is neat. 
Then the canonical morphism
\[
\res(\CV) \longto \deg(\CV) 
\]
is an isomorphism. 
\end{MThm}

\begin{Proof}
Let $s \in \{1, \ldots, r \}$.
According to Corollary~\ref{4Ab}~(a), the object
\[
R(k_{Q_s})_* k_{Q_s}^* \ldots R(k_{Q_r})_* k_{Q_r}^* \CV
\]
is supported in $\overline{e^K ( Q_s \cap \ldots \cap Q_r,P(\BA_f) )}$. We claim: $(*)$ 
under our hypotheses, its restriction to 
$\overline{e^K ( Q_s \cap \ldots \cap Q_r,P(\BA_f) )}$ equals 
$\CV_{\tei \overline{e^K ( Q_s \cap \ldots \cap Q_r,P(\BA_f) )}}$.

In order to establish claim $(*)$, we apply descending induction on $s$. For $s=r$, use
Proposition~\ref{4B} (for $Q = Q_r$) and Corollary~\ref{1C}~(b). 
Let $1 \le s < r$, and assume claim $(*)$ is true for 
$R(k_{Q_{s+1}})_* k_{Q_{s+1}}^* \ldots R(k_{Q_r})_* k_{Q_r}^* \CV$. The restriction
\[
k_{Q_s}^* R(k_{Q_{s+1}})_* k_{Q_{s+1}}^* \ldots R(k_{Q_r})_* k_{Q_r}^* \CV
\] 
is therefore supported in
\[
e^K \bigl( Q_s,P(\BA_f) \bigr)' \cap 
              \overline{e^K \bigl( Q_{s+1} \cap \ldots \cap Q_r,P(\BA_f) \bigr)} 
= e^K \bigl( Q_s \cap \ldots \cap Q_r,P(\BA_f) \bigr)'
\]
(Proposition~\ref{4Aa}), and its
restriction to $e^K ( Q_s \cap \ldots \cap Q_r,P(\BA_f) )'$ equals the restriction of $\CV$.
Now apply Proposition~\ref{4B} (for $Q = Q_s \cap \ldots \cap Q_r$) and Corollary~\ref{1C}~(b).

It remains to observe that for $s=1$, claim $(*)$
implies that for any 
$\uQ = (Q_1 \prec Q_2 \prec \ldots \prec Q_r)$, 
the adjunction
\[
k_{\bigcap \uQ}^* \CV \longto 
k_{\bigcap \uQ}^* \bigl( R(k_{Q_1})_* k_{Q_1}^* \ldots R(k_{Q_r})_* k_{Q_r}^* \CV \bigr) 
\]
is an isomorphism. 
\end{Proof}


\bigskip

%
%

\section{Second Main Theorem: exceptional inverse images in the
Borel--Serre compactification}
\label{5}



Fix mixed Shimura data  $(P,\FX)$ satisfying hypotheses 
$(+)$ and $(U=0)$, an open compact sub-group $K$ of $P(\BA_f)$, 
a proper parabolic sub-group $Q$ of $P$, and 
an open sub-set $e^K ( Q,P(\BA_f) )^o$ of $\overline{e^K ( Q,P(\BA_f) )}$
containing $e^K ( Q,P(\BA_f) )$. Consider the immersions
\[
k_Q^o: e^K \bigl( Q,P(\BA_f) \bigr)^o 
\longinto \partial M^K  (P,\FX) (\BC)^{BS} 
\]
and 
\[
j_Q: \partial M^K  (P,\FX) (\BC)^{BS} - \overline{e^K \bigl( Q,P(\BA_f) \bigr)} 
\longinto \partial M^K  (P,\FX) (\BC)^{BS} 
\]
(the second of which is open).

\begin{MThm} \label{5MT}
Let $\CV \in D^+(\partial M^K  (P,\FX) (\BC)^{BS})$. Denote
by $\j$ the open immersion of $e^K ( Q,P(\BA_f))$ into 
$e^K ( Q,P(\BA_f) )^o$,
and by 
\[
\i: \partial e^K \bigl( Q,P(\BA_f) \bigr)^o 
      := e^K \bigl( Q,P(\BA_f) \bigr)^o - e^K \bigl( Q,P(\BA_f) \bigr)
\longinto e^K \bigl( Q,P(\BA_f) \bigr)^o
\]
the closed immersion complementary to $\j$. \\[0.1cm]
(a)~If $Q$ is maximal proper, then 
Construction~\ref{1Con} yields a canonical morphism of exact triangles
\[
\vcenter{\xymatrix@R-10pt{
\j_! \CV_{\tei \! e^K ( Q,P(\BA_f) )} \ar[r] \ar[d]_-{\ad} & 
 \CV_{\tei \! e^K ( Q,P(\BA_f) )^o} \ar[r] \ar@{=}[d] & 
 \i_* \CV_{\tei \! \partial e^K ( Q,P(\BA_f) )^o} \ar[r]^-{[1]} 
                                                \ar[d]^-{\i_* \i^* (k_Q^o)^*(\ad)} & 
\quad \\
R (k_Q^o)^! \CV \ar[r] & 
(k_Q^o)^* \CV \ar[r] & 
(k_Q^o)^* R (j_Q)_* j_Q^* \CV \ar[r]^-{[1]} &  
\quad
\\}}
\] 
in $D^+(e^K ( Q,P(\BA_f) )^o)$. It is an isomorphism if $K$ is neat,
and if the cohomology objects of $\CV$ are locally constant.\\[0.1cm]
(b)~Assume that $K$ is neat,
and that the cohomology objects of $\CV$ are locally constant.
If $Q$ is the intersection of $r$ distinct maximal proper sub-groups of $P$, then
\[
R (k_Q^o)^! \CV \cong \j_! \CV_{\tei e^K ( Q,P(\BA_f) )}[-(r-1)]
\]
in $D^+(e^K ( Q,P(\BA_f) )^o)$. \\[0.1cm] 
(c)~Assume that $K$ is neat,
and that the cohomology objects of $\CV$ are locally constant.
If $Q$ is not maximal proper, then the exact triangle
\[
R (k_Q^o)^! \CV \stackrel{\alpha}{\longto} (k_Q^o)^* \CV 
 \longto (k_Q^o)^* R (j_Q)_* j_Q^* \CV \stackrel{[1]}{\longto}                    
\]
is split: the morphism $\alpha$ is zero. 
\end{MThm}

The application of Main Theorem~\ref{5MT}, that will turn out to be relevant
in the sequel, concerns the immersion $k_Q^o = k_Q$ of 
\[
e^K \bigl( Q,P(\BA_f) \bigr)^o = e^K \bigl( Q,P(\BA_f) \bigr)' \; .
\]

\begin{Proofof}{Main Theorem~\ref{5MT}}
We may assume that $e^K ( Q,P(\BA_f) )^o$ is equal to $\overline{e^K ( Q,P(\BA_f) )}$;
the general result follows
from this special case by restriction. Our notation thus 
concerns
\[
k_Q^o: e^K \bigl( Q,P(\BA_f) \bigr)^o = \overline{e^K \bigl( Q,P(\BA_f) \bigr)}
\longinto \partial M^K  (P,\FX) (\BC)^{BS} \; ,
\]
\[
\j: e^K \bigl( Q,P(\BA_f) \bigr) \longinto \overline{e^K \bigl( Q,P(\BA_f) \bigr)} \; ,
\]
and  
\[
\i: \partial \overline{e^K \bigl( Q,P(\BA_f) \bigr)} 
      = \overline{e^K \bigl( Q,P(\BA_f) \bigr)} - e^K \bigl( Q,P(\BA_f) \bigr)
\longinto \overline{e^K \bigl( Q,P(\BA_f) \bigr)} \; .
\]
(a): in order to apply Proposition~\ref{1D}~(b) to 
\[
X = \partial M^K  (P,\FX) (\BC)^{BS} \; ,
\]
$j = j_Q$, and $Z^0 = e^K ( Q,P(\BA_f) )$ (noting that $Z^0$ is indeed open in
$X$ as $Q$ is assumed maximal proper), we need to check that the immersion
of 
\[
U = \partial M^K  (P,\FX) (\BC)^{BS} - \overline{e^K \bigl( Q,P(\BA_f) \bigr)}
\]
into 
\[
X - Z^0 = \partial M^K  (P,\FX) (\BC)^{BS} - e^K \bigl( Q,P(\BA_f) \bigr)
\]
is contractible. But this is precisely the case $R = P$ of \cite[Prop.~5.5~(b)]{W1}.

\noindent (b): we apply induction on $r$, the case $r=1$ resulting from (a).
If $r \ge 2$, let $Q_1,\ldots,Q_r$ be the distinct maximal proper sub-groups of $P$
containing $Q$; we thus have $Q = Q_1 \cap \ldots \cap Q_r \,$. Define
$\tilde{Q} := Q_2 \cap \ldots \cap Q_r \,$,
\[
k_{\tilde{Q}}^o: \overline{e^K \bigl( \tilde{Q},P(\BA_f) \bigr)}
\longinto \partial M^K  (P,\FX) (\BC)^{BS} \; ,
\]
and
\[
\tilde{\j}: e^K \bigl( \tilde{Q},P(\BA_f) \bigr) 
\longinto \overline{e^K \bigl( \tilde{Q},P(\BA_f) \bigr)} \; .
\]
By the induction hypothesis, we have
\[
R (k_{\tilde{Q}}^o)^! \CV \cong \tilde{\j}_! \CV_{\tei e^K ( \tilde{Q},P(\BA_f) )}[-(r-2)] \; .
\]
Now $R (k_Q^o)^! = R (k')^! \circ R (k_{\tilde{Q}}^o)^! \,$, for the closed immersion $k'$ of
$\overline{e^K ( Q,P(\BA_f) )}$ into $\overline{e^K ( \tilde{Q},P(\BA_f) )}$.
Letting
\[
\j': \overline{e^K \bigl( \tilde{Q},P(\BA_f) \bigr)} - \overline{e^K \bigl( Q,P(\BA_f) \bigr)}
\longinto \overline{e^K \bigl( \tilde{Q},P(\BA_f) \bigr)}
\]
denote the open immersion complementary to $k'$, we get $\tilde{\j} = \j' \circ \j''$, where
$\j''$ is the open immersion of $e^K ( \tilde{Q},P(\BA_f) )$ into
$\overline{e^K ( \tilde{Q},P(\BA_f) )} - \overline{e^K ( Q,P(\BA_f) )}$. As $\j'$ and
$k'$ are complementary to each other, we have
\[
R (k')^! \circ \j'_! = (k')^* \circ R \j'_* [-1] 
\]
\cite[(1.4.6.4)]{BBD}. From what was said so far, we conclude that $R (k_Q^o)^! \CV$ 
is isomorphic to
\[
R (k')^! \circ R (k_{\tilde{Q}}^o)^! \CV 
\cong (k')^* \circ R \j'_* \bigl( \j_!'' \CV_{\tei e^K ( \tilde{Q},P(\BA_f) )} \bigr)[-(r-1)]
\; .
\]
Let us compute $R \j'_* ( \j_!'' \CV_{\tei e^K ( \tilde{Q},P(\BA_f) )} )$.
The object $\CW:= \j_!'' \CV_{\tei e^K ( \tilde{Q},P(\BA_f) )}$ is part of
the localization exact triangle
\[
\CW
\longto \CV_{\tei \overline{e^K ( \tilde{Q},P(\BA_f) )} - \overline{e^K ( Q,P(\BA_f) )}}
\longto \i''_* 
      \CV_{\tei \partial \overline{e^K ( \tilde{Q},P(\BA_f) )} - \overline{e^K ( Q,P(\BA_f) )}}
\longto \CW [1] \; ,
\]
where $\i''$ is the closed immersion of 
$\partial \overline{e^K ( \tilde{Q},P(\BA_f) )} - \overline{e^K ( Q,P(\BA_f) )}$.
Now $\i''$ and $\j'$ occur in the commutative diagram
\[
\vcenter{\xymatrix@R-10pt{
\partial \overline{e^K \bigl(\tilde{Q},P(\BA_f)\bigr)} - \overline{e^K \bigl(Q,P(\BA_f)\bigr)}
   \ar@{^(->}[r]^-{\i''} \ar@{_(->}[d]_-{\partial \j'} & 
\overline{e^K \bigl(\tilde{Q},P(\BA_f)\bigr)} - \overline{e^K \bigl(Q,P(\BA_f)\bigr)} 
   \ar@{_(->}[d]_-{\j'} &  \\
\partial \overline{e^K \bigl(\tilde{Q},P(\BA_f)\bigr)} - e^K \bigl(Q,P(\BA_f)\bigr)
   \ar@{^(->}[r]^-I & 
\overline{e^K \bigl(\tilde{Q},P(\BA_f)\bigr)} 
\\}}
\] 
(the immersions $\partial \j'$ and $I$ are defined by the diagram).
Both $\j'$ and $\partial \j'$ are contractible \cite[Prop.~5.5~(a) and (b)]{W1}
(note that $\tilde{Q}$ contains $Q$ as a maximal proper parabolic). 
According to Proposition~\ref{1B}, we therefore have
\[
R \j'_* \CV_{\tei \overline{e^K ( \tilde{Q},P(\BA_f) )} - \overline{e^K ( Q,P(\BA_f) )}}
= \CV_{\tei \overline{e^K ( \tilde{Q},P(\BA_f) )}}   \quad\quad\quad\quad (\ast)
\]
and
\[
R (\partial \j)'_* \CV_{\tei \overline{e^K ( \tilde{Q},P(\BA_f) )} - \overline{e^K ( Q,P(\BA_f) )}}
= \CV_{\tei \overline{e^K ( \tilde{Q},P(\BA_f) )} - e^K (Q,P(\BA_f))} \; ,
\]
hence 
\[
R \j'_* \i''_* 
      \CV_{\tei \partial \overline{e^K ( \tilde{Q},P(\BA_f) )} - \overline{e^K ( Q,P(\BA_f) )}}
= I_* \CV_{\tei \overline{e^K ( \tilde{Q},P(\BA_f) )} - e^K (Q,P(\BA_f))} \; .
\]   
The desired object $R \j'_* ( \j_!'' \CV_{\tei e^K ( \tilde{Q},P(\BA_f) )} ) = R \j'_* \CW$ 
sits in the image under $R \j'_*$ of the above localization exact triangle, which
by what was just said, equals
\[
R \j'_* \CW  
\longto \CV_{\tei \overline{e^K ( \tilde{Q},P(\BA_f) )}}
\longto I_* \CV_{\tei \overline{e^K ( \tilde{Q},P(\BA_f) )} - e^K (Q,P(\BA_f))}
\longto R \j'_* \CW [1] \; ,
\]
Restriction under $k'$ of the latter gives the exact triangle
\[
(k')^* R \j'_* \CW  
\longto \CV_{\tei \overline{e^K ( Q,P(\BA_f) )}}
\longto \i_* \CV_{\tei \partial \overline{e^K ( Q,P(\BA_f) )}}
\longto (k')^* R \j'_* \CW [1] \; ,
\] 
whence 
\[
(k')^* R \j'_* (\j_!'' \CV_{\tei e^K ( \tilde{Q},P(\BA_f) )} )) = (k')^* R \j'_* \CW
\cong \j_! \CV_{\tei e^K ( Q,P(\BA_f) )} \; .
\]
(c): keep the notation from the proof of (b), case $r \ge 2$, in particular
\[
\j': \overline{e^K \bigl( \tilde{Q},P(\BA_f) \bigr)} - \overline{e^K \bigl( Q,P(\BA_f) \bigr)}
\longinto \overline{e^K \bigl( \tilde{Q},P(\BA_f) \bigr)} \; ,
\]
and 
\[
k': \overline{e^K \bigl( Q,P(\BA_f) \bigr)}
\longinto \overline{e^K \bigl( \tilde{Q},P(\BA_f) \bigr)} \; .
\]
We have 
\[
j_Q \circ k_{\tilde{Q}}' = k_{\tilde{Q}}^o \circ \j': 
\overline{e^K \bigl( \tilde{Q},P(\BA_f) \bigr)} - \overline{e^K \bigl( Q,P(\BA_f) \bigr)}
\longinto \partial M^K  (P,\FX) (\BC)^{BS} \; ,
\]
where $k_{\tilde{Q}}'$ denotes the open immersion of 
$\overline{e^K ( \tilde{Q},P(\BA_f) )} - \overline{e^K ( Q,P(\BA_f) )}$ into 
$\partial M^K  (P,\FX) (\BC)^{BS} - \overline{e^K ( Q,P(\BA_f) )}$.
Consider the composition of adjunctions
\[
(k_Q^o)^* \CV \longto 
(k_Q^o)^* R (j_Q)_* j_Q^* \CV \longto
(k')^* R \j'_* (\j')^* (k_{\tilde{Q}}^o)^* \CV \; .
\]
According to formula $(\ast)$ from the proof of (b), 
\[
R \j'_* (\j')^* (k_{\tilde{Q}}^o)^* \CV = (k_{\tilde{Q}}^o)^* \CV \; ,
\]
and hence, 
\[
(k')^* R \j'_* (\j')^* (k_{\tilde{Q}}^o)^* \CV = (k_Q^o)^* \CV \; .
\]
Therefore, the adjunction 
$(k_Q^o)^* \CV \longto (k_Q^o)^* R (j_Q)_* j_Q^* \CV$ is split monomorphic.
\end{Proofof}

\begin{Rem}
The isomorphism
\[
R (k_Q^o)^! \CV \isoto \j_! \CV_{\tei e^K ( Q,P(\BA_f) )}[-(r-1)]
\]
from Main~Theorem~\ref{5MT}~(b) depends on the choice of an order $Q_1,\ldots,Q_r$
of the distinct maximal proper sub-groups containing $Q$. Changing the order
by a permutation $\sigma$ gives rise to another isomorphism. A careful analysis
of the proof shows that the two isomorphisms only differ by a sign, 
which equals $\sgn(\sigma)$. 
\end{Rem}

The rest of this section is concerned with compatibility of the isomorphisms
from Main Theorems~\ref{4MT} and \ref{5MT}. We keep our
mixed Shimura data  $(P,\FX)$ satisfying hypotheses 
$(+)$ and $(U=0)$, as well as 
the open compact sub-group $K$ of $P(\BA_f)$. Fix in addition
an element
$\uQ = (Q_1 \prec Q_2 \prec \ldots \prec Q_r)$ of $\CC_{(P,\FX)}$, such that
$Q_1$ is maximal proper, and $r \ge 2$. We denote by $\i$ the closed immersion
of $\partial e^K ( Q_1,P(\BA_f) )'$ into $e^K ( Q_1,P(\BA_f) )'$. \\

Let $\CV \in D^+(\partial M^K  (P,\FX) (\BC)^{BS})$.
On the one hand, we have a morphism
\[
\i_* \i^* k_{Q_1}^*(\ad): \i_* \CV_{\tei \! \partial e^K ( Q_1,P(\BA_f) )'} 
\longto k_{Q_1}^* R (j_{Q_1})_* j_{Q_1}^* \CV   
\]
in $D^+(e^K ( Q_1,P(\BA_f) )')$ (Main Theorem~\ref{5MT}~(a), applied to $Q = Q_1$
and $e^K ( Q_1,P(\BA_f) )^o = e^K ( Q_1,P(\BA_f) )'$). \\

On the other hand, if $K$ is neat,
and if the cohomology objects of $\CV$ are locally constant,
then the component $\uQ$ of the canonical morphism
\[
\res(\CV) \longto \deg(\CV) 
\] 
in $\CC_{(P,\FX)}^{K,BS}$ from Section~\ref{4} is equal to the adjunction
\[
\ad: k_{\bigcap \uQ}^* \CV
\longto k_{\bigcap \uQ}^* \bigl( R(k_{Q_1})_* k_{Q_1}^* R(k_{Q_2})_* k_{Q_2}^* \ldots R(k_{Q_r})_* k_{Q_r}^* \CV \bigr) \; .
\]
Writing $k_{\uQ,Q_1}$ for the
immersion of $e^K ( \bigcap \uQ,P(\BA_f) )'$ into 
$e^K ( Q_1,P(\BA_f) )'$
(recall that by Proposition~\ref{2C}, we have $\adm(\bigcap \uQ) = Q_1$), 
this adjunction induces a morphism $(k_{\uQ,Q_1})_*(\ad)$
in $D^+(e^K ( Q_1,P(\BA_f) )')$. \\

The respective sources and targets of $\i_* \i^* k_{Q_1}^*(\ad)$ and $(k_{\uQ,Q_1})_*(\ad)$
are related by morphisms: indeed,
\[
\i_* \CV_{\tei \! \partial e^K ( Q_1,P(\BA_f) )'} \longto (k_{\uQ,Q_1})_* k_{\bigcap \uQ}^* \CV
\]
is adjoint to the identity
\[
k_{\uQ,Q_1}^* \bigl( \i_* \CV_{\tei \! \partial e^K ( Q_1,P(\BA_f) )'} \bigr)
\longto k_{\bigcap \uQ}^* \CV
\]
(note that $e^K ( \bigcap \uQ,P(\BA_f) )' \subset \partial e^K ( Q_1,P(\BA_f) )'$
as $\bigcap \uQ$ is properly contained in $Q_1$, thanks to our assumption $r \ge 2$).
The morphism 
\[
k_{Q_1}^* R (j_{Q_1})_* j_{Q_1}^* \CV 
\longto 
(k_{\uQ,Q_1})_* k_{\bigcap \uQ}^* \bigl( R(k_{Q_1})_* k_{Q_1}^* R(k_{Q_2})_* k_{Q_2}^* \ldots 
             R(k_{Q_r})_* k_{Q_r}^* \CV \bigr) 
\]
is adjoint to a morphism
\[
\beta: k_{\bigcap \uQ}^* R (j_{Q_1})_* j_{Q_1}^* \CV 
\longto 
k_{\bigcap \uQ}^* \bigl( R(k_{Q_1})_* k_{Q_1}^* R(k_{Q_2})_* k_{Q_2}^* \ldots 
             R(k_{Q_r})_* k_{Q_r}^* \CV \bigr) \; ,
\]
defined as follows: since 
$e^K ( \bigcap \uQ,P(\BA_f) )' \subset e^K ( Q_1,P(\BA_f) )'$, the target
of $\beta$ equals
\[
k_{\bigcap \uQ}^* \bigl( R(k_{Q_2})_* k_{Q_2}^* \ldots 
             R(k_{Q_r})_* k_{Q_r}^* \CV \bigr) \; .
\]
Now the object $k_{Q_2}^* \ldots R(k_{Q_r})_* k_{Q_r}^* \CV$
has support in $e^K ( Q_2 \cap \ldots \cap Q_r,P(\BA_f) )'$ (Corollary~\ref{4Ab}~(a)),
and the latter sub-set of $\partial M^K  (P,\FX) (\BC)^{BS}$ is contained in
\[
e^K ( Q_2,P(\BA_f) )' \subset
\partial M^K  (P,\FX) (\BC)^{BS} - \overline{e^K \bigl( Q_1,P(\BA_f) \bigr)} 
\]
\cite[Prop.~5.7]{W1}.
It follows that the adjunction 
\[
R(k_{Q_2})_* k_{Q_2}^* \ldots R(k_{Q_r})_* k_{Q_r}^* \CV 
\longto R (j_{Q_1})_* j_{Q_1}^* R(k_{Q_2})_* k_{Q_2}^* \ldots R(k_{Q_r})_* k_{Q_r}^* \CV
\]
is an isomorphism. Consequently, the adjunction 
\[
\CV \longto R(k_{Q_2})_* k_{Q_2}^* \ldots R(k_{Q_r})_* k_{Q_r}^* \CV
\]
induces a morphism
\[
\gamma : R (j_{Q_1})_* j_{Q_1}^* \CV 
\longto R(k_{Q_2})_* k_{Q_2}^* \ldots R(k_{Q_r})_* k_{Q_r}^* \CV \; .
\]
Define $\beta := k_{\bigcap \uQ}^* (\gamma)$. Altogether, we end up with a diagram
\[
\vcenter{\xymatrix@R-10pt{
\i_* \CV_{\tei \! \partial e^K ( Q_1,P(\BA_f) )'} \ar[r]^-{\ad}
                                      \ar[d]^-{\i_* \i^* k_{Q_1}^*(\ad)} & 
(k_{\uQ,Q_1})_* k_{\bigcap \uQ}^* \CV  
                                                        \ar[d]_-{(k_{\uQ,Q_1})_*(\ad)} \\
k_{Q_1}^* R (j_{Q_1})_* j_{Q_1}^* \CV \ar[r]^-{\ad(\beta)} &  
(k_{\uQ,Q_1})_* k_{\bigcap \uQ}^* 
\bigl( R(k_{Q_1})_* k_{Q_1}^* R(k_{Q_2})_* k_{Q_2}^* \ldots R(k_{Q_r})_* k_{Q_r}^* \CV \bigr)
\\}}
\] 
in $D^+(e^K ( Q_1,P(\BA_f) )')$.

\begin{Prop} \label{5B}
Assume that $K$ is neat,
and that the cohomology objects of $\CV$ are locally constant.
The above diagram is commutative. 
\end{Prop}

Recall that according to Main Theorems~\ref{5MT}~(a) and \ref{4MT},
the vertical arrows of the diagram are isomorphisms.

\medskip

\begin{Proofof}{Proposition~\ref{5B}}
Write $\utQ := (Q_2 \prec \ldots \prec Q_r)$. According to Proposition~\ref{4Aa}, we have
\[
e^K \bigl( Q_1,P(\BA_f) \bigr)' \cap 
              \overline{e^K \bigl( \bigcap \utQ,P(\BA_f) \bigr)} 
              = e^K \bigl( \bigcap \uQ,P(\BA_f) \bigr)' \; .
\]
Applying Proposition~\ref{1D}~(a) to the closed immersion $f$ of
\[
\overline{e^K \bigl( \bigcap \utQ,P(\BA_f) \bigr)}
- \Bigl( 
\overline{e^K \bigl( Q_1,P(\BA_f) \bigr)} - e^K \bigl( Q_1,P(\BA_f) \bigr)'
\Bigr) 
\]
into 
\[
\partial M^K  (P,\FX) (\BC)^{BS} 
- \Bigl( 
\overline{e^K \bigl( Q_1,P(\BA_f) \bigr)} - e^K \bigl( Q_1,P(\BA_f) \bigr)'
\Bigr) \; ,
\] 
we get the commutativity of the diagram
\[
\vcenter{\xymatrix@R-10pt{
\i_* \CV_{\tei \! \partial e^K ( Q_1,P(\BA_f) )'} \ar[r]^-{\ad}
                                      \ar[d]_-{\i_* \i^* k_{Q_1}^*(\ad)} & 
(k_{\uQ,Q_1})_* k_{\bigcap \uQ}^* \CV  
                               \ar[d]^-{(k_{\uQ,Q_1})_*k_{\bigcap \uQ}^*(\ad)f^*} \\
k_{Q_1}^* R (j_{Q_1})_* j_{Q_1}^* \CV \ar[r] &  
(k_{\uQ,Q_1})_* k_{\bigcap \uQ}^* 
\bigl( R j_*' (j')^* \CV \bigr)
\\}} \; ,
\] 
where $j'$ denotes the immersion of 
\[
U' = \Bigl( \partial M^K  (P,\FX) (\BC)^{BS} - \overline{e^K \bigl( Q_1,P(\BA_f) \bigr)} \Bigr)
\cap \overline{e^K \bigl( \bigcap \utQ,P(\BA_f) \bigr)}
\]
into $\partial M^K  (P,\FX) (\BC)^{BS}$.
According to \cite[Prop.~5.7]{W1}, the set $U'$ contains $e^K ( \bigcap \utQ,P(\BA_f) )'$.
This allows (\emph{cf.} the definition of the morphism
$\gamma$) to complete the diagram commutatively to
\[
\vcenter{\xymatrix@R-10pt{
\i_* \CV_{\tei \! \partial e^K ( Q_1,P(\BA_f) )'} \ar[r]^-{\ad}
                                      \ar[d]_-{\i_* \i^* k_{Q_1}^*(\ad)} & 
(k_{\uQ,Q_1})_* k_{\bigcap \uQ}^* \CV  
                               \ar[d]^-{(k_{\uQ,Q_1})_*k_{\bigcap \uQ}^*(\ad)f^*} \\
k_{Q_1}^* R (j_{Q_1})_* j_{Q_1}^* \CV \ar[r] \ar[dr]&  
(k_{\uQ,Q_1})_* k_{\bigcap \uQ}^* 
\bigl( R j_*' (j')^* \CV \bigr) \ar[d] \\
&
(k_{\uQ,Q_1})_* k_{\bigcap \uQ}^* 
R(k_{Q_1})_* k_{Q_1}^* R(k_{Q_2})_* k_{Q_2}^* \ldots R(k_{Q_r})_* k_{Q_r}^* \CV 
\\}} \; .
\]  
We leave it to the reader to verify that the morphisms 
\[
k_{Q_1}^* R (j_{Q_1})_* j_{Q_1}^* \CV
\longto (k_{\uQ,Q_1})_* k_{\bigcap \uQ}^* 
R(k_{Q_1})_* k_{Q_1}^* R(k_{Q_2})_* k_{Q_2}^* \ldots R(k_{Q_r})_* k_{Q_r}^* \CV
\]
and
\[
(k_{\uQ,Q_1})_* k_{\bigcap \uQ}^* \CV
\longto (k_{\uQ,Q_1})_* k_{\bigcap \uQ}^* 
R(k_{Q_1})_* k_{Q_1}^* R(k_{Q_2})_* k_{Q_2}^* \ldots R(k_{Q_r})_* k_{Q_r}^* \CV
\]
in this latter diagram are equal to $\ad(\beta)$ and $(k_{\uQ,Q_1})_*(\ad)$, respectively.
\end{Proofof}


\bigskip

%
%

\section{Consequences for the
Baily--Borel compactification}
\label{6}



In this section, we reformulate 
the main results from the two preceding sections
in terms of the Baily--Borel compactification 
(Theorems~\ref{6MT1} and \ref{6MT2} and their corollaries). 
Fix pure Shimura data  $(P,\FX) = (G,\FX)$ satisfying hypothesis 
$(+)$, and an open compact sub-group $K$ of $G(\BA_f)$. 

\begin{Def} \label{6A}
Define the \emph{boundary of the Baily--Borel compactification} of $M^K  (G,\FX)$ to be
the closed complement 
\[
\partial M^K  (G,\FX)^* := M^K  (G,\FX)^* - M^K  (G,\FX) \; ,
\] 
of $M^K  (G,\FX)$ in $M^K  (G,\FX)^*$, equipped with the reduced scheme structure.
Denote by $i$ the closed immersion of $\partial M^K  (G,\FX)^*$ into $M^K  (G,\FX)^*$.
\end{Def}

Recall that $M^K  (G,\FX)$ is the unique open stratum of the canonical
stratification of $M^K  (G,\FX)^*$ (Definition~\ref{3D}), and that 
for each admissible parabolic sub-group $Q_s$ of $G$, the immersion of the stratum
$M^K (Q_s,\FX)$ into $M^K  (G,\FX)^*$ is denoted by $i_{Q_s}$. \\

The following is the analogue of Construction~\ref{4Con1} for the Baily--Borel
compactification.

\begin{Cons} \label{6Con}
Let $\CW \in D^+(\partial M^K  (G,\FX)^* (\BC))$. Define an object
$\deg(\CW)$ of $\CC_{(G,\FX)}^{K,*}$ as follows.
\[
\deg(\CW) = \bigl( \CK_{\uQ} , a_{\uQ \subset \utQ} \bigr)_{\uQ, \uQ \subset \utQ} \; ,
\]
where for $\uQ = (Q_1 \prec Q_2 \prec \ldots \prec Q_r) \in \CC_{(G,\FX)}$,
(hence $b(\uQ) = Q_1 \,$,)
\[
\CK_{\uQ} := i_{Q_1}^* R(i_{Q_2})_* i_{Q_2}^* \ldots R(i_{Q_r})_* i_{Q_r}^* \CW 
\in D^+ \bigl(M^K (Q_1,\FX)(\BC) \bigr) \; ,
\]
and for $\uQ \subset \utQ$,
\[
a_{\uQ \subset \utQ}: R (i_{b(\utQ)})_* \CK_{\utQ} \longto R (i_{b(\uQ)})_* \CK_{\uQ}
\]
is the adjunction.
\end{Cons}

Note that all immersions $i_{Q_s}$ occurring in Construction~\ref{6Con} are associated
to proper parabolic sub-groups of $P$. Therefore, their images are contained in
$\partial M^K  (G,\FX)^*$. We shall therefore 
use the same symbol $i_{Q_s}$ to denote
the immersion of $M^K (Q_s,\FX)$ into 
the boundary $\partial M^K  (G,\FX)^*$ as well. \\

We leave it to the reader to complete Construction~\ref{6Con} to a functor
\[
\deg: D^+ \bigl( \partial M^K  (G,\FX)^* (\BC) \bigr) \longto \CC_{(G,\FX)}^{K,*} \; .
\]

\begin{Rem} \label{6Rem}
The data 
\[
a_{\uQ \subset \utQ}: R (i_{b(\utQ)})_* \CK_{\utQ} \longto R (i_{b(\uQ)})_* \CK_{\uQ}
\]
occurring in Construction~\ref{6Con} are equivalent to their adjoints
\[
\ad(a_{\uQ \subset \utQ}): i_{b(\uQ)}^*R (i_{b(\utQ)})_* \CK_{\utQ} \longto  \CK_{\uQ} \; .
\]
If $b(\uQ) = b(\utQ)$, then $\ad(a_{\uQ \subset \utQ})$ is simply a morphism
$\CK_{\utQ} \to  \CK_{\uQ}$.
Else, we have $\utQ = (Q_2 \prec \ldots \prec Q_r)$ and 
$Q_1:= b(\uQ) \prec Q_2 \,$.
Putting $\uQ':= (Q_1 \prec Q_2 \prec \ldots \prec Q_r)$, we have
$\uQ \subset \uQ' \subset \utQ$, and the equality
\[
\ad(a_{\uQ \subset \utQ}) = \ad(a_{\uQ \subset \uQ'}) \; .
\]
This shows that $\deg(\CW) \in \CC_{(G,\FX)}^{K,*}$ is completely determined by 
the collection of all
$i^*_{Q_1} \deg(\CW)$ (Definition~\ref{3E}(c)),
where $Q_1$ runs through the proper admissible parabolic sub-groups of $G$.
\end{Rem}

Recall the functor
\[
Rp^K_*: \CC_{(G,\FX)}^{K,BS} \longto \CC_{(G,\FX)}^{K,*} 
\]
from Proposition~\ref{3F}, which is associated to
the continuous map 
\[
p^K: M^K  (G,\FX) (\BC)^{BS} \to M^K (G,\FX)^* (\BC)
\]
from \cite[Cor.~4.18]{W1}. 
The latter having compact source and target
(\cite[Thm.~9.3]{BS}, \cite[Thm.~3.14]{W1}, \cite[Chap.~II, Thm.~2]{AMRT}), it is proper.
This fact will be frequently employed in the sequel, as will the following.

\begin{Prop}[{\cite[Thm.~5.10~(a)]{W1}}] \label{6Aa}
Let $Q_s$ be an admissible parabolic sub-group of $G$.
Then the pre-image of 
$M^K (Q_s,\FX)(\BC) \subset M^K (G,\FX)^* (\BC)$ under $p^K$ equals 
$e^K ( Q_s,G(\BA_f) )' \subset M^K  (G,\FX) (\BC)^{BS}$.
\end{Prop}

\begin{Prop} \label{6B}
Proper base change induces a canonical isomorphism
\[
Rp^K_* \circ \deg \isoto \deg \circ R p^K_*  
\]
of functors
$D^+(\partial M^K  (G,\FX) (\BC)^{BS}) \to \CC_{(G,\FX)}^{K,*}$. 
\end{Prop}

\begin{Proof}
Write 
\[
\deg(\CV) = \bigl( \CK_{\uQ} , a_{\uQ \subset \utQ} \bigr)_{\uQ, \uQ \subset \utQ}
\]
and 
\[
\deg \bigl( Rp^K_* \CV \bigr) 
= \bigl( \CL_{\uQ} , b_{\uQ \subset \utQ} \bigr)_{\uQ, \uQ \subset \utQ} \; .
\]
Let $\uQ = (Q_1 \prec Q_2 \prec \ldots \prec Q_r) \in \CC_{(G,\FX)}$.
We have
\[
\CK_{\uQ} =  
k_{\bigcap \uQ}^* \bigl( R(k_{Q_1})_* k_{Q_1}^* R(k_{Q_2})_* k_{Q_2}^* \ldots R(k_{Q_r})_* k_{Q_r}^* \CV \bigr) 
\]
(Construction~\ref{4Con1}), 
and by Corollary~\ref{4Ab}~(b), the adjunction
\[
\ad : k_{Q_1}^* \ldots R(k_{Q_r})_* k_{Q_r}^* \CV \longto (k_{\uQ,Q_1})_* \CK_{\uQ}
\]
is an isomorphism. The restriction of $p^K$ to $e^K ( \bigcap \uQ,G(\BA_f) )'$
factors through $k_{\uQ,Q_1} \,$; hence
\[
Rp^K_*(\ad) : Rp^K_* \bigl( k_{Q_1}^* \ldots R(k_{Q_r})_* k_{Q_r}^* \CV \bigr)
\isoto Rp^K_* \bigl( \CK_{\uQ} \bigr) = \CL_{\uQ} \; .
\]
But by proper base change \cite[Thm.~6.2]{I} and
Proposition~\ref{6Aa}, we have 
$Rp^K_* k_{Q_s}^* \isoto i_{Q_s}^* Rp^K_*$ canonically, for all
$s = 1,\ldots,r$. 
\end{Proof}

\begin{Thm} \label{6MT1}
Let $\CV \in D^+(\partial M^K  (G,\FX) (\BC)^{BS})$, and assume that its
cohomology objects are locally constant.
Assume also that $K$ is neat.
Then there is a cano\-nical isomorphism
\[
Rp^K_* (\res(\CV)) \isoto \deg \bigl( Rp^K_*\CV \bigr)  
\]
in $\CC_{(G,\FX)}^{K,*}$. 
\end{Thm}

\begin{Proof}
Combine Proposition~\ref{6B} and Main Theorem~\ref{4MT}.
\end{Proof}

The application of Theorem~\ref{6MT1} that we have in mind, concerns
objects of $D^+(\partial M^K  (G,\FX) (\BC)^{BS})$ coming from the open stratum
\[
M^K  (G,\FX) (\BC) = e^K ( G,G(\BA_f)) 
\]
of the canonical stratification of $M^K  (G,\FX) (\BC)^{BS}$ (Definition~\ref{3A}).
Recall that $k_G$ denotes
the open immersion of $M^K  (G,\FX) (\BC) (= e^K ( G,G(\BA_f))')$ into
$M^K  (G,\FX) (\BC)^{BS}$.

\begin{Prop} \label{6Ba}
Let $\CV \in D^+(M^K  (G,\FX) (\BC))$, and assume that its
cohomology objects are locally constant.
Assume also that $K$ is neat.
Then the 
cohomology objects of 
\[
R(k_G)_* \CV \in D^+ \bigl( M^K  (G,\FX) (\BC)^{BS} \bigr)
\]
are locally constant as well.
\end{Prop}

\begin{Proof}
Since $K$ is assumed neat, the open immersion $k_G$ is contractible
\cite[Prop.~5.4~(c)]{W1}. Now apply Corollary~\ref{1C}~(c).
\end{Proof}

\begin{Cor} \label{6C}
Let $\CV \in D^+(M^K  (G,\FX) (\BC))$, and assume that its
cohomology objects are locally constant.
Assume also that $K$ is neat.
Denote by $\CV^{\partial BS}$ the restriction to $\partial M^K  (G,\FX) (\BC)^{BS}$
of the direct image $R(k_G)_* \CV$. 
Then there is a cano\-nical isomorphism
\[
Rp^K_* \bigl( \res \bigr( \CV^{\partial BS} \bigr) \bigr) 
\isoto \deg \bigl( i^*R(i_G)_* \CV \bigr)  
\]
in $\CC_{(G,\FX)}^{K,*}$. 
\end{Cor}

\begin{Proof}
By Proposition~\ref{6Ba}, the 
cohomology objects of $R(k_G)_* \CV$
are locally constant. Hence the same is true for $\CV^{\partial BS}$.
According to Theorem~\ref{6MT1},
\[
Rp^K_* \bigl( \res \bigr( \CV^{\partial BS} \bigr) \bigr) 
\isoto \deg \bigl( Rp^K_* \CV^{\partial BS} \bigr) \; . 
\]
But by proper base change \cite[Thm.~6.2]{I},
\[
Rp^K_* \CV^{\partial BS} \isoto i^*R(i_G)_* \CV \; .
\]
\end{Proof}

\begin{Thm} \label{6MT2}
Let $\CV \in D^+(\partial M^K  (G,\FX) (\BC)^{BS})$,
and $Q_1 \ne G$ an admissible parabolic sub-group of $G$.
Denote (\emph{cf.} Main Theorem~\ref{5MT})
by $\j$ the open immersion of $e^K ( Q_1,G(\BA_f))$ into 
$e^K ( Q_1,G(\BA_f) )'$, and by 
\[
\i: \partial e^K \bigl( Q_1,G(\BA_f) \bigr)' \longinto e^K \bigl( Q_1,G(\BA_f) \bigr)'
\]
the closed immersion complementary to $\j$. Also, let
\[
h_{Q_1}: \partial M^K  (G,\FX)^* - \overline{M^K \bigl( Q_1,\FX \bigr)} 
\longinto \partial M^K  (G,\FX)^*
\]
denote
the open immersion of the complement of $\overline{M^K \bigl( Q_1,\FX \bigr)}$. \\[0.1cm]
(a)~Assume that $Q_1$ is maximal proper. Then the application of $Rp_*^K$ to the morphism
from Main Theorem~\ref{5MT}~(a) yields  
a canonical morphism of exact triangles from
\[
Rp^K_* \bigl( \j_! \CV_{\tei \! e^K ( Q_1,G(\BA_f) )} \bigr)  \to
Rp^K_* \bigl( \CV_{\tei \! e^K ( Q_1,G(\BA_f) )'} \bigr) \to 
Rp^K_* \bigl( \i_* \CV_{\tei \! \partial e^K ( Q_1,G(\BA_f) )'} \bigr) \stackrel{[1]}{\to} 
\]
to 
\[
R i_{Q_1}^! Rp^K_* \CV \longto
i_{Q_1}^* Rp^K_* \CV \longto
i_{Q_1}^* R (h_{Q_1})_* h_{Q_1}^* Rp^K_* \CV \stackrel{[1]}{\longto}  
\] 
in $D^+(M^K ( Q_1,\FX )(\BC))$. It is an isomorphism if $K$ is neat,
and if the cohomology objects of $\CV$ are locally constant.\\[0.1cm]
(b)~Assume that $K$ is neat,
and that the cohomology objects of $\CV$ are locally constant.
If $Q_1$ is the intersection of $r$ distinct maximal proper sub-groups of $G$, then
\[
R i_{Q_1}^! Rp^K_* \CV 
\cong Rp^K_* \bigl( \j_! \CV_{\tei \! e^K ( Q_1,G(\BA_f) )} \bigr)[-(r-1)]
\]
in $D^+(M^K ( Q_1,\FX )(\BC))$. \\[0.1cm] 
(c)~Assume that $K$ is neat,
and that the cohomology objects of $\CV$ are locally constant.
If $Q_1$ is not maximal proper, then the exact triangle
\[
R i_{Q_1}^! Rp^K_* \CV \stackrel{\alpha}{\longto} i_{Q_1}^* Rp^K_* \CV 
\longto i_{Q_1}^* R (h_{Q_1})_* h_{Q_1}^* Rp^K_* \CV \stackrel{[1]}{\longto}                    
\]
is split: the morphism $\alpha$ is zero. 
\end{Thm}

\begin{Proof}
First, apply the functor $Rp_*^K$ to the isomorphisms from Main Theorem~\ref{5MT}.
Then, use proper base change \cite[Thm.~6.2]{I} and Proposition~\ref{6Aa}.
\end{Proof}

\begin{Cor} \label{6D}
Let $\CV \in D^+(M^K  (G,\FX) (\BC))$, and assume that its
cohomology objects are locally constant.
Assume also that $K$ is neat.
Denote by $\CV^{\partial BS}$ the restriction to $\partial M^K  (G,\FX) (\BC)^{BS}$
of $R(k_G)_* \CV$. Let $Q_1 \ne G$ an admissible parabolic sub-group of $G$,
and keep the notations $\j$, $\i$, and $h_{Q_1}$ from Theorem~\ref{6MT2}. \\[0.1cm]
(a)~Assume that $Q_1$ is maximal proper. Then the isomorphism
from Theorem~\ref{6MT2}~(a) induces 
a canonical isomorphism of exact triangles from
\[
Rp^K_* \bigl( \j_! \CV^{\partial BS}_{\tei \! e^K ( Q_1,G(\BA_f) )} \bigr) \to 
Rp^K_* \bigl( \CV^{\partial BS}_{\tei \! e^K ( Q_1,G(\BA_f) )'} \bigr) \to 
Rp^K_* \bigl( \i_* \CV^{\partial BS}_{\tei \! \partial e^K ( Q_1,G(\BA_f) )'} \bigr) \stackrel{[1]}{\to} 
\]
to 
\[
R i_{Q_1}^! i^*R(i_G)_* \CV \longto
i_{Q_1}^* R(i_G)_* \CV \longto
i_{Q_1}^* R (h_{Q_1})_* h_{Q_1}^* i^*R(i_G)_* \CV \stackrel{[1]}{\longto} 
\] 
in $D^+(M^K ( Q_1,\FX )(\BC))$. \\[0.1cm]
(b)~If $Q_1$ is the intersection of $r$ distinct maximal proper sub-groups of $G$, then
\[
R i_{Q_1}^! i^*R(i_G)_* \CV
\cong Rp^K_* \bigl( \j_! \CV^{\partial BS}_{\tei \! e^K ( Q_1,G(\BA_f) )} \bigr)[-(r-1)]
\]
in $D^+(M^K ( Q_1,\FX )(\BC))$. \\[0.1cm] 
(c)~If $Q_1$ is not maximal proper, then the exact triangle
\[
R i_{Q_1}^! i^*R(i_G)_* \CV \stackrel{\alpha}{\longto} i_{Q_1}^* R(i_G)_* \CV 
\longto i_{Q_1}^* R (h_{Q_1})_* h_{Q_1}^* i^*R(i_G)_* \CV 
\stackrel{[1]}{\longto}                    
\]
is split: the morphism $\alpha$ is zero.
\end{Cor}

\begin{Proof}
Apply Proposition~\ref{6Ba},
Theorem~\ref{6MT2}, and proper base change \cite[Thm.~6.2]{I}.
\end{Proof}

\begin{Rem}
One word of caution concerning our (ab)use of the notation $i_{Q_1}$
is in order. By definition, the symbol denotes the immersion
\[
M^K (Q_1,\FX) \longinto M^K  (G,\FX)^* \; .
\]
Since $Q_1$ is assumed to be proper, $M^K (Q_1,\FX)$ is actually contained
in the boundary $\partial M^K  (G,\FX)^*$, and we agreed to use the symbol $i_{Q_1}$
to denote the immersion into $\partial M^K  (G,\FX)^*$ as well.
The two notations being thus related by the ``formula'' $i_{Q_1} = i \circ i_{Q_1}$,
the composition of inverse images $i_{Q_1}^* i^*$ can be simplified to $i_{Q_1}^*$,
and this is what we did in the second term of the lower row of Corollary~\ref{6D}~(a).
By contrast, the first term
\[
R i_{Q_1}^! i^*R(i_G)_* \CV
\]
cannot be simplified: first, take the \emph{usual} inverse image under $i$ of the object
$R(i_G)_* \CV$; then, take the \emph{exceptional} inverse image of the result under 
\[
M^K (Q_1,\FX) \longinto \partial M^K  (G,\FX)^* \; .
\]
\end{Rem}
\forget{
\begin{Cor} \label{6E}
Keep the notations of Corollary~\ref{6D}.
Assume that $K$ is neat,
and that the cohomology objects of $\CV \in D^+(M^K  (G,\FX) (\BC))$ are locally constant.
If $Q_1$ is (not maximal proper, and equal to)
the intersection of $r$ distinct maximal proper sub-groups of $G$, 
with $r \ge 2$, then there is an isomorphism from $i_{Q_1}^* R (h_{Q_1})_* h_{Q_1}^* i^* R(i_G)_* \CV$ to
\[ 
Rp^K_* \bigl( \CV^{\partial BS}_{\tei \! e^K ( Q_1,P(\BA_f) )'} \bigr) \oplus 
Rp^K_* \bigl( \j_! \CV^{\partial BS}_{\tei \! e^K ( Q_1,G(\BA_f) )} \bigr)[-(r-2)]
\]
in $D^+(M^K ( Q_1,\FX )(\BC))$, identifying the adjunction
\[
Rp^K_* \bigl( \CV^{\partial BS}_{\tei \! e^K ( Q_1,P(\BA_f) )'} \bigr) \stackrel{\ref{6D}}{\cong}
i_{Q_1}^* R(i_G)_* \CV \longto i_{Q_1}^* R (h_{Q_1})_* h_{Q_1}^* i^*R(i_G)_* \CV
\]
with the inclusion of the first component, and the boundary
\[
i_{Q_1}^* R (h_{Q_1})_* h_{Q_1}^* R(i_G)_* \CV 
\to R i_{Q_1}^! i^*R(i_G)_* \CV [1] \stackrel{\ref{6D}}{\cong}
Rp^K_* \bigl( \j_! \CV^{\partial BS}_{\tei \! e^K ( Q_1,G(\BA_f) )} \bigr)[-(r-2)]
\]
with the projection onto the second component.
\end{Cor}
}

To conclude this section, fix an element
$\uQ = (Q_1 \prec Q_2 \prec \ldots \prec Q_r)$ of $\CC_{(G,\FX)}$, such that
$Q_1$ is maximal proper, and $r \ge 2$. 
As before, 
\[
h_{Q_1}: \partial M^K  (G,\FX)^* - \overline{M^K \bigl( Q_1,\FX \bigr)}
\longinto \partial M^K  (G,\FX)^* 
\]
denotes
the open immersion of the complement of $\overline{M^K \bigl( Q_1,\FX \bigr)}$. \\

Let $\CV \in D^+(M^K  (G,\FX) (\BC))$, and assume that its
cohomology objects are locally constant.
Assume also that $K$ is neat.
Denote by $\CV^{\partial BS}$ the restriction to $\partial M^K  (G,\FX) (\BC)^{BS}$
of $R(k_G)_* \CV$. On the one hand, Corollary~\ref{6D}~(a) provides us with an isomorphism
\[
Rp^K_* \bigl( \CV^{\partial BS}_{\tei \! \partial e^K ( Q_1,G(\BA_f) )'} \bigr)
\isoto i_{Q_1}^* R (h_{Q_1})_* h_{Q_1}^* i^*R(i_G)_* \CV \; .
\]
On the other hand, the component $\uQ$ of the canonical isomorphism from Corollary~\ref{6C}
reads
\[
Rp^K_* \bigl( \CV^{\partial BS}_{\tei \! e^K ( \bigcap \uQ ,G(\BA_f) )'} \bigr)
\isoto i_{Q_1}^* R(i_{Q_2})_* i_{Q_2}^* \ldots R(i_{Q_r})_* i_{Q_r}^* i^*R(i_G)_* \CV \; .
\] 

\begin{Prop} \label{6F}
Together with the morphisms induced by the restriction
\[
\bigl( \CV^{\partial BS}_{\tei \! \partial e^K ( Q_1,G(\BA_f) )'} \bigr)_{\tei \! e^K ( \bigcap \uQ ,G(\BA_f) )'}
\isoto \CV^{\partial BS}_{\tei \! e^K ( \bigcap \uQ ,G(\BA_f) )'}
\]
and the adjunction
\[
R (h_{Q_1})_* h_{Q_1}^* \longto R(i_{Q_2})_* i_{Q_2}^* \ldots R(i_{Q_r})_* i_{Q_r}^* \; ,
\]
respectively, the above isomorphisms fit into a commutative diagram
\[
\vcenter{\xymatrix@R-10pt{
Rp^K_* \bigl( \CV^{\partial BS}_{\tei \! \partial e^K ( Q_1,G(\BA_f) )'} \bigr)
   \ar[r] \ar[d]^-{\ref{6D}~(a)}_-{\cong} & 
Rp^K_* \bigl( \CV^{\partial BS}_{\tei \! e^K ( \bigcap \uQ ,G(\BA_f) )'} \bigr) 
                                                        \ar[d]_-{\ref{6C}}^-{\cong} \\
i_{Q_1}^* R (h_{Q_1})_* h_{Q_1}^* i^*R(i_G)_* \CV 
   \ar[r] &  
i_{Q_1}^* R(i_{Q_2})_* i_{Q_2}^* \ldots R(i_{Q_r})_* i_{Q_r}^* i^*R(i_G)_* \CV
\\}} \; .
\] 
\end{Prop}

\begin{Proof}
Apply $R p_*^K$ to the commutative diagram from Proposition~\ref{5B}.
\end{Proof}

\begin{Rem} \label{6G}
Breaking down the construction of the diagram from Proposition~\ref{5B},
one sees that the horizontal morphisms
\[
Rp^K_* \bigl( \CV^{\partial BS}_{\tei \! \partial e^K ( Q_1,G(\BA_f) )'} \bigr)
\longto Rp^K_* \bigl( \CV^{\partial BS}_{\tei \! e^K ( \bigcap \uQ ,G(\BA_f) )'} \bigr)
\]
and
\[
i_{Q_1}^* R (h_{Q_1})_* h_{Q_1}^* i^*R(i_G)_* \CV 
\longto i_{Q_1}^* R(i_{Q_2})_* i_{Q_2}^* \ldots R(i_{Q_r})_* i_{Q_r}^* i^*R(i_G)_* \CV
\]
in Proposition~\ref{6F}
are equal to the adjunctions associated to the inclusions
\[
e^K \bigl( \bigcap \uQ ,G(\BA_f) \bigr)' \subset \partial e^K \bigl( Q_1,G(\BA_f) \bigr)'
\]
and
\[
M^K \bigl( Q_s,\FX \bigr) \subset 
\partial M^K  (G,\FX)^* - \overline{M^K \bigl( Q_1,\FX \bigr)} \; , \; s = 2,\ldots,r \; ,
\]
respectively.
We leave the verification to the reader.
\end{Rem}


\bigskip

%
%

\section{Reformulation in terms of group cohomology. I}
\label{7}



The purpose of this section is to translate Theorem~\ref{6MT1}
(through Corollary~\ref{6C})
into group cohomology (Main Theorem~\ref{7MT} and its Variants~\ref{7Var1}
and \ref{7Var2}). 
We fix pure Shimura data  $(P,\FX) = (G,\FX)$ satisfying hypothesis 
$(+)$, and an open compact sub-group $K$ of $G(\BA_f)$, which we assume
to be neat. We also fix 
a proper admissible parabolic sub-group $Q_1$ of $G$, and $g \in G (\BA_f)$. \\

By definition, the space of $\BC$-valued points of the stratum $M^K (Q_1,\FX)$ 
of $M^K (G,\FX)^*$ equals
\[
Q_1 (\BQ) \backslash \Bigl( \coprod \FX_1 / W_1 \times G (\BA_f) / K \Bigr) \; ,
\]
where $(P_1,\FX_1)$ runs over the finitely many rational boundary components 
of $(G,\FX)$ associated to $Q_1$ ($W_1:=$ the unipotent radical of $P_1$). 
The formulae we shall obtain in Main Theorem~\ref{7MT} and Variant~\ref{7Var1}
concern the (open and closed) sub-space 
$pr_{gK}(\coprod \FX_1/W_1)$ of
$M^K (Q_1,\FX)(\BC)$ defined as the image 
of $\coprod \FX_1/W_1 \times \{gK\}$ under the projection 
\[
\FX^* \times G (\BA_f) / K 
\longonto G (\BQ) \backslash \bigl( \FX^* \times G (\BA_f) / K \bigr) = M^K (G,\FX)^* (\BC) \; .
\]
As a quotient of $\coprod \FX_1 / W_1 = \coprod \FX_1 / W_1 \times \{gK\}$,
\[
pr_{gK} \bigl( \coprod \FX_1/W_1 \bigr) = H_1 \backslash \bigl( \coprod \FX_1/W_1 \bigr) \; ,
\]
where $H_1:= H_1(gK):= Q_1 (\BQ) \cap gKg^{-1}$. Let us denote by $i_{Q_1,g}$
the immersion of $pr_{gK}(\coprod \FX_1/W_1)$ into $\partial M^K  (G,\FX)^*(\BC)$.

\begin{Def}[{\cite[Def.~4.9]{W1}}] \label{7A}
Denote by $\pi_{Q_1}: Q_1 \onto \bar{Q}_1$ the ca\-no\-ni\-cal epimorphism
of $Q_1$ to its maximal reductive quotient.
Define the closed connected normal sub-group $C_1$ of $Q_1$ as
\[ 
C_1 := \{ q \in Q_1 \; , \; \pi_{Q_1}(q) \in \Cent_{\bar{Q}_1}(\pi_{Q_1}(P_1)) \}^0 \; .
\]
\end{Def}

Define
\[
H_C:=H_C(gK):= C_1 (\BQ) \cap gKg^{-1} \; .
\]
Let us start to set up the cohomological data necessary for the statement of 
Main Theorem~\ref{7MT}.

\begin{Def} \label{7C}
(a)~Define
\[
R \Gamma (H_C, \argdot): D^+ (\Rep H_1) \longto D^+ \bigl( \Rep (H_1/H_C) \bigr)
\]
as the derived functor of the $H_C$-invariants $\Gamma (H_C, \argdot) = (\argdot)^{H_C}$
on $\Rep H_1 \,$. \\[0.1cm]
(b)~Let $Q$ be a parabolic sub-group of $G$ contained in $Q_1 \,$. Define
\[
{}_Q \! \argdot: \Rep H_1 \longto \Rep H_1
\]
as the functor associating to $\BV$ the representation 
\[
{}_Q \! \BV := \App \bigl( Q_1(\BQ) / Q(\BQ) , \BV \bigr) 
\]
on (set-theoretical) applications $Q_1(\BQ) / Q(\BQ) \to \BV$, where 
by definition $h \in H_1$ acts as 
\[
f \longmapsto 
           \Bigl( hf: q_1 Q(\BQ) \longmapsto h \bigl( f(h^{-1} q_1 Q(\BQ)) \bigr) \Bigr) \; , 
\]
for all $f: Q_1(\BQ) / Q(\BQ) \to \BV$.  
\end{Def}

In part~(b), we leave the definition of the effect 
on morphisms to the reader; a similar procedure will be applied 
to some of the functorial
constructions that are to follow. There is a variant of the functor
${}_{Q_1} \! \argdot$ on the level of representations of a certain
quotient of $H_1 \,$: since the action of 
the subgroup $W_1(\BQ)$ of $Q_1(\BQ)$ on the set $Q_1(\BQ) / Q(\BQ)$
is trivial, setting
\[
H_W := H_W(gK) := W_1(\BQ) \cap H_1 (= W_1(\BQ) \cap H_C) \; ,
\]
there is a functor
\[
{}_Q \! \argdot: \Rep (H_1/H_W) \longto \Rep (H_1/H_W)
\]
defined by the same formula. \\

Note that ${}_{Q_1} \! \argdot$ is the identity. There are
natural transformations ${}_{\tQ} \! \argdot \to {}_Q \! \argdot$ for each pair $(Q,\tQ)$
of parabolics of $Q_1$ such that $Q \subset \tQ$, and these transformations
satisfy the cocycle relation. All functors ${}_Q \! \argdot$ are exact, and therefore
derive trivially to give triangulated endo-functors of $D^+ (\Rep H_1)$ and of
$D^+ (\Rep (H_1/H_W))$, respectively, still denoted by the same symbols.

\begin{Lem} \label{7inj}
Let $Q$ be a parabolic sub-group of $G$ contained in $Q_1 \,$.
The endo-functor ${}_Q \! \argdot$ maps injective objects
to injective objects.
\end{Lem}

\begin{Proof}
The functor ${}_Q \! \argdot$ admits a left adjoint 
$\argdot \otimes_\BZ \BZ[Q_1(\BQ) / Q(\BQ)]$,
and this left adjoint is exact. 
\end{Proof}

\begin{Prop} \label{7D}
Let $Q$ be a parabolic sub-group of $G$ contained in $Q_1 \,$, and
$\Omega \subset Q_1(\BQ)$ a set of representatives of 
$H_C \backslash Q_1(\BQ) / Q(\BQ)$. \\[0.1cm]
(a)~The functor 
\[
{}_Q \! \BV \mapsto \oplus_{q_1 \in \Omega} \BV \; , \;
f \longmapsto \bigl( f(q_1 Q(\BQ)) \bigr)_{q_1 \in \Omega}
\]
induces an isomorphism
\[
\Gamma (H_C, {}_Q \! \argdot) 
   \isoto \bigoplus_{q_1 \in \Omega} \Gamma \bigl(H_C \cap q_1 Q(\BQ) q_1^{-1}, \argdot \bigr)
\]
of functors on $\Rep H_1 \,$. \\[0.1cm]
(b)~The isomorphism from (a) extends to an isomorphism
\[
R \Gamma (H_C, {}_Q \! \argdot) \isoto 
\bigoplus_{q_1 \in \Omega} R \Gamma \bigl(H_C \cap q_1 Q(\BQ) q_1^{-1}, \argdot \bigr)
\]
of functors on $D^+ (\Rep H_1)$, where for each $q_1 \in \Omega$, 
\[
R \Gamma \bigl(H_C \cap q_1 Q(\BQ) q_1^{-1}, \argdot \bigr)
\]
denotes the derived functor of 
$\Gamma (H_C \cap q_1 Q(\BQ) q_1^{-1}, \argdot )$. \\[0.1cm]
(c)~Let $\tQ$ be a parabolic sub-group of $Q_1$ containing $Q$, and
$\tilde{\Omega} \subset Q_1(\BQ)$ a set of representatives of 
$H_C \backslash Q_1(\BQ) / \tQ(\BQ)$. Then the isomorphisms from (b),
for $Q$ and for $\tQ$, fit into a commutative diagram
\[
\vcenter{\xymatrix@R-10pt{
R \Gamma (H_C, {}_{\tQ} \! \argdot) \ar[r]^-\cong \ar[d] & 
\bigoplus_{\tilde{q_1} \in \tilde{\Omega}} 
 R \Gamma \bigl(H_C \cap \tilde{q_1} \tQ(\BQ) \tilde{q_1}^{-1}, \argdot \bigr) 
           \ar[d]^-{\Res^{\tilde{\Omega}}_\Omega} \\
R \Gamma (H_C, {}_Q \! \argdot) \ar[r]^-\cong & 
\bigoplus_{q_1 \in \Omega} R \Gamma \bigl(H_C \cap q_1 Q(\BQ) q_1^{-1}, \argdot \bigr)
\\}} \quad ,
\] 
where the natural transformation 
$R \Gamma (H_C, {}_{\tQ} \! \argdot) \to R \Gamma (H_C, {}_Q \! \argdot)$ is
induced by ${}_{\tQ} \! \argdot \to {}_Q \! \argdot$. 
The natural transformation
\[
\Res^{\tilde{\Omega}}_\Omega: \bigoplus_{\tilde{q_1} \in \tilde{\Omega}} 
 R \Gamma \bigl(H_C \cap \tilde{q_1} \tQ(\BQ) \tilde{q_1}^{-1}, \argdot \bigr)
\longto 
 \bigoplus_{q_1 \in \Omega} R \Gamma \bigl(H_C \cap q_1 Q(\BQ) q_1^{-1}, \argdot \bigr) 
 \]
is defined as follows: for each $\tilde{q_1} \in \tilde{\Omega}$, put 
\[
\Omega_{\tilde{q_1}} := \Omega \cap H_C \tilde{q_1} \tQ(\BQ) \subset \Omega \; .
\]
On $R \Gamma (H_C \cap \tilde{q_1} \tQ(\BQ) \tilde{q_1}^{-1}, \argdot )$, the
transformation $\Res^{\tilde{\Omega}}_\Omega$ is then defined as the direct sum over
all $q_1 \in \Omega_{\tilde{q_1}}$ of the restriction
\[
R \Gamma \bigl(H_C \cap \tilde{q_1} \tQ(\BQ) \tilde{q_1}^{-1}, \argdot \bigr)
\longto 
R \Gamma \bigl(H_C \cap q_1 Q(\BQ) q_1^{-1}, \argdot \bigr) 
\]  
from $H_C \cap \tilde{q_1} \tQ(\BQ) \tilde{q_1}^{-1}$ to the sub-group
$H_C \cap q_1 Q(\BQ) q_1^{-1}$.
\end{Prop}

\emph{A priori}, the targets of the isomorphisms from Proposition~\ref{7D}~(a) and (b)
are functors with values in the category $\Ab$ of Abelian groups, and in $D^+(\Ab)$. The
additional action of the quotient $H_1/H_C$ is obtained \emph{a posteriori} from
the sources of the isomorphisms. We leave it to the reader to make explicit that action.

\medskip

\begin{Proofof}{Proposition~\ref{7D}}
(a): left to the reader.

\noindent (b): the functor ${}_Q \! \argdot$ is exact, and according to Lemma~\ref{7inj},
it respects injectives. Therefore,
any injective resolution $\BI$ of $\BV$ induces an injective resolution 
${}_Q \BI$ of ${}_Q \! \BV$. Altogether,
\[
R \Gamma (H_C, {}_Q \! \BV) = R \Gamma (H_C, {}_Q \BI) \stackrel{(a)}{\isoto} 
    \bigoplus_{q_1 \in \Omega} R \Gamma \bigl(H_C \cap q_1 Q(\BQ) q_1^{-1}, \BI \bigr) \; ,
\]
and the latter object equals 
$\oplus_{q_1 \in \Omega} R \Gamma (H_C \cap q_1 Q(\BQ) q_1^{-1}, \BV )$,
as the restriction to a sub-group respects injectives.

(c): left to the reader.
\end{Proofof} 

\begin{Rem} \label{7Db}
(a)~The proof of Proposition~\ref{7D}~(b) is a variant of the proof of Shapiro's 
Lemma. \\[0.1cm]
(b)~There is a variant of Proposition~\ref{7D} for the functors
\[
\Gamma (H_C/H_W, {}_Q \! \argdot)  
\] 
and 
\[
R \Gamma (H_C/H_W, {}_Q \! \argdot)  
\]
on $\Rep (H_1/H_W)$ and on $D^+(\Rep (H_1/H_W))$, respectively,  
involving the same set of representatives $\Omega$
(but note that the canonical surjection
\[
H_C \backslash Q_1(\BQ) / Q(\BQ) \longto 
(H_C/H_W) \backslash \bar{Q}_1(\BQ) / ((Q/W_1)(\BQ))
\]
is a bijection),
and the groups $(H_C \cap q_1 Q(\BQ) q_1^{-1})/H_W$
instead of $H_C \cap q_1 Q(\BQ) q_1^{-1}$. 
We leave the statement and the proof
to the reader. \\[0.1cm]
(c)~\emph{Via} $q_1 \mapsto C_j \cap q_1 Q q_1^{-1}$, 
the double quotient 
\[
H_C \backslash Q_j(\BQ) / Q(\BQ)
\] 
is in bijection with the set of orbits
under the action of $H_C$ on the parabolic subgroups of $C_j$ that are
$Q_j(\BQ)$-conjugate to $C_j \cap Q$ \cite[Rem.~6.22~(a)]{W1}. Note that
the set of such parabolics of $C_j$ may consist of several $C_j(\BQ)$-conjugation classes 
(depending on the index of $C_j(\BQ) Q(\BQ)$ in $Q_j(\BQ)$).
\end{Rem}

The functors 
\[
{}_Q \! \argdot: \Rep H_1 \longto \Rep H_1
\]
and 
\[
{}_Q \! \argdot: \Rep (H_1/H_W) \longto \Rep (H_1/H_W)
\]
are related by the following result; Lemma~\ref{7inj}
provides the essential ingredient of its proof, which we leave to the reader.

\begin{Prop} \label{7Da}
Let $Q$ be a parabolic sub-group of $G$ contained in $Q_1 \,$. \\[0.1cm]
(a)~We have
\[
\Gamma ( H_W, {}_Q \! \argdot ) = {}_Q \bigl( \Gamma ( H_W, \argdot ) \bigr)
\]
as functors from $\Rep H_1$ to $\Rep (H_1/H_W)$. \\[0.1cm]
(b)~The equality from (a) extends to a canonical isomorphism
\[
R \Gamma ( H_W, {}_Q \! \argdot ) 
\isoto {}_Q \bigl( R \Gamma ( H_W, \argdot ) \bigr)
\]
of functors from $D^+(\Rep H_1)$ to $D^+ (\Rep (H_1/H_W))$
($R \Gamma ( H_W, \argdot ) :=$ the derived functor of
$\Gamma ( H_W, \argdot )$).   
\end{Prop}

\begin{Def} \label{7E}
Define a category $\CR_{(Q_1,\FX) \tei g}^K$ as follows: objects are of the form
\[
\bigl( \BW_{\uQ} , a_{\uQ \subset \utQ} \bigr)_{\uQ, \utQ \in \CC_{(G,\FX) \tei Q_1}} \; ,
\]
where $\BW_{\uQ} \in D^+ ( \Rep (H_1/H_C) )$,
for each $\uQ \in \CC_{(G,\FX) \tei Q_1}$, and 
\[
a_{\uQ \subset \utQ}: \BW_{\utQ} \longto \BW_{\uQ}
\]
is a morphism in $D^+ ( \Rep (H_1/H_C) )$, for each relation $\uQ \subset \utQ$
in $\CC_{(G,\FX) \tei Q_1}$. 
Morphisms between $(\BW_{\uQ} , a_{\uQ \subset \utQ})_{\uQ,\utQ}$
and $(\BX_{\uQ} , b_{\uQ \subset \utQ})_{\uQ, \utQ}$ are of the form
$( \alpha_{\uQ} )_{\uQ}$,
where $\alpha_{\uQ}: \BW_{\uQ} \to \BX_{\uQ}$ is a morphism in 
$D^+ ( \Rep (H_1/H_C) )$, such that for each relation $\uQ \subset \utQ$, one has 
\[
\alpha_{\uQ} \circ a_{\uQ \subset \utQ} = b_{\uQ \subset \utQ} \circ \alpha_{\utQ}:
   \BW_{\utQ} \longto \BX_{\uQ} \; .
\]  
\end{Def}

\begin{Cons} \label{7Con}
(a)~Let $\BV \in D^+ (\Rep H_1)$. Define an object
$\coh(\BV)$ of $\CR_{(Q_1,\FX) \tei g}^K$ as follows.
\[
\coh(\BV) = \bigl( \BW_{\uQ} , a_{\uQ \subset \utQ} \bigr)_{\uQ, \utQ} \; ,
\]
where for $\uQ \in \CC_{(G,\FX) \tei Q_1}$, 
\[
\BW_{\uQ} := R \Gamma \bigl( H_C, {}_{\bigcap \uQ} \! \BV \bigr)
\in D^+ \bigl( \Rep (H_1/H_C) \bigr) \; ,
\]
and for $\uQ \subset \utQ$,
\[
a_{\uQ \subset \utQ}: R \Gamma \bigl( H_C, {}_{\bigcap \utQ} \! \BV \bigr) 
\longto R \Gamma \bigl( H_C, {}_{\bigcap \uQ} \! \BV \bigr)
\]
is defined as being induced by the natural transformation 
${}_{\bigcap \utQ} \! \argdot \to {}_{\bigcap \uQ} \! \argdot$. \\[0.1cm]
(b)~Let $\BV \in D^+ (\Rep (H_1/H_W))$. Define an object
\[
\overline{\coh}(\BV) = \bigl( \BW_{\uQ} , a_{\uQ \subset \utQ} \bigr)_{\uQ, \utQ} 
\]
of $\CR_{(Q_1,\FX) \tei g}^K$, 
where for $\uQ \in \CC_{(G,\FX) \tei Q_1}$, 
\[
\BW_{\uQ} := R \Gamma \bigl( H_C/H_W, {}_{\bigcap \uQ} \! \BV \bigr)
\in D^+ \bigl( \Rep (H_1/H_C) \bigr) 
\]
($R \Gamma ( H_C/H_W, \argdot ) :=$ the derived functor of
\[
\Gamma ( H_C/H_W, \argdot ) : \Rep (H_1/H_W) \longto \Rep (H_1/H_C) \; ).
\]
The morphisms $a_{\uQ \subset \utQ}$, $\uQ \subset \utQ$,
are defined as in part~(a).
\end{Cons}

We leave it to the reader to complete Construction~\ref{7Con} to give functors
\[
\coh: D^+ (\Rep H_1) \longto \CR_{(Q_1,\FX) \tei g}^K 
\]
and
\[
\overline{\coh}: D^+ (\Rep H_1/H_W) \longto \CR_{(Q_1,\FX) \tei g}^K \; .
\]

\begin{Prop} \label{7Cona}
There is a canonical isomorphism
\[
\coh \isoto \overline{\coh} \circ R \Gamma ( H_W, \argdot ) 
\]
of functors from $D^+ (\Rep H_1)$ to $\CR_{(Q_1,\FX) \tei g}^K$.
\end{Prop}

\begin{Proof}
This follows from the definitions, and from Proposition~\ref{7Da}~(b).
\end{Proof}

\begin{Def} \label{7F}
(a)~Define a category $\CC_{(Q_1,\FX) \tei g}^K$ as follows: objects are of the form
\[
\bigl( \CK_{\uQ} , a_{\uQ \subset \utQ} \bigr)_{\uQ, \utQ \in \CC_{(G,\FX) \tei Q_1}} \; ,
\]
where $\CK_{\uQ} \in D^+(pr_{gK}(\coprod \FX_1/W_1))$,
for each $\uQ \in \CC_{(G,\FX) \tei Q_1} \,$, and 
\[
a_{\uQ \subset \utQ}: R (i_{Q_1,g})_* \CK_{\utQ} \longto R (i_{Q_1,g})_* \CK_{\uQ}
\]
is a morphism in $D^+(M^K  (G,\FX)^* (\BC))$, for each relation $\uQ \subset \utQ$
in $\CC_{(G,\FX) \tei Q_1}$.
Morphisms are defined as in $\CC_{(G,\FX)}^K$ (Definition~\ref{3E}~(a)). \\[0.1cm]
(b)~Define
\[
i^*_{Q_1,g}: \CC_{(G,\FX)}^{K,*} \longto \CC_{(Q_1,\FX) \tei g}^K
\]
as $i^*_{Q_1}$ (Definition~\ref{3E}~(c)), followed by
the restriction from $M^K  (Q_1,\FX)(\BC)$ to $pr_{gK}(\coprod \FX_1/W_1)$.
\end{Def}

Since $K$ is assumed neat, the action of $G(\BQ)$ on $\FX \times G (\BA_f) / K$
is free \cite[Lemma~1.3]{BW}.

\begin{Def} \label{7B}
Define the \emph{canonical construction}
\[
\mu_{K}: \Rep \bigl( G(\BQ) \bigr) \longto \Loc \bigl( M^K (G,\FX) (\BC) \bigr)
\] 
(of level $K$)
from the category $\Rep(G(\BQ))$ of representations of $G(\BQ)$ in Abelian groups
to the category $\Loc(M^K (G,\FX) (\BC))$ of local systems
on the space $M^K (G,\FX) (\BC)$ as the functor associating to each representation $\BV$
the sheaf of continuous sections of
\[
G (\BQ) \backslash \bigl( \BV \times \FX \times G (\BA_f) / K \bigr) 
\longonto G (\BQ) \backslash \bigl( \FX \times G (\BA_f) / K \bigr)
= M^K (G,\FX) (\BC) 
\]
(with respect to the discrete topology on $\BV$).
\end{Def}

Again since $K$ is neat, the action of the normal sub-group $H_C$
of $H_1$ on $\coprod \FX_1 / W_1$ is trivial, and the induced action 
on $\coprod \FX_1 / W_1$
of the quotient $H_1/H_C$ is free \cite[Prop.~7.22~(b)]{W1}. 
In analogy to the construction from Definition~\ref{7B}, define 
\[
\bar{\mu}: \Rep(H_1/H_C) \longto \Loc \bigl( pr_{gK} \bigl(\coprod \FX_1/W_1 \bigr) \bigr) 
\] 
as the functor associating to each representation $\BV$
the sheaf of continuous sections of
\[
(H_1/H_C) \backslash \bigl( \BV \times \coprod \FX_1/W_1  \bigr) 
\longonto (H_1/H_C) \backslash \bigl( \coprod \FX_1/W_1  \bigr)
= pr_{gK} \bigl(\coprod \FX_1/W_1 \bigr) \; .
\]
As before, this functor induces
\[
\bar{\mu}: D^+ \bigl( \Rep(H_1/H_C) \bigr) \longto 
D^+ \bigl( pr_{gK} \bigl(\coprod \FX_1/W_1 \bigr) \bigr) \; .
\]
We define
\[
\nu_g : \CR_{(Q_1,\FX) \tei g}^K \longto \CC_{(Q_1,\FX) \tei g}^K
\]
as the componentwise extension of $\bar{\mu}$ to $\CR_{(Q_1,\FX) \tei g}^K \,$. \\

Recall that $i_G$ denotes the open immersion of $M^K  (G,\FX)$ into $M^K  (G,\FX)^*$, 
and $i$ the complementary immersion of $\partial M^K  (G,\FX)^*$. 

\begin{MThm} \label{7MT}
(a)~There is a natural commutative diagram
\[
\vcenter{\xymatrix@R-10pt{
D^+ \bigl( \Rep \bigl(G(\BQ) \bigr) \bigr) \ar[rr]^-{\mu_K} \ar[dd]_-{\Res^{G(\BQ)}_{H_1}} 
         &&
                  D^+ \Bigl( \Loc \bigl( M^K (G,\FX) (\BC) \bigr) \Bigr) \ar[d] \\
         && D^+ \bigl( M^K  (G,\FX) (\BC) \bigr) \ar[d]^-{i^* R(i_G)_*} \\
         D^+ ( \Rep H_1 ) \ar[dd]_-{\coh} 
         &&  D^+ \bigl( \partial M^K  (G,\FX)^* (\BC) \bigr) \ar[d]^-{\deg} \\
         && \CC_{(G,\FX)}^{K,*} \ar[d]^-{i^*_{Q_1,g}} \\
\CR_{(Q_1,\FX) \tei g}^K \ar[rr]^-{\nu_g} && \CC_{(Q_1,\FX) \tei g}^K
\\}}
\]
($\Res^{G(\BQ)}_{H_1} :=$ the restriction from $G(\BQ)$ to $H_1$). \\[0.1cm]
(b)~The functor 
\[
\deg \circ i^* R(i_G)_* \circ \mu_K: 
                    D^+ \bigl( \Rep \bigl(G(\BQ) \bigr) \bigr) \longto \CC_{(G,\FX)}^{K,*}
\]
takes values in the sub-category of objects 
$( \CK_{\uQ} , a_{\uQ \subset \utQ})_{\uQ, \uQ \subset \utQ}$
of $\CC_{(G,\FX)}^{K,*}$, all of whose components
\[
\CK_{\uQ} \in D^+ \bigl( M^K (b(\uQ),\FX)(\BC) \bigr)
\] 
can be represented by complexes of local systems.
\end{MThm}

\begin{Rem} \label{7Rem}
(a)~Main Theorem~\ref{7MT} 
yields comparison statements indexed by $\uQ \in \CC_{(G,\FX) \tei Q_1} \,$, 
and under these comparisons, it identifies the transition 
$a_{\uQ \subset \utQ}$
for each relation $\uQ \subset \utQ$ in $\uQ \in \CC_{(G,\FX) \tei Q_1} \,$. 
The simplest case concerns $\uQ = (Q_1)$ (the singleton consisting of 
the parabolic $Q_1$). Given the definitions of the functors $\coh$
(Construction~\ref{7Con}(a)) and $\deg$ (Construction~\ref{6Con}),
Main Theorem~\ref{7MT}~(a) implies that
\[
i_{Q_1,g}^* \circ i^* R(i_G)_* \circ \mu_K \cong 
\bar{\mu} \circ R \Gamma (H_C,\argdot) \circ \Res^{G(\BQ)}_{H_1} 
\]
as functors from $D^+ ( \Rep (G(\BQ) ) )$ to $D^+(pr_{gK}(\coprod \FX_1/W_1))$.
This result is well known: see \emph{e.g.} 
\cite[(6.2)]{LR}. \\[0.1cm]
(b)~The formula $\coh \cong \overline{\coh} \circ R \Gamma ( H_W, \argdot )$ from
Proposition~\ref{7Cona} provides an alternative description of the
left half of the commutative diagram from Main Theorem~\ref{7MT}~(a).
This will be exploited in the context of algebraic representations
of $G$ (see Variants~\ref{7Var1} and \ref{7Var2}).
\end{Rem}

The proof of Main Theorem~\ref{7MT} involves 
the Borel--Serre compactification $M^K (G,\FX) (\BC)^{BS}$. 
We need an auxiliary construction
on the level of $M^K (G,\FX) (\BC)^{BS}$.

\begin{Aux} \label{7Aux}
(a)~The functor $\mu_K$ admits an extension
\[
\mu_{K}^{BS}: \Rep \bigl( G(\BQ) \bigr) \longto \Loc \bigl( M^K (G,\FX) (\BC)^{BS} \bigr) \; ,
\]
equal to $R(k_G)_* \circ \mu_K$ (\cite[Prop.~5.4~(c)]{W1} and Corollary~\ref{1C}~(c)).
Note: the action of $G (\BQ)$ on $\FX^{BS} \times G (\BA_f) / K$ being free
\cite[9.5]{BS}, the functor $\mu_{K}^{BS}$ can be defined directly, 
by associating to each representation $\BV$
the sheaf of continuous sections of
\[
G (\BQ) \backslash \bigl( \BV \times \FX^{BS} \times G (\BA_f) / K \bigr) 
\longonto M^K (G,\FX) (\BC)^{BS} \; .
\]
It restricts to a functor
\[
\mu_{K}^{\partial BS}: \Rep \bigl( G(\BQ) \bigr) 
                    \longto \Loc \bigl( M^K (G,\FX) (\BC)^{\partial BS} \bigr) \; .
\]
By the very definition of $\CV \mapsto \CV^{\partial BS}$ (see Corollary~\ref{6C}),
we have 
\[
\mu_{K}^{\partial BS} = (\CV \mapsto \CV^{\partial BS}) \circ \mu_K \; .
\]
The functor $\mu_{K}^{\partial BS}$ being exact, it derives trivially, to give
\[
\mu_{K}^{\partial BS}: D^+ \bigl( \Rep \bigl(G(\BQ) \bigr) \bigr)
     \longto D^+ \Bigl( \Loc \bigl( M^K (G,\FX) (\BC)^{\partial BS} \bigr) \Bigr) \; .
\]
(b)~In analogy to Definition~\ref{7F}, define a category $\CC_{(Q_1,\FX) \tei g}^{K,BS}$ 
as follows: objects are of the form
\[
\bigl( \CK_{\uQ} , a_{\uQ \subset \utQ} \bigr)_{\uQ, \utQ \in \CC_{(G,\FX) \tei Q_1}} \; ,
\]
where
\[
\CK_{\uQ} \in 
 D^+ \Bigl( \bigl( p^K \bigr)^{-1} \bigl(pr_{gK} \bigl( \coprod \FX_1/W_1 \bigr) \bigr) 
    \cap e^K \bigl( \bigcap \uQ,G(\BA_f) \bigr)' \Bigr) \; ,
\]
for each $\uQ \in \CC_{(G,\FX) \tei Q_1} \,$, and 
\[
a_{\uQ \subset \utQ}: R (k_{\bigcap \utQ,g})_* \CK_{\utQ} 
                                        \longto R (k_{\bigcap \uQ,g})_* \CK_{\uQ}
\]
is a morphism in $D^+(M^K  (G,\FX) (\BC)^{BS})$, for each relation $\uQ \subset \utQ$
in $\CC_{(G,\FX) \tei Q_1}$ ($k_{\bigcap \uQ,g} :=$ the immersion of
\[
\bigl( p^K \bigr)^{-1} \bigl(pr_{gK} \bigl(\coprod \FX_1/W_1 \bigr) \bigr) 
              \cap e^K \bigl( \bigcap \uQ,G(\BA_f) \bigr)'
\]
into $M^K  (G,\FX) (\BC)^{BS}$).
Morphisms are defined as in $\CC_{(G,\FX)}^{K,BS}$. Also, 
define
\[
k^*_{Q_1,g}: \CC_{(G,\FX)}^{K,BS} \longto \CC_{(Q_1,\FX) \tei g}^{K,BS}
\]
to be the restriction. \\[0.1cm]
(c)~The continuous map
\[
p^K: M^K (G,\FX) (\BC)^{BS} \longto M^K (G,\FX)^* (\BC)
\]
is induced by the $G(\BQ)$-equivariant continuous map
\[
p \times \id_{G (\BA_f) / K} : 
      \FX^{BS} \times G (\BA_f) / K \longto \FX^* \times G (\BA_f) / K \; ,
\]
where $p: \FX^{BS} \to \FX^*$ is defined in 
\cite[Constr.~4.13]{W1} (see also \cite[Sect.~(3.7)--(3.11)]{Z}).
The sub-space $p^{-1} ( \coprod \FX_1/W_1 ) \times \{gK\}$ of
$\FX^{BS} \times G (\BA_f) / K$ is stable under $H_1 \subset G (\BQ)$.
According to \cite[Prop.~7.22~(a)]{W1},
the pre-image under $p^K$ of $pr_{gK}(\coprod \FX_1/W_1)$ equals
\[
H_1 \backslash \Bigl( p^{-1} \bigl( \coprod \FX_1/W_1 \bigr) \times \{gK\} \Bigr)
\subset G (\BQ) \backslash \bigl( \FX^{BS} \times G (\BA_f) / K \bigr) \; .
\]
The action of $G (\BQ)$ on
$\FX^{BS} \times G (\BA_f) / K$ being free, the same is true for the (induced)
action of $H_1$ on $p^{-1} ( \coprod \FX_1/W_1 ) \times \{gK\}$.
We thus get a functor
\[
\Rep ( H_1 ) \longto 
\Loc \Bigl( \bigl( p^K \bigr)^{-1} \bigl( pr_{gK} \bigl(\coprod \FX_1/W_1 \bigr) \bigr) \Bigr) \; ,
\]
by associating to each representation the sheaf of continuous sections,
in complete ana\-logy to $\nu_g$.
Define 
\[
\eta_g: D^+ ( \Rep H_1 ) \longto \CC_{(Q_1,\FX) \tei g}^{K,BS}
\]
to be equal to the composition of its (trivially) derived functor, 
and the restrictions to the 
\[
\bigl( p^K \bigr)^{-1} \bigl(pr_{gK} \bigl(\coprod \FX_1/W_1 \bigr) \bigr) 
    \cap e^K \bigl( \bigcap \uQ,G(\BA_f) \bigr)' \; ,
\]
for $\uQ \in \CC_{(G,\FX) \tei Q_1}$. The transition morphisms
$a_{\uQ \subset \utQ}$, for $\uQ \subset \utQ$, are given by adjunction.
\end{Aux}

\begin{Prop} \label{7G}
There is a natural commutative diagram  
\[
\vcenter{\xymatrix@R-10pt{
D^+ \bigl( \Rep \bigl(G(\BQ) \bigr) \bigr) 
                   \ar[rr]^-{\mu_K^{\partial BS}} \ar[ddd]_-{\Res^{G(\BQ)}_{H_1}} &&
            D^+ \Bigl( \Loc \bigl( \partial M^K  (G,\FX) (\BC)^{BS} \bigr) \Bigr)\ar[d] \\
         && D^+_{\Loc} \bigl( \partial M^K  (G,\FX) (\BC)^{BS} \bigr)  
                                  \ar[d]^-{\res} \\
         && \CC_{(G,\FX)}^{K,BS} \ar[d]^-{k^*_{Q_1,g}} \\
D^+ ( \Rep H_1 ) \ar[rr]^-{\eta_g} && \CC_{(Q_1,\FX) \tei g}^{K,BS} 
\\}}
\]
($D^+_{\Loc} ( \partial M^K  (G,\FX) (\BC)^{BS} ) :=$ the full (triangulated) sub-category
of objects of $D^+ ( \partial M^K  (G,\FX) (\BC)^{BS} )$, whose cohomology
objects are locally constant). 
\end{Prop}

\begin{Proof} 
This is clear from our definitions.
\end{Proof} 

The proof of Main Theorem~\ref{7MT}, requires a
decription of the sub-space $p^{-1} ( \coprod \FX_1/W_1 )$ of $\FX^{BS}$
occurring in Auxiliary Construction~\ref{7Aux}~(c).

\begin{Def}[{\cite[Def.~5.3]{W1}}] \label{7H}
Let $Q$ be a parabolic sub-group of $P$ satisfying $\adm(Q) = Q_1$
(in other words, $P_1 \subset Q \subset Q_1$). Define 
\[
e(Q)' := \coprod_{P_1 \subset R \subset Q} e(R) \subset \FX^{BS} \; ,
\]
where the union runs over all parabolic sub-groups $R$ contained in $Q \,$,
and containing $P_1$.
\end{Def}

\begin{Prop}[{\cite[Comp.~4.15]{W1}}] \label{7I}
Under the $G(\BQ)$-equivariant continuous map $p: \FX^{BS} \to \FX^*$, we have
\[
p^{-1} \bigl( \coprod \FX_1 / W_1 \bigr) = e(Q_1)' \subset \FX^{BS} \; .
\]
\end{Prop}

Recall \cite[Rem.~3.12]{W1} that the induced action of $G(\BQ)$ 
on the set of faces of $\FX^{BS}$
is given by $\gamma: e(R) \mapsto e(\gamma R \gamma^{-1})$. 
In particular, the sub-space $e(Q_1)'$ is stable under 
$Q_1(\BQ) \subset G(\BQ)$, 
and for a parabolic sub-group $Q$ of $Q_1$ satisfying $\adm(Q) = Q_1 \,$, 
the face $e(Q) \subset e(Q_1)'$ is mapped to $e(\gamma Q \gamma^{-1}) \subset e(Q_1)'$ 
whenever $\gamma \in Q_1(\BQ)$,
and similarly for $e(Q)' \subset e(Q_1)'$.
The restriction of the action of $Q_1(\BQ)$ on $e(Q_1)'$ to $H_1 \subset Q_1(\BQ)$ equals the action
on 
\[
p^{-1} \bigl( \coprod \FX_1 / W_1 \bigr) = p^{-1} \bigl( \coprod \FX_1 / W_1 \bigr) \times \{gK\}
\]
from Auxiliary Construction~\ref{7Aux}~(c). \\

In order to conclude the preparation of the proof of
Main Theorem~\ref{7MT},
we recall a number of principles concerning group cohomology, and its relation to cohomology of topological spaces. 

\begin{Scho} \label{7Scho}
Let $X$ and $Z$ be two topological spaces, together
with actions of a group $H'$ on the product $X \times Z$ and on $Z$, such that 
the projection onto the second factor $p^X: X \times Z \longonto Z$
is $H'$-equivariant. Let $H$ be a normal sub-group of $H'$. We make the following assumptions:
\begin{enumerate}
\item[(i)]~the induced action of $H'$ on the set of connected components $\pi_0(X \times Z)$
is diagonal with respect to the identification 
\[
\pi_0(X \times Z) = \pi_0(X) \times \pi_0(Z) \; .
\]
In other words: for any $x \in X$ and any $h \in H'$, there is a connected
component $X^0$ of $X$ such that 
\[
h(x,z) \in X^0 \times Z
\]
whenever $z \in Z$, 
\item[(ii)]~the
action of $H'$ on $Z$ factors through the quotient $H'/ H$. In other words:
the restriction to $H$ of the action on $X \times Z$
is diagonal, with the trivial action on $Z$, 
\item[(iii)]~the actions of $H'$ on $X \times Z$, and of $H'/ H$ on $Z$ are both free.
\end{enumerate}
We get exact functors
\[
\mu^{X \times Z}_{H'}: \Rep H' \longto \Loc \bigl( H' \backslash (X \times Z) \bigr) 
\] 
and 
\[
\mu^{X \times Z}_{H'}: D^+ ( \Rep H' ) 
                                 \longto D^+ \bigl( H' \backslash (X \times Z) \bigr) \; ,
\]
\[
\mu^Z_{H'/ H}: \Rep (H'/ H) \longto \Loc (H' \backslash Z) 
\] 
and 
\[
\mu^Z_{H'/ H}: D^+ \bigl( \Rep (H'/ H) \bigr)  \longto D^+ (H' \backslash Z) 
\]
as usual, by associating to each representation $\BV$ the sheaf of continuous sections of
\[
H' \backslash ( \BV \times X \times Z) \longonto H' \backslash ( X \times Z) 
\]
and of
\[
(H'/ H) \backslash ( \BV \times Z) \longonto H' \backslash Z \; ,
\]
respectively. Denote by
\[
p^{X,H'}: H' \backslash (X \times Z) \longonto H' \backslash Z
\]
the map induced by $p^X$.
We then have: 
\begin{enumerate}
\item[(1)]~there is a canonical isomorphism of functors
\[
\mu^Z_{H'/ H} \Gamma \bigl( H , \App( \pi_0(X) , \argdot ) \bigr) \isoto R^0p^{X,H'}_* \circ \mu^{X \times Z}_{H'} 
\] 
on $\Rep H'$ ($\App( \pi_0(X) , \argdot ):=$ applications whose source is
$\pi_0(X)$, together with the action
mapping $h \in H'$ to
\[
f \longmapsto h \circ f \circ h^{-1} \; )
\]
(use the very definition of $\mu^X$),
\item[(2)]~if all connected components of $X$ are contractible, then (1) extends to a canonical
isomorphism
\[
\mu^Z_{H'/ H} R \Gamma \bigl( H , \App( \pi_0(X) , \argdot ) \bigr) \isoto Rp^{X,H'}_* \circ \mu^{X \times Z}_{H'} 
\] 
of functors 
on $D^+ ( \Rep H' )$ (use Shapiro's Lemma \cite[Prop.~2]{AW} to reduce
to the case where $X$ is connected), 
\item[(3)]~if $X_1$ is a sub-space of $X$,
stable under the action of $H'$, and such that the action of $H'$  on $\pi_0(X_1 \times Z)$
remains diagonal, then the isomorphisms from (1) fit into a commutative
diagram
\[
\vcenter{\xymatrix@R-10pt{
\mu^Z_{H'/ H} \Gamma \bigl( H , \App( \pi_0(X) , \argdot ) \bigr)  
                         \ar[r]^-{\cong} \ar[d]_-{\mu^Z_{H'/ H} \Gamma(H,\Res^X_{X_1})} &  
         R^0p^{X,H'}_* \circ \mu^{X \times Z}_{H'} \ar[d]^-{R^0p^{X_1,H'}_* \Res^X_{X_1}}  \\
\mu^Z_{H'/ H} \Gamma \bigl( H , \App( \pi_0(X_1) , \argdot ) \bigr)  \ar[r]^-{\cong} &  
         R^0p^{X_1,H'}_* \circ \mu^{X_1 \times Z}_{H'}
\\}}
\]
($\Res^X_{X_1}:=$ the restrictions from $\pi_0(X)$ to $\pi_0(X_1)$, 
and from $H' \backslash (X \times Z)$ to $H' \backslash (X_1 \times Z)$, respectively). 
If furthermore all connected components of $X$ 
and of $X_1$ are contractible, then the isomorphisms from (2) fit into a commutative
diagram
\[
\vcenter{\xymatrix@R-10pt{
\mu^Z_{H'/ H} R \Gamma \bigl( H , \App( \pi_0(X) , \argdot ) \bigr)  
                         \ar[r]^-{\cong} \ar[d]_-{\mu^Z_{H'/ H} R \Gamma(H,\Res^X_{X_1})} &  
         Rp^{X,H'}_* \circ \mu^{X \times Z}_{H'} \ar[d]^-{R p^{X_1,H'}_* \Res^X_{X_1}}  \\
\mu^Z_{H'/ H} R \Gamma \bigl( H , \App( \pi_0(X_1) , \argdot ) \bigr)  \ar[r]^-{\cong} &  
         Rp^{X_1,H'}_* \circ \mu^{X_1 \times Z}_{H'}
\\}}
\quad .
\]
\end{enumerate}
\end{Scho}

\begin{Proofof}{Main Theorem~\ref{7MT}}
Part~(b) is a formal consequence of (a) (for all $Q_1$ and all $g$).

The first step in the proof of part~(a) is to apply Main Theorem~\ref{6MT1} 
through Corollary~\ref{6C}.
According to the latter, we need to establish 
the existence of a natural commutative diagram
\[
\vcenter{\xymatrix@R-10pt{
D^+ \bigl( \Rep \bigl(G(\BQ) \bigr) \bigr) 
                             \ar[rr]^-{\mu_K} \ar[dd]_-{\Res^{G(\BQ)}_{H_1}} &&
                  D^+ \Bigl( \Loc \bigl( M^K (G,\FX) (\BC) \bigr) \Bigr) 
                                  \ar[d]^-{\CV \mapsto \CV^{\partial BS}} \\
         && D^+_{\Loc} \bigl( \partial M^K  (G,\FX) (\BC)^{BS} \bigr)  
                                  \ar[d]^-{\res} \\
         D^+ ( \Rep H_1 ) \ar[dd]_-{\coh} 
         && \CC_{(G,\FX)}^{K,BS} \ar[d]^-{Rp^K_*} \\
         && \CC_{(G,\FX)}^{K,*} \ar[d]^-{i^*_{Q_1,g}} \\
\CR_{(Q_1,\FX) \tei g}^K \ar[rr]^-{\nu_g} && \CC_{(Q_1,\FX) \tei g}^K
\\}}
\quad .
\]
Next, according to Auxiliary Construction~\ref{7Aux}~(a), the upper right corner 
of our desired diagram equals
\[
\vcenter{\xymatrix@R-10pt{
D^+ \bigl( \Rep \bigl(G(\BQ) \bigr) \bigr) \ar[rr]^-{\mu_K^{\partial BS}} &&
            D^+ \Bigl( \Loc \bigl( \partial M^K  (G,\FX) (\BC)^{BS} \bigr) \Bigr)\ar[d] \\
         && D^+_{\Loc} \bigl( \partial M^K  (G,\FX) (\BC)^{BS} \bigr)  
\\}}
\quad .
\]
The functor $Rp^K_*$ induces a functor, denoted by the same symbol
\[
Rp^K_*: \CC_{(Q_1,\FX) \tei g}^{K,BS} \longto \CC_{(Q_1,\FX) \tei g}^K \; ,        
\]
and we have
\[
i^*_{Q_1,g} \circ Rp^K_* = Rp^K_* \circ k^*_{Q_1,g}: 
                \CC_{(G,\FX)}^{K,BS} \longto \CC_{(Q_1,\FX) \tei g}^K \; .
                \]
The right contour of our desired diagram thus transforms into  
\[
\vcenter{\xymatrix@R-10pt{
D^+ \bigl( \Rep \bigl(G(\BQ) \bigr) \bigr) \ar[rr]^-{\mu_K^{\partial BS}} &&
            D^+ \Bigl( \Loc \bigl( \partial M^K  (G,\FX) (\BC)^{BS} \bigr) \Bigr)\ar[d] \\
         && D^+_{\Loc} \bigl( \partial M^K  (G,\FX) (\BC)^{BS} \bigr)  
                                  \ar[d]^-{\res} \\
         && \CC_{(G,\FX)}^{K,BS} \ar[d]^-{k^*_{Q_1,g}} \\
         && \CC_{(Q_1,\FX) \tei g}^{K,BS} \ar[d]^-{Rp^K_*} \\
         && \CC_{(Q_1,\FX) \tei g}^K
\\}}
\quad .
\]
Given Proposition~\ref{7G}, it remains to establish natural commutativity of 
\[
\vcenter{\xymatrix@R-10pt{
D^+ ( \Rep H_1 ) \ar[rr]^-{\eta_g} \ar[d]_-{\coh} && 
                                        \CC_{(Q_1,\FX) \tei g}^{K,BS} \ar[d]^-{Rp^K_*} \\
\CR_{(Q_1,\FX) \tei g}^K \ar[rr]^-{\nu_g} && \CC_{(Q_1,\FX) \tei g}^K
\\}}
\quad .
\]
Let us apply Scholie~\ref{7Scho}
to the groups $H = H_C$ and $H' = H_1 \,$, with 
varying $X = X_Q \,$, depending on parabolic sub-groups $Q$ satisfying  $\adm(Q) = Q_1 \,$.

In order to define the $X_Q \,$, let
\[
T_{Q_1} := e(Q_1)' \subset \FX^{BS} \; ,
\]
and 
\[
T_Q := \bigcup_{q_1 \in Q_1(\BQ) / Q(\BQ)} e(q_1 Q q_1^{-1})' \subset T_{Q_1}
\]
whenever $Q$ is parabolic, with $\adm(Q) = Q_1 \,$.
We leave it to the reader to show, using \cite[Prop.~3.7]{W1}, that the 
union ``$\bigcup_{q_1 \in Q_1\BQ) / Q(\BQ)}$''
is actually disjoint, \emph{i.e.}, we have 
\[
T_Q = \coprod_{q_1 \in Q_1(\BQ) / Q(\BQ)} e(q_1 Q q_1^{-1})' \; .
\]
The sub-spaces $T_Q$ of $\FX^{BS}$ are stable under the action of $Q_1(\BQ)$: indeed,
for $q_1 , h \in Q_1(\BQ)$, the element $h$ of $Q_1(\BQ)$
maps $e(q_1 Q q_1^{-1})'$ to $e((h q_1) Q (h q_1)^{-1})'$.
Defining $pr_{gK}^{BS}$ by mapping $x \in \FX^{BS}$ to the class of 
$(x,gK)$ in $M^K  (G,\FX) (\BC)^{BS}$, the diagram
\[
\vcenter{\xymatrix@R-10pt{
T_Q \ar@{->>}[r]^-{pr_{gK}^{BS}} \ar@{^{ (}->}[d] & 
\bigl( p^K \bigr)^{-1} \bigl(pr_{gK} \bigl(\coprod \FX_1/W_1 \bigr) \bigr) 
              \cap e^K \bigl( Q,G(\BA_f) \bigr)' \ar@{^{ (}->}[d] \\
T_{Q_1}  \ar@{->>}[r]^-{pr_{gK}^{BS}} \ar@{->>}[d]_-{p} & 
(p^K)^{-1} \bigl( pr_{gK} \bigl( \coprod \FX_1/W_1 \bigr) \bigr) \ar@{->>}[d]^-{p^K} \\
\coprod \FX_1/W_1 \ar@{->>}[r]^-{pr_{gK}} & 
pr_{gK} \bigl( \coprod \FX_1/W_1 \bigr) 
\\}}
\]
is commutative, and its upper half is Cartesian. By 
\cite[Prop.~7.22~(a)]{W1} and Proposition~\ref{7I},
the maps $pr_{gK}^{BS}$ identify their targets with the quotients of their sources
by the (free) action of $H_1$. 

Write $Z := \coprod \FX_1/W_1 \,$, and fix an element $z_0 \in Z$. Define
\[
X_Q := p^{-1}(z_0) \cap T_Q \; ,
\]
for every parabolic $Q$ satisfying $\adm(Q) = Q_1 \,$. In particular
(Proposition~\ref{7I}), we have $X_{Q_1} = p^{-1}(z_0)$.
According to \cite[Thm.~7.20]{W1} and Proposition~\ref{7I}, 
there is an $H_C$-equivariant isomorphism
\[
k: X_{Q_1} \times Z = p^{-1}(z_0) \times Z \isoto T_{Q_1} \; ,
\]
that restricts to give
\[
X_Q \times Z \isoto T_Q \; ,
\]
for any parabolic $Q \in \adm^{-1}(Q_1)$, the $H_C$-action on the product being dia\-gonal,
with the trivial action on $Z$. Furthermore, under $k$,
the map $p$ is identified with projection $p^{X_{Q_1}}$ to the second factor.

Transport of structure \emph{via} $k$ provides for an action of $H_1$ on $X_{Q_1}  \times Z$,
extending the diagonal action of $H_C$.
With respect to this action,
\[
p^{X_{Q_1}}: X_{Q_1}  \times Z \longonto Z
\]
is $H_1$-equivariant (the action on the target factorizing through $H_1/H_C$).

For $Q \in \adm^{-1}(Q_1)$, consider the projection $\pi_{0,Q}$ from
\[
X_Q = p^{-1}(z_0) \cap T_Q
    = \coprod_{q_1 \in Q_1(\BQ) / Q(\BQ)} \bigl( p^{-1}(z_0) \cap 
                                                    e(q_1 Q q_1^{-1})' \bigr)
\]
to $Q_1(\BQ) / Q(\BQ)$. By Proposition~\ref{7I}, we have
\[ 
p^{-1}(z_0) \cap e(q_1 Q q_1^{-1})' 
        = p^{-1}(z_0) \cap \overline{e(q_1 Q q_1^{-1})}
\]
($\overline{e(q_1 Q q_1^{-1})}:=$ the closure of $e(q_1 Q q_1^{-1})$),
for any $q_1 \in Q_1(\BQ)$. According to \cite[Cor.~6.12~(b)]{W1},
each intersection $p^{-1}(z_0) \cap \overline{e(q_1 Q q_1^{-1})}$ is
contractible. In parti\-cu\-lar, the projection $\pi_{0,Q}$ induces an isomorphism between
$\pi_0(X_Q)$ and $Q_1(\BQ) / Q(\BQ)$. For an inclusion $Q \subset \tilde{Q}$
in $\adm^{-1}(Q_1)$, the isomorphisms induced by $\pi_{0,Q}$ and $\pi_{0,\tilde{Q}}$ 
fit into a commutative diagram 
\[
\vcenter{\xymatrix@R-10pt{
\pi_0(X_Q) \ar[r]^-{\cong} \ar[d] & Q_1(\BQ) / Q(\BQ) \ar@{^{>>}}[d] \\
\pi_0(X_{\tilde{Q}}) \ar[r]^-{\cong} & Q_1(\BQ) / \tilde{Q}(\BQ)
\\}}
\]
(with the obvious vertical maps).

It remains to verify
condition~\ref{7Scho}~(i). Thus, let 
\[
x \in p^{-1}(z_0) \cap e(q_1 Q q_1^{-1})' \subset X_Q \; , \;
\]
for some $q_1 \in Q_1(\BQ)$, and $h \in H_1$. We claim that
\[
h(x,z) \in \bigl( p^{-1}(z_0) \cap e((hq_1) Q (hq_1)^{-1})' \bigr) \times Z
\]
whenever $z \in Z =  \coprod \FX_1/W_1 \,$. Indeed, the point $k(x,z) \in T_Q$
belongs to $e(q_1 Q q_1^{-1})'$ \cite[Thm.~7.20~(b)]{W1}. Therefore,
\[
h (k(x,z)) \in e((hq_1) Q (hq_1)^{-1})' \; .
\]
Again by \cite[Thm.~7.20~(b)]{W1}, the image of the latter under $k^{-1}$ is
in 
\[
\bigl( p^{-1}(z_0) \cap e((hq_1) Q (hq_1)^{-1})' \bigr) \times Z \; ,
\]
as claimed. In other words, the action of $H_1$ on
\[
\pi_0(X_Q \times Z) = \pi_0(X_Q) \times \pi_0(Z) 
                    = \bigl( Q_1(\BQ) / Q(\BQ) \bigr) \times \pi_0(Z)
\]
is diagonal; more precisely, when $\pi_0(X_Q)$ is identified with $Q_1(\BQ) / Q(\BQ)$,
then the action of $H_1$ is by multiplication from the left.

Altogether, all conditions from Scholie~\ref{7Scho} are satisfied. To conclude the proof of part~(a),
apply principle~\ref{7Scho}~(3) to $Q = \bigcap \uQ$, for varying
$\uQ \in \CC_{(G,\FX) \tei Q_1}$.
\forget{
The sub-group $K$ of $G(\BA_f)$ being neat, we may apply \cite[Prop.~7.22~(c)]{W1}, 
to conclude that the induced diagram
\[
\vcenter{\xymatrix@R-10pt{
H_C \backslash T_Q \ar@{->>}[r]^-{pr_{gK}^{BS}} \ar@{^{ (}->}[d] & 
\bigl( p^K \bigr)^{-1} \bigl(pr_{gK} \bigl(\coprod \FX_1/W_1 \bigr) \bigr) 
              \cap e^K \bigl( Q,G(\BA_f) \bigr)' \ar@{^{ (}->}[d] \\
H_C \backslash T_{Q_1}  \ar@{->>}[r]^-{pr_{gK}^{BS}} \ar@{->>}[d]_-{p} & 
(p^K)^{-1} \bigl( pr_{gK} \bigl( \coprod \FX_1/W_1 \bigr) \bigr) \ar@{->>}[d]^-{p^K} \\
\coprod \FX_1/W_1 \ar@{->>}[r]^-{pr_{gK}} & 
pr_{gK} \bigl( \coprod \FX_1/W_1 \bigr) 
\\}}
\]
is Cartesian (including its lower half). 

Let $\BV \in D^+ (\Rep H_1)$. We have
$\coh(\BV) = ( \BW_{\uQ} , a_{\uQ \subset \utQ} )_{\uQ, \uQ \subset \utQ} \,$, with
\[
\BW_{\uQ} = R \Gamma \bigl( H_C, {}_{\bigcap \uQ} \! \BV \bigr) \; , \;
\]
(Construction~\ref{7Con}). Fix $\uQ \in \CC_{(G,\FX) \tei Q_1}$. 
The map
\[
p^K : \bigl( p^K \bigr)^{-1} \bigl(pr_{gK} \bigl(\coprod \FX_1/W_1 \bigr) \bigr) 
              \cap e^K \bigl( \bigcap \uQ,G(\BA_f) \bigr)'
      \longto pr_{gK} \bigl(\coprod \FX_1/W_1 \bigr)         
\]
is a locally trivial fibration \cite[Comp.~7.25]{W1}. As for its fibres,
according to \cite[Thm.~6.22]{W1}
\? \?
\[
p^{-1} \bigl( \coprod \FX_1/W_1 \bigr) 
= \bigcup_{P_1 \subset R \subset Q_1} e(R) \subset \FX^{BS} \; ,
\]
where the disjoint union runs over the 
parabolics $R$ of $G$ contained in $Q_1$ and containing
$P_1 \,$, the canonical normal sub-group of $Q_1$ 
\cite[Comp.~4.15]{W1}. Defining $e(Q_1)':= \cup_{P_1 \subset R \subset Q_1} e(R)$,
}
\end{Proofof}

Let $F$ be a field of characteristic zero. 
The canonical construction $\mu_K$ admits an $F$-linear, algebraic variant
\[
\Rep_F G \longto \Loc \bigl( M^K (G,\FX) (\BC) \bigr) 
\]
($\Rep_F G =$ the
category of algebraic representations of $G$ in finite dimensional $F$-vector spaces),
denoted by the same symbol $\mu_K$, and obtained by composing
\[
\mu_{K}: \Rep \bigl(G(\BQ) \bigr) \longto \Loc \bigl( M^K (G,\FX) (\BC) \bigr)
\] 
with the restriction
\[
\Res^G_{G(\BQ)}: \Rep_F G \longto \Rep \bigl(G(\BQ) \bigr) \; .
\]
It induces a triangulated functor between derived categories, which will still
be denoted by the symbol $\mu_K$. 
Using Proposition~\ref{7Cona} to factorize the functor $\coh$,
the natural commutative diagram from 
Main Theorem~\ref{7MT} restricts to
\[
\vcenter{\xymatrix@R-10pt{
D^+ \bigl( \Rep_F G \bigr) \ar[rr]^-{\mu_K} \ar[d]_-{\Res^G_{G(\BQ)}} &&
                  D^+ \Bigl( \Loc \bigl( M^K (G,\FX) (\BC) \bigr) \Bigr) \ar[d] \\
D^+ \bigl( \Rep \bigl(G(\BQ) \bigr) \bigr) \ar[d]_-{\Res^{G(\BQ)}_{H_1}} && 
                  D^+ \bigl( M^K  (G,\FX) (\BC) \bigr) \ar[d]^-{i^* R(i_G)_*} \\
D^+ ( \Rep H_1 ) \ar[d]_-{R \Gamma ( H_W, \argdot )} &&
                  D^+ \bigl( \partial M^K  (G,\FX)^* (\BC) \bigr) \ar[d]^-{\deg} \\
D^+ \bigl( \Rep H_1/H_W \bigr) \ar[d]_-{\overline{\coh}} && 
                  \CC_{(G,\FX)}^{K,*} \ar[d]^-{i^*_{Q_1,g}} \\
\CR_{(Q_1,\FX) \tei g}^K \ar[rr]^-{\nu_g} && \CC_{(Q_1,\FX) \tei g}^K
\\}}
\quad ,
\]
where as before $H_W$ denotes the sub-group $W_1(\BQ) \cap H_1$. 
Denote by 
\[
R \Gamma ( W_1, \argdot ) : 
  D^+ \bigl( \Rep_F Q_1 \bigr) \longto D^+ \bigl( \Rep_F \bar{Q}_1 \bigr)  
\]
the derived functor of
$\Gamma ( W_1, \argdot ) : \Rep_F Q_1 \longto \Rep_F \bar{Q}_1 \,$. 
Here as before, the symbol
$\bar{Q}_1$ denotes the maximal reductive quotient of $Q_1 \,$, \emph{i.e.} 
\cite[proof of Lemma~4.8]{P}, the quotient of $Q_1$ by $W_1 \,$.

\begin{Var} \label{7Var1}
The diagram from Main Theorem~\ref{7MT}~(a) restricts to give a natural
commutative diagram
\[
\vcenter{\xymatrix@R-10pt{
D^+ \bigl( \Rep_F G \bigr)  \ar[rr]^-{\mu_K} \ar[d]_-{\Res^G_{Q_1}} &&
                  D^+ \Bigl( \Loc \bigl( M^K (G,\FX) (\BC) \bigr) \Bigr) \ar[d] \\
D^+ \bigl( \Rep_F Q_1 \bigr) \ar[d]_-{R \Gamma ( W_1, \argdot )} && 
                  D^+ \bigl( M^K  (G,\FX) (\BC) \bigr) \ar[d]^-{i^* R(i_G)_*} \\
D^+ \bigl( \Rep_F \bar{Q}_1 \bigr) \ar[d]_-{\Res^{\bar{Q}_1}_{H_1/H_W}} &&
                  D^+ \bigl( \partial M^K  (G,\FX)^* (\BC) \bigr) \ar[d]^-{\deg} \\
D^+ \bigl( \Rep H_1/H_W \bigr) \ar[d]_-{\overline{\coh}} &&  
                  \CC_{(G,\FX)}^{K,*} \ar[d]^-{i^*_{Q_1,g}} \\
\CR_{(Q_1,\FX) \tei g}^K \ar[rr]^-{\nu_g} && \CC_{(Q_1,\FX) \tei g}^K
\\}}
\]
($\Res^G_{Q_1}$ and $\Res^{\bar{Q}_1}_{H_1/H_W} :=$ 
the restrictions from $G$ to $Q_1$ and from $\bar{Q}_1$ to $H_1/H_W$, respectively). 
\end{Var}

\begin{Proof}
Use the diagram preceding our statement, together with the fact that as functors
\[
D^+ \bigl( \Rep_F Q_1 \bigr) \longto D^+ \bigl( \Rep H_1/H_W \bigr) \; ,
\]
the compositions
$\Res^{\bar{Q}_1}_{H_1/H_W} \circ R \Gamma ( W_1, \argdot )$ and
$R \Gamma ( H_W, \argdot ) \circ \Res^{Q_1}_{H_1}$
($\Res^{Q_1}_{H_1} :=$ the restriction from $Q_1$ to $H_1$)
are canonically isomorphic, as $H_W$ is an arithmetic sub-group of the unipotent
group $W_1$.
\end{Proof}

\forget{
By Construction~\ref{7Con}~(b), the functor $\overline{\coh}$ encodes 
cohomology of the quotient $H_C/H_W$, and it will be useful to study the
behaviour of this quotient when $g \in G(\BA_f)$ varies. Recall that by definition,
\[
H_C = H_C(gK) = C_1 (\BQ) \cap gKg^{-1} 
\]
and 
\[
H_W = H_W(gK) = W_1(\BQ) \cap gKg^{-1} \; .
\]
As suggested by the notation, $H_C$ and $H_W$ remain invariant when $g$ is replaced
by an element in $gK$; the same invariance therefore holds for $H_C/H_W$.}

\begin{Def} \label{7K}
(a)~Define $M^K (Q_1;P_1(\BA_f)gK,\FX) \subset M^K (Q_1,\FX)$ as the image 
of $\coprod \FX_1/W_1 \times P_1(\BA_f)gK/K$ under the projection 
\[
\FX^* \times G (\BA_f) / K 
\longonto G (\BQ) \backslash \bigl( \FX^* \times G (\BA_f) / K \bigr) = M^K (G,\FX)^* (\BC) \; .
\]
(b)~Define $i_{Q_1;P_1(\BA_f)gK}$ to be the
immersion of $M^K (Q_1;P_1(\BA_f)gK,\FX)$ into $\partial M^K  (G,\FX)^*$. \\[0.1cm]
(c)~Define a category $\CC_{(Q_1,\FX) \tei P_1(\BA_f)g}^K$ as follows: 
objects are of the form
\[
\bigl( \CK_{\uQ} , a_{\uQ \subset \utQ} \bigr)_{\uQ, \utQ \in \CC_{(G,\FX) \tei Q_1}} \; ,
\]
where $\CK_{\uQ} \in D^+(M^K (Q_1;P_1(\BA_f)gK,\FX)(\BC))$,
for each $\uQ \in \CC_{(G,\FX) \tei Q_1} \,$, and 
\[
a_{\uQ \subset \utQ}: R (i_{Q_1;P_1(\BA_f)gK})_* \CK_{\utQ} 
                                    \longto R (i_{Q_1;P_1(\BA_f)gK})_* \CK_{\uQ}
\]
is a morphism in $D^+(M^K  (G,\FX)^* (\BC))$, for each relation $\uQ \subset \utQ$
in $\CC_{(G,\FX) \tei Q_1}$.
Morphisms are defined as in $\CC_{(G,\FX)}^K$ (Definition~\ref{3E}~(a)). \\[0.1cm]
(d)~Define
\[
i^*_{Q_1;P_1(\BA_f)gK}: \CC_{(G,\FX)}^{K,*} \longto \CC_{(Q_1,\FX) \tei P_1(\BA_f)g}^K
\]
as $i^*_{Q_1}$ (Definition~\ref{3E}~(c)), followed by
the restriction from $M^K (Q_1,\FX)(\BC)$ to $M^K (Q_1;P_1(\BA_f)gK,\FX)(\BC)$. \\[0.1cm]
(e)~Define 
\[
H_1' := Q_1(\BQ) \cap \bigl( P_1(\BA_f) \cdot gKg^{-1} \bigr) 
( = Q_1(\BQ) \cap \bigl( P_1(\BA_f) \cdot (Q_1(\BA_f) \cap gKg^{-1}) \bigr) ) \; .
\]
\end{Def}

Note that $H_1'$ is the stabilizer in $G(\BQ)$ of 
$\coprod \FX_1/W_1 \times P_1(\BA_f)gK/K$. Therefore, the quotient of the action
of $\coprod \FX_1/W_1 \times P_1(\BA_f)gK/K$ by $H_1'$ equals
$M^K (Q_1;P_1(\BA_f)gK,\FX)(\BC)$.
The space $M^K (Q_1;P_1(\BA_f)gK,\FX)(\BC)$ contains
the set $pr_{gK}(\coprod \FX_1/W_1)$ considered so far. Actually, we have
\[
M^K (Q_1;P_1(\BA_f)gK,\FX)(\BC) = \coprod_{p_1} pr_{(p_1g)K}(\coprod \FX_1/W_1) \; ,
\]
where $p_1$ runs through a set of representatives of the (finite) double co-set
\[
H_1' \backslash P_1(\BA_f) / \bigl( P_1(\BA_f) \cap gKg^{-1} \bigr) \; .
\]
Here, $H_1' \subset P_1(\BA_f) \cdot (Q_1(\BA_f) \cap gKg^{-1})$
acts on the quotient 
\[
P_1(\BA_f) / \bigl( P_1(\BA_f) \cap gKg^{-1} \bigr)
\]
by translation from the left, once the latter is identified with
\[
\bigl( P_1(\BA_f) \cdot (Q_1(\BA_f) \cap gKg^{-1}) \bigr) / 
( Q_1(\BA_f) \cap gKg^{-1} )  \; .
\]
The space $M^K (Q_1;P_1(\BA_f)gK,\FX)(\BC)$ is open and closed in 
$M^K (Q_1,\FX)(\BC)$.

\begin{Rem} \label{7L}
(a)~By definition, the space $M^K (Q_1;P_1(\BA_f)gK,\FX)(\BC)$ is the union of the
images $M^K (Q_1;\FX_1/W_1;P_1(\BA_f)gK,\FX)(\BC)$
of the space $\FX_1/W_1 \times P_1(\BA_f)gK/K$ under 
\[
\FX^* \times G (\BA_f) / K 
\longonto G (\BQ) \backslash \bigl( \FX^* \times G (\BA_f) / K \bigr) = M^K (G,\FX)^* (\BC) \; ,
\]
for the individual boundary components $(P_1,\FX_1)$ associated to $Q_1$  
(some of the $M^K (Q_1;\FX_1/W_1;P_1(\BA_f)gK,\FX)(\BC)$ may actually be identified
under the action of $H_1'$).
For each $(P_1,\FX_1)$, we have
\[
M^K (Q_1;\FX_1/W_1;P_1(\BA_f)gK,\FX)(\BC) = 
\Delta_1 \backslash M^{\pi_1(K_f^1)} ((P_1,\FX_1)/W_1)(\BC)
\]
in the notation introduced in \cite[Sect.~6.3]{P}
(put $K_f := K$ and $p_f:= g$ in \loccit). 
The same object is denoted by $M_1^K$ in \cite{BW}. \\[0.1cm]
(b)~\emph{A priori}, Definition~\ref{7K}~(a) concerns an open and closed
subset 
\[
M^K (Q_1;P_1(\BA_f)gK,\FX)(\BC)
\] 
of the space of $\BC$-valued
points $M^K (Q_1,\FX)(\BC)$ of $M^K (Q_1,\FX)$. 
Given the observation from (a),
it follows from \cite[Main Theorem~12.3 for the Baily--Borel compactification]{P} that
this subset is indeed identified with the set of $\BC$-valued points of an
open and closed subscheme $M^K (Q_1;P_1(\BA_f)gK,\FX)$ of $M^K (Q_1,\FX)$. 
\end{Rem}

The second variant of Main Theorem~\ref{7MT}~(a) will give a natural commutative
diagram completing
\[
\vcenter{\xymatrix@R-10pt{
D^+ \bigl( \Rep_F G \bigr)  \ar[rr]^-{\mu_K} &&
                  D^+ \Bigl( \Loc \bigl( M^K (G,\FX) (\BC) \bigr) \Bigr) \ar[d] \\
                  && 
                  D^+ \bigl( M^K  (G,\FX) (\BC) \bigr) \ar[d]^-{i^* R(i_G)_*} \\
                  &&
                  D^+ \bigl( \partial M^K  (G,\FX)^* (\BC) \bigr) \ar[d]^-{\deg} \\
                  && \CC_{(G,\FX)}^{K,*} \ar[d]^-{i^*_{Q_1;P_1(\BA_f)gK}} \\
                  && \CC_{(Q_1,\FX) \tei P_1(\BA_f)g}^K
\\}} \quad ,
\]
the latter half diagram being identical, up to the last functor
$i^*_{Q_1;P_1(\BA_f)gK} \,$, to the right contour of the diagram from
Variant~\ref{7Var1}. 

\begin{Def} \label{7M}
(a)~Define 
\[
H_C' := C_1(\BQ) \cap \bigl( W_1(\BA_f) \cdot gKg^{-1} \bigr) 
( = C_1(\BQ) \cap \bigl( W_1(\BA_f) \cdot (C_1(\BA_f) \cap gKg^{-1}) \bigr) ) \; .
\]
(b)~Define a category $\CR_{(Q_1,\FX) \tei P_1(\BA_f)g}^K$ as follows: 
objects are of the form
\[
\bigl( \BW_{\uQ} , a_{\uQ \subset \utQ} \bigr)_{\uQ, \utQ \in \CC_{(G,\FX) \tei Q_1}} \; ,
\]
where $\BW_{\uQ} \in D^+ ( \Rep (H_1'/H_C') )$,
for each $\uQ \in \CC_{(G,\FX) \tei Q_1}$, and 
\[
a_{\uQ \subset \utQ}: \BW_{\utQ} \longto \BW_{\uQ}
\]
is a morphism in $D^+ ( \Rep (H_1'/H_C') )$, for each relation $\uQ \subset \utQ$
in $\CC_{(G,\FX) \tei Q_1}$. 
Morphisms are defined as in Definition~\ref{7E}.
\end{Def}

Construction~\ref{7Con}~(b) can be applied \emph{mutatis mutandis}, to yield a
functor
\[
\overline{\coh} : D^+ \bigl( \Rep \bigl( H_1'/(W_1(\BQ)) \bigr) \bigr)
 \longto \CR_{(Q_1,\FX) \tei P_1(\BA_f) g}^K 
\]
(for $\BV \in D^+ ( \Rep ( H_1'/(W_1(\BQ)) ) )$ the component $\BW_{\uQ}$
of
\[
\overline{\coh}(\BV) 
  = \bigl( \BW_{\uQ} , a_{\uQ \subset \utQ} \bigr)_{\uQ, \utQ \in \CC_{(G,\FX) \tei Q_1}} 
\]
is defined as $\BW_{\uQ}:= R \Gamma ( H_C'/(W_1(\BQ)), {}_{\bigcap \uQ} \! \BV )$)
($R \Gamma ( H_C'/(W_1(\BQ)), \argdot ) :=$ the derived functor of
\[
\Gamma \bigl( H_C'/(W_1(\BQ) ), \argdot \bigr) : \Rep \bigl( H_1'/(W_1(\BQ)) \bigr) \longto \Rep (H_1'/H_C') \; ).
\]
The sub-group $W_1(\BQ)$ of $H_1'$ is normal, and acts trivially on
$\coprod \FX_1/W_1 \,$. Therefore, the action of $H_1'$ on 
$\coprod \FX_1/W_1 \times P_1(\BA_f)gK/K$ induces an action on
\[
\coprod \FX_1/W_1 \times \bigl( W_1(\BA_f) \backslash P_1(\BA_f)gK/K \bigr)
\]
(whose restriction to $W_1(\BQ)$ is trivial). Strong approximation for unipotent
groups implies that the quotient remains unchanged, \emph{i.e.}, 
\[
H_1' \backslash 
\Bigl( \coprod \FX_1/W_1 \times \bigl( W_1(\BA_f) \backslash P_1(\BA_f)gK/K \bigr) \Bigl)
= M^K (Q_1;P_1(\BA_f)gK,\FX)(\BC) \; .
\]

\begin{Rem} \label{7Ma}
Define $\pi_1: P_1 \onto G_1$ as the canonical epimorphism of $P_1$ to its
maximal reductive quotient (\emph{i.e.}, as the restriction of $\pi_{Q_1}$
to $P_1$), and 
\[
K_1 := P_1(\BA_f) \cap gKg^{-1} \; .
\]
The map $p_1K_1 \mapsto p_1gK$ is a bijection
\[
P_1(\BA_f)/K_1 \longto P_1(\BA_f)gK/K \; ;
\]
therefore, the product 
\[
\coprod \FX_1/W_1 \times \bigl( W_1(\BA_f) \backslash P_1(\BA_f)gK/K \bigr)
\]
is identified with
\[
\coprod \FX_1/W_1 \times \bigl( W_1(\BA_f) \backslash P_1(\BA_f)/K_1 \bigr) \; ,
\]
which thanks to strong approximation for unipotent groups is in turn identified with
\[
\coprod \FX_1/W_1 \times \bigl( G_1(\BA_f)/\pi_1(K_1) \bigr) \; ,
\]
\emph{i.e.}, with the disjoint union of the spaces covering the Shimura varieties
of level $\pi_1(K_1)$ associated to the quotient Shimura data $(P_1,\FX_1)/W_1$
(cmp.\ Remark~\ref{7L}~(a)).
\end{Rem}

\begin{Prop} \label{7N}
(a)~We have 
\[
H_C' = W_1(\BQ) \cdot H_C 
\]  
and 
\[
H_C/H_W = H_C'/(W_1(\BQ)) \; .
\]
(b)~The action of the normal sub-group $H_C'$
of $H_1'$ on 
\[
\coprod \FX_1/W_1 \times \bigl( W_1(\BA_f) \backslash P_1(\BA_f)gK/K \bigr)
\]
is trivial. \\[0.1cm] 
(c)~The induced action of the quotient $H_1'/H_C'$
on 
\[
\coprod \FX_1/W_1 \times \bigl( W_1(\BA_f) \backslash P_1(\BA_f)gK/K \bigr)
\]
is free, and the quotient of this action equals 
\[
M^K (Q_1;P_1(\BA_f)gK,\FX)(\BC) \subset \partial M^K  (G,\FX)^* (\BC) \; .
\]
\end{Prop}

\begin{Proof}
(a): by definition, and by strong approximation for unipotent groups,
\[
H_C' = C_1(\BQ) \cap \bigl( W_1(\BQ) \cdot gKg^{-1} \bigr) \; ,
\]
which equals
\[     
       W_1(\BQ) \cdot \bigl( C_1(\BQ) \cap gKg^{-1} \bigr)
     = W_1(\BQ) \cdot H_C \; .
\]
Consequently, the inclusion of $H_C$ into $H_C'$ induces a canonical isomorphism
\[
H_C/H_W = H_C/ \bigl( W_1(\BQ) \cap H_C \bigr) \isoto H_C'/(W_1(\BQ)) \; .
\]
(b): given (a), it suffices to show that $H_C$ acts trivially on 
$\coprod \FX_1/W_1 \,$, and on $W_1(\BA_f) \backslash P_1(\BA_f)gK/K$.
As for $\coprod \FX_1/W_1 \,$, we refer to \cite[Prop.~7.22~(b)]{W1} (this
was already used before). Let $\gamma \in H_C$ and $p_1 \in P_1(\BA_f)$.
Since $\gamma \in C_1(\BQ)$, we have $\gamma p_1 = wp_1\gamma$, for some 
$w \in W_1(\BA_f)$.
Since $\gamma \in gKg^{-1}$, we have 
\[
\gamma p_1gK = wp_1\gamma gK = wp_1gK \; .
\]
(c): first, let us analyze the stabilizers of the action of $H_1'$
on
\[
\coprod \FX_1/W_1 \times \bigl( W_1(\BA_f) \backslash P_1(\BA_f)gK/K \bigr) \; .
\]
Let $z \in \FX_1/W_1 =: \FH_1 \,$, $p_1 \in P_1 (\BA_f)$, and
$\gamma \in H_1'$, and suppose that 
\[
\gamma \cdot (z,[p_1])] = (z,[p_1])
\]
in $\FH_1 \times ( W_1(\BA_f) \backslash P_1(\BA_f)gK/K )$. 
Write $\gamma = p_2 k$, with $p_2 \in P_1(\BA_f)$ and 
$k \in Q_1(\BA_f) \cap gKg^{-1}$. 
Thus, 
\begin{enumerate}
\item[(1)] $\gamma \in Q_1(\BQ)$ stabilizes $z$.
\item[(2)] There is $k_1 \in Q_1(\BA_f) \cap gKg^{-1}$ such that
$p_1 = \gamma p_1 k^{-1} k_1 \mod W_1 (\BA_f)$. 
Modulo $W_1(\BQ)$, the element $\gamma$ thus belongs to
$(p_1 g) K' (p_1 g)^{-1}$, which is a neat sub-group of
$G(\BA_f)$.
\end{enumerate}
Choose a complement $G_2$ of $G_1= P_1/W_1$ in $\bar{Q}_1 = Q_1/W_1 \,$.
The groups $G_1$
and $G_2$ centralize each other. Denote by $\Pi(\gamma)$ the image of
$\gamma$ in $Q/G_2$. We identify this group with
$G_1/G_1 \cap G_2$. Because of (1), the element $\Pi(\gamma)$
stabilizes a point in the space $\FH_1 / G_1 \cap G_2$
belonging to the quotient Shimura data
\[
(Q/G_2, \FH_1 / G_1 \cap G_2) := (G_1,\FH_1) / G_1 \cap G_2 \; .
\]
By \cite[Lemma~1.3]{BW}, the element $\Pi(\gamma)$ is of finite order. Because
of (2), it must be trivial. We conclude:
\begin{enumerate}
\item[(3)] The element $\gamma$ belongs to the group of rational points of
\[ 
C_1' := \{ q \in Q_1 \; , \; \pi_{Q_1}(q) \in \Cent_{\bar{Q}_1}(G_1) \} \; .
\]
\end{enumerate}
Since $G_1(\BR)$ acts transitively on $\FH_1$, (1) and (3) imply
that $\gamma$ acts trivially on $\FH_1$.
We conclude:
\begin{enumerate}
\item[(4)] The kernel $H_C''$ of the action of $H_1'$
on 
\[
\coprod \FX_1/W_1 \times \bigl( W_1(\BA_f) \backslash P_1(\BA_f)gK/K \bigr)
\]
is equal to any of its stabilizers.
\end{enumerate}
Second, we get from (2) and (3) that 
\[
p_1 k^{-1} k_1 = \gamma^{-1} p_1 =
p_1 \gamma^{-1} \mod W_1 (\BA_f) \; .
\] 
Thus,
\begin{enumerate}
\item[(5)] $H_C'' 
              = \Cent_{Q_1(\BQ)} (\coprod \FX_1/W_1) \cap (W_1(\BA_f) \cdot gKg^{-1})$.
\end{enumerate}
According to (3), the group $H_C''$ is contained in
\[
H_C''' := C_1'(\BQ) \cap \bigl( W_1(\BA_f) \cdot gKg^{-1} \bigr) \; .
\] 
Strong approximation shows (as for (a)) that
\[
H_C''' = W_1(\BQ) \cdot \bigl( C_1'(\BQ) \cap gKg^{-1} \bigr) \; . 
\]
But $K$ being neat, 
\[
C_1'(\BQ) \cap gKg^{-1} = C_1(\BQ) \cap gKg^{-1} = H_C
\]
(recall that by Definiton~\ref{7A}, the group $C_1$ is the neutral connected
component of $C_1'$). We conclude:
\begin{enumerate}
\item[(6)] The inclusions $H_C' \subset H_C'' \subset H_C'''$
are equalities. 
\end{enumerate}
In particular, by~(4) and (6), the action
of the quotient $H_1'/H_C'$
on 
\[
\coprod \FX_1/W_1 \times \bigl( W_1(\BA_f) \backslash P_1(\BA_f)gK/K \bigr)
\]
is without fixed points.

In order to show that the action is free, it suffices to establish that it is
properly discontinuous. According to \cite[Prop.~1.1~(a)]{BW}
(using (6), and taking care to adapt the notation of \loccit \
\emph{via} the second statement of \cite[Prop.~1.1~(b)]{BW}),
the sub-group $P_1(\BQ) H_C'$ of $H_1'$ is of finite index.
Hence so is the sub-group $G_1(\BQ)$ of $H_1'/H_C'$.
But $G_1(\BQ)$ acts properly discontinuously on 
\[
\coprod \FX_1/W_1 \times \bigl( G_1(\BA_f)/\pi_1(K_1) \bigr) \stackrel{\ref{7Ma}}{=}
\coprod \FX_1/W_1 \times \bigl( W_1(\BA_f) \backslash P_1(\BA_f)gK/K \bigr)
\]
\cite[Lemma~1.3]{BW}.
\end{Proof}

Part~(c) of Proposition~\ref{7N} allows to define 
\[
\mu_{\pi_1(K_1)}: \Rep(H_1'/H_C') \longto 
         \Loc \bigl( M^K (Q_1;P_1(\BA_f)gK,\FX)(\BC) \bigr) 
\] 
as the functor associating to each representation $\BV$
the sheaf of continuous sections of the projection from
\[
(H_1'/H_C') \backslash 
              \Bigl( \BV \times \coprod \FX_1/W_1  
                   \times \bigl( W_1(\BA_f) \backslash P_1(\BA_f)gK/K \bigr) \Bigr) 
\]
to
\[
(H_1'/H_C') \backslash 
              \Bigl( \coprod \FX_1/W_1  
                   \times \bigl( W_1(\BA_f) \backslash P_1(\BA_f)gK/K \bigr) \Bigr) \; ,
\]
the latter space being identified with $M^K (Q_1;P_1(\BA_f)gK,\FX)(\BC)$.
As usual, this functor induces
\[
\mu_{\pi_1(K_1)}: D^+ \bigl( \Rep(H_1'/H_C') \bigr) \longto 
D^+ \bigl( M^K (Q_1;P_1(\BA_f)gK,\FX)(\BC) \bigr) \; .
\]
Define
\[
\nu_{P_1(\BA_f)g}: \CR_{(Q_1,\FX) \tei P_1(\BA_f)g}^K
     \longto \CC_{(Q_1,\FX) \tei P_1(\BA_f)g}^K \; .
\]
as the componentwise extension of $\mu_{\pi_1(K_1)}$ to 
$\CR_{(Q_1,\FX) \tei P_1(\BA_f)g}^K \,$. 

\begin{Var} \label{7Var2}
(a)~There is a natural commutative diagram
\[
\vcenter{\xymatrix@R-10pt{
D^+ \bigl( \Rep_F G \bigr)  \ar[rr]^-{\mu_K} \ar[d]_-{\Res^G_{Q_1}} &&
                  D^+ \Bigl( \Loc \bigl( M^K (G,\FX) (\BC) \bigr) \Bigr) \ar[d] \\
D^+ \bigl( \Rep_F Q_1 \bigr) \ar[d]_-{R \Gamma ( W_1, \argdot )} && 
                  D^+ \bigl( M^K  (G,\FX) (\BC) \bigr) \ar[d]^-{i^* R(i_G)_*} \\
D^+ \bigl( \Rep_F \bar{Q}_1 \bigr) \ar[d]_-{\Res^{\bar{Q}_1}_{H_1'/(W_1(\BQ))}} &&
                  D^+ \bigl( \partial M^K  (G,\FX)^* (\BC) \bigr) \ar[d]^-{\deg} \\
D^+ \bigl( \Rep H_1'/(W_1(\BQ)) \bigr) \ar[d]_-{\overline{\coh}} &&  
                  \CC_{(G,\FX)}^{K,*} \ar[d]^-{i^*_{Q_1;P_1(\BA_f)gK}} \\
\CR_{(Q_1,\FX) \tei P_1(\BA_f)g}^K \ar[rr]^-{\nu_{P_1(\BA_f)g}} && 
                  \CC_{(Q_1,\FX) \tei P_1(\BA_f)g}^K
\\}}
\]
($\Res^{\bar{Q}_1}_{H_1'/(W_1(\BQ))} :=$ 
the restriction from $Q_1$ to $H_1'/(W_1(\BQ))$). \\[0.1cm]
(b)~Let $\BV \in D^+ ( \Rep_F G )$. There are isomorphisms $\alpha_{\uQ}$
from
\[
\CK_{\uQ}:=
i^*_{Q_1;P_1(\BA_f)gK} R(i_{Q_2})_* i_{Q_2}^* \ldots R(i_{Q_r})_* i_{Q_r}^* i^* R(i_G)_* 
\circ \mu_K (\BV)
\]
to
\[
\CL_{\uQ}:= \mu_{\pi_1(K_1)} \circ 
R \Gamma \bigl( H_C/H_W, {}_{\bigcap \uQ} R \Gamma ( W_1, \Res^G_{Q_1}\BV) \bigr)
\]
in the category $D^+ ( M^K (Q_1;P_1(\BA_f)gK,\FX)(\BC) )$, functorial in $\BV$ and
indexed by $\uQ = (Q_1 \prec Q_2 \prec \ldots \prec Q_r) \in \CC_{(G,\FX) \tei Q_1}$,
such that for each relation $\uQ \subset \utQ$ in $\CC_{(G,\FX) \tei Q_1}$, the diagram
\[
\vcenter{\xymatrix@R-10pt{
\CK_{\utQ} \ar[r]^-{a_{\uQ \subset \utQ}} 
                                   \ar[d]_-{\alpha_{\utQ}} 
  & \CK_{\uQ} \ar[d]^-{\alpha_{\uQ}}  \\
\CL_{\utQ} \ar[r]^-{b_{\uQ \subset \utQ}}  & \CL_{\uQ}
\\}}
\] 
commutes. Here, the transition morphisms $a_{\uQ \subset \utQ}$ are
given by adjunctions as in Construction~\ref{6Con}, 
and the $b_{\uQ \subset \utQ}$ are induced by the natural transformations
${}_{\bigcap \utQ} \! \argdot \to {}_{\bigcap \uQ} \! \argdot$
as in Construction~\ref{7Con}~(b).
\end{Var}

\begin{Proof}
Part~(b) is merely a reformulation of (a), given the definitions,
and the identification of $H_C/H_W$ with $H_C'/(W_1(\BQ))$ from Proposition~\ref{7N}~(a).
 
The space $M^K (Q_1;P_1(\BA_f)gK,\FX)(\BC)$ is covered by 
open and closed sub-sets of the form 
$pr_{p_1gK}(\coprod \FX_1/W_1)$, for $p_1 \in P_1(\BA_f)$. It is therefore sufficient to
prove commutativity of the diagram of (a) after applying the restriction from
$\CC_{(Q_1,\FX) \tei P_1(\BA_f)g}^K$ to $\CC_{(Q_1,\FX) \tei p_1g}^K$
induced by the inclusion 
\[
pr_{p_1gK}(\coprod \FX_1/W_1) \into M^K (Q_1;P_1(\BA_f)gK,\FX)(\BC) \; ,
\]
for any $p_1 \in P_1(\BA_f)$. 

Let us identify this restricted diagram: defining
\[
H_1 (p_1gK) := Q_1(\BQ) \cap (p_1g)K(p_1g)^{-1} \; ,
\]
\[
H_W (p_1gK) := W_1(\BQ) \cap H_1 (p_1gK) \; ,
\]
and noting that $H_1 (p_1gK)$ is the stabilizer of
$\coprod \FX_1/W_1 \times \{ p_1gK \}$ in $H_1' \,$, 
we obtain a diagram
\[
\vcenter{\xymatrix@R-10pt{
D^+ \bigl( \Rep_F G \bigr)  \ar[rr]^-{\mu_K} \ar[d]_-{\Res^G_{Q_1}} &&
                  D^+ \Bigl( \Loc \bigl( M^K (G,\FX) (\BC) \bigr) \Bigr) \ar[d] \\
D^+ \bigl( \Rep_F Q_1 \bigr) \ar[d]_-{R \Gamma ( W_1, \argdot )} && 
                  D^+ \bigl( M^K  (G,\FX) (\BC) \bigr) \ar[d]^-{i^* R(i_G)_*} \\
D^+ \bigl( \Rep_F \bar{Q}_1 \bigr) \ar[d]_-{\Res^{\bar{Q}_1}_{H_1(p_1gK)/H_W(p_1gK)}} &&
                  D^+ \bigl( \partial M^K  (G,\FX)^* (\BC) \bigr) \ar[d]^-{\deg} \\
D^+ \bigl( \Rep H_1(p_1gK)/H_W(p_1gK) \bigr) \ar[d]_-{\overline{\coh}'} &&  
                  \CC_{(G,\FX)}^{K,*} \ar[d]^-{i^*_{Q_1,p_1g}} \\
\CR_{(Q_1,\FX) \tei p_1g}^K \ar[rr]^-{\nu_{p_1g}} && \CC_{(Q_1,\FX) \tei p_1g}^K
\\}}
\]
that is identical to the one from Variant~\ref{7Var1},
where $g$ is replaced by $p_1g$, up to the
functor $\overline{\coh}'$. Indeed, the latter involves cohomology
of the quotient $H_C'/W_1(\BQ)$, while the functor $\overline{\coh}$ from
Variant~\ref{7Var1}
involves cohomology of $H_C (p_1gK)/H_W (p_1gK)$, where
\[
H_C (p_1gK) := C_1(\BQ) \cap (p_1g)K(p_1g)^{-1} \; .
\]
We leave it to the reader to show that (1)~$H_C (p_1gK)$ is a sub-group of $H_C'$,
(2)~$H_C'$ does not change when $g$ is replaced by $p_1g$, \emph{i.e.},
\[
H_C' = C_1(\BQ) \cap \bigl( W_1(\BA_f) \cdot (p_1g)K(p_1g)^{-1} \bigr) \; .
\] 
In order to conclude \emph{via} Variant~\ref{7Var1}, it suffices to note that
the inclusion $H_C (p_1gK) \into H_C'$ induces an isomorphism
\[
H_C (p_1gK)/H_W (p_1gK) \isoto H_C'/W_1(\BQ) 
\] 
(Proposition~\ref{7N}~(a), applied to $p_1g$ instead of $g$).
\end{Proof}

\begin{Rem} \label{7O}
The diagram from Variant~\ref{7Var2}~(a) restricts to give na\-tural commutative
diagrams concerning 
\[
M^K (Q_1;\FX_1/W_1;P_1(\BA_f)gK,\FX) \subset 
M^K (Q_1;P_1(\BA_f)gK,\FX) 
\]
(see Remark~\ref{7L}~(a)),
for each individual boundary component $(P_1,\FX_1)$ associated to $Q_1$.
For the correct equivariance statement, the group $H_1'$ needs to replaced by the
stabilizer in $G(\BQ)$ of 
$\FX_1/W_1 \times P_1(\BA_f)gK/K$, \emph{i.e.}, by
\[
\Stab_{Q_1(\BQ)}(\FX_1/W_1) \cap \bigl( P_1(\BA_f) \cdot gKg^{-1} \bigr) \subset H_1' 
\]
(this stabilizer is denoted by $H_Q$ in \cite{BW}).
We leave the details to the reader.
\end{Rem}

\begin{Cor} \label{7P}
Let $\BV \in D^+ ( \Rep_F G )$. \\[0.1cm]
(a)~There are canonical and
functorial $E_2$-spectral sequences $(E^\argast_{\uQ})$, 
indexed by $\uQ = (Q_1 \prec Q_2 \prec \ldots \prec Q_r) \in \CC_{(G,\FX) \tei Q_1} \,$,
\[
E^{p,s}_{2,\uQ} \Longrightarrow \CH^{p+s} \bigl( i^*_{Q_1;P_1(\BA_f)gK} R(i_{Q_2})_* 
i_{Q_2}^* \ldots R(i_{Q_r})_* i_{Q_r}^* i^* R(i_G)_* \bigr) \circ \mu_K (\BV) \; ,
\]
where
\[
E^{p,s}_{2,\uQ} := \mu_{\pi_1(K_1)} \circ 
     H^p \bigl( H_C/H_W, {}_{\bigcap \uQ} H^s ( W_1, \Res^G_{Q_1}\BV) \bigr) \; .
\]
(b)~For each relation $\uQ \subset \utQ$ in $\CC_{(G,\FX) \tei Q_1}$, there
is a canonical and functorial morphism of spectral sequences
$(E^\argast_{\utQ}) \to (E^\argast_{\uQ})$. 
On the $E_2$-terms, it is induced by the natural transformation
${}_{\bigcap \utQ} \! \argdot \to {}_{\bigcap \uQ} \! \argdot$
as in Construction~\ref{7Con}~(b).
On the end terms, it is given by adjunction as in Construction~\ref{6Con}.
\end{Cor}

\begin{Thm} \label{7Q}
The spectral sequences of Corollary~\ref{7P} degenerate and split canonically
in a compatible way. More precisely, let $\BV \in D^+ ( \Rep_F G )$. \\[0.1cm]
(a)~For any $n \in \BZ$, and any $\uQ = (Q_1 \prec Q_2 \prec \ldots \prec Q_r) \in \CC_{(G,\FX) \tei Q_1} \,$,
there are canonical and
functorial isomorphisms of local systems on the space
$M^K (Q_1;P_1(\BA_f)gK,\FX)$ between
\[
\CH^{n} \bigl( i^*_{Q_1;P_1(\BA_f)gK} R(i_{Q_2})_* 
i_{Q_2}^* \ldots R(i_{Q_r})_* i_{Q_r}^* i^* R(i_G)_* \bigr) \circ \mu_K (\BV)
\]
and 
\[
\bigoplus_{p+s = n} \mu_{\pi_1(K_1)} \circ 
     H^p \bigl( H_C/H_W, {}_{\bigcap \uQ} H^s ( W_1, \Res^G_{Q_1}\BV) \bigr) \; .
\]
(b)~For each relation 
$\uQ \subset \utQ = (Q_1=\tQ_1 \prec \tQ_2 \prec \ldots \prec \tQ_t)$ in 
$\CC_{(G,\FX) \tei Q_1} \,$, 
and any $n \in \BZ$, the adjunction from
\[
\CH^{n} \bigl( i^*_{Q_1;P_1(\BA_f)gK} R(i_{\tQ_2})_* 
i_{\tQ_2}^* \ldots R(i_{\tQ_t})_* i_{\tQ_t}^* i^* R(i_G)_* \bigr) \circ \mu_K (\BV)
\]
to 
\[
\CH^{n} \bigl( i^*_{Q_1;P_1(\BA_f)gK} R(i_{Q_2})_* 
i_{Q_2}^* \ldots R(i_{Q_r})_* i_{Q_r}^* i^* R(i_G)_* \bigr) \circ \mu_K (\BV)
\]
is identified, under the isomorphisms from~(a) (for $\tQ$ and for $Q$),
with the morphism (respecting the direct sum $\oplus_{p+s = n}$) 
induced by the natural transformation
${}_{\bigcap \utQ} \! \argdot \to {}_{\bigcap \uQ} \! \argdot$.
\end{Thm}

\begin{Proof}
The argument is actually well known (see \emph{e.g.} 
the proofs of \cite[Prop.~(5.2.1)]{P2} or \cite[Thm.~2.9]{BW}): 
since $\bar{Q}_1$ is reductive,
the category $\Rep_F \bar{Q}_1$ is semi-simple. Its derived category
is therefore canonically
equivalent to the category of graded objects in $\Rep_F \bar{Q}_1$.
Therefore,
for any $\BX \in D^+ (\Rep_F Q_1)$, there is a canonical and 
functorial isomorphism in $D^+ (\Rep_F \bar{Q}_1)$
\[
R\Gamma(W_1, \BX) \isoto
                  \bigoplus_{s \in \BZ} H^s (W_1, \BX)[-s] \; .
\]
\end{Proof}

\begin{Rem}
According to Proposition~\ref{7D}~(b), Remark~\ref{7Db}~(b) and 
Proposition~\ref{7N}~(a), the terms 
\[
H^p \bigl( H_C/H_W, {}_{\bigcap \uQ} H^s ( W_1, \Res^G_{Q_1}\BV) \bigr) 
\]
occurring in Theorem~\ref{7Q}, for
$\uQ \in \CC_{(G,\FX) \tei Q_1} \,$,
are identified with
\[
\bigoplus_{\bar{q}_1 \in \Omega}  H^p \Bigl(
    H_C/H_W \cap \bar{q}_1 \bigl( (\bigcap \uQ/W_1)(\BQ) \bigr) \bar{q}_1^{-1}, 
               H^s ( W_1, \Res^G_{Q_1}\BV)  \Bigr) \; ,
\]
where $\Omega$ is a set of representatives of  
$(H_C/H_W) \backslash \bar{Q}_1(\BQ) / ((\bigcap \uQ/W_1)(\BQ))$.
The vector space $H^s ( W_1, \Res^G_{Q_1}\BV)$ underlies an algebraic representation
of $\bar{Q}_1 \,$, hence of its sub-groups 
$C_1/W_1 \cap \bar{q}_1 \bigl( \bigcap \uQ/W_1 \bigr) \bar{q}_1^{-1}$,
$\bar{q}_1 \in \Omega$. It follows as in the proof of
Theorem~\ref{7Q} that each
\[
H^p \Bigl(
    H_C/H_W \cap \bar{q}_1 \bigl( (\bigcap \uQ/W_1)(\BQ) \bigr) \bar{q}_1^{-1}, 
               H^s ( W_1, \Res^G_{Q_1}\BV)  \Bigr) 
\]
admits a canonical decomposition into a direct sum $\oplus_{q+t = p} E^{q,t}_2$, where
\[
E^{q,t}_2 := H^q \Bigl( \overline{H_{C,\bar{q}_1}} , H^t \bigl( 
           \Rad^u( C_1/W_1 \cap \bar{q}_1 \bigl( \bigcap \uQ/W_1 \bigr) \bar{q}_1^{-1} ) ,
               H^s ( W_1, \Res^G_{Q_1}\BV) \bigr) \Bigr) \; .
\]
Here, $\overline{H_{C,\bar{q}_1}}$ denotes the image of the arithmetic sub-group
\[
H_{C,\bar{q}_1} := 
    H_C/H_W \cap \bar{q}_1 \bigl( (\bigcap \uQ/W_1)(\BQ) \bigr) \bar{q}_1^{-1}
\]
of $(C_1/W_1 \cap \bar{q}_1 ( \bigcap \uQ/W_1 ) \bar{q}_1^{-1})(\BQ)$
under the canonical epimorphism 
to the quotient of 
$C_1/W_1 \cap \bar{q}_1 ( \bigcap \uQ/W_1 ) \bar{q}_1^{-1}$
by its unipotent radical, denoted by
$\Rad^u( C_1/W_1 \cap \bar{q}_1 ( \bigcap \uQ/W_1 ) \bar{q}_1^{-1} )$.
\end{Rem}

\begin{Rem}
Thanks to the algebraic nature of the canonical strata $M^K (Q_j,\FX)$, and of $M^K (Q_1;P_1(\BA_f)gK,\FX)$,
the local systems 
\[
\CH^{n} \bigl( i^*_{Q_1;P_1(\BA_f)gK} R(i_{Q_2})_* 
i_{Q_2}^* \ldots R(i_{Q_r})_* i_{Q_r}^* i^* R(i_G)_* \bigr) \circ \mu_K (\BV)
\]
from Theorem~\ref{7Q} underly \emph{variations of Hodge structure}. It appears reasonable to expect
Theorem~\ref{7Q} itself to admit a Hodge theoretic variant, generalizing the main result from \cite{BW}. 
A preliminary step would be to give a Hodge theoretic interpretation of the objects  
\[
\mu_{\pi_1(K_1)} \circ H^p \bigl( H_C/H_W, {}_{\bigcap \uQ} H^s ( W_1, \Res^G_{Q_1}\BV) \bigr) \; .
\]
For the proof of the Hodge theoretic version of Theorem~\ref{7Q}, the ne\-cessary approach will certainly be
very different from the present one, as the Borel--Serre compactification is not algebraic. 

An analogous remark holds for \emph{lisse $\ell$-adic sheaves}, where one may expect a generalization 
of the main result from \cite{P2}.
\end{Rem} 

\begin{Ex} \label{7R}
We consider the Shimura data $(G,\FX)$, whose associated Shimura varieties are
\emph{Siegel threefolds} (see \cite[Ex.~2.7]{P} and \cite[Sect.~1]{W2}).
Fix a four-dimensio\-nal $\BQ$-vector space $V$, together with a $\BQ$-valued
non-degenerate symplectic bilinear form $J$, and define
\[
G := GSp(V,J) \subset GL (V) 
\]
as the group of symplectic similitudes of $V$.
The similitude norm defines a canonical morphism
\[
\lambda: G \longto \Gm \; .
\]
The center $Z(G)$ equals $\Gm \subset GL(V)$
(inclusion of scalar automorphisms).

The analytic space $\FX$ is defined as the sub-space of $M_2 (\BC)$ of those
complex $2 \times 2$-matrices, which are symmetrical, and
whose imaginary part is (positive or negative) definite.
The group of real points $G(\BR)$ acts on $\FX$
by analytical automorphisms. 

In fact,
$(G,\FX)$ \emph{are} pure Shimura data \cite[Ex.~2.7]{P}. 
They satisfy hypo\-the\-sis $(+)$ since $Z(G) = \Gm \,$. 
\\[0.1cm]
(a)~According to \cite[Ex.~4.25]{P},
the proper admissible parabolic sub-groups of $G$
--- equivalently, its maximal proper parabolic sub-groups, since $G^{\ad}$ is simple ---
correspond bijectively to the totally isotropic sub-spaces
of $V$ of strictly positive dimension. More precisely, $V' \subset V$
being a totally isotropic sub-space of dimension one or two, the associated
parabolic equals $\Stab_G(V') \subset G$.
There is a unique boundary component associated to each admissible parabolic
sub-group, meaning that in the notation used so far, we have
\[
\coprod \FX_1/W_1 = \FX_1/W_1 \; .
\]
We fix a symplectic basis $(e_1,e_2,e_3,e_4)$ of $V$,
in which $J$ acquires the $4 \times 4$-matrix
\[
\bigl( \begin{array}{cc}
0 & I_2 \\
-I_2 & 0
\end{array} \bigr) \; ,
\]  
equally denoted by $J$
($I_2:=$ the $2 \times 2$-matrix representing the identity). We use this basis
to identify $G$ with the sub-group $GSp_4$ of 
$GL_4$ of matrices $g$
satisfying the relation
\[
^{t}\! g J g = \lambda(g) \cdot J \; .
\]
The sub-spaces $V_1'$ and $V_2'$ generated
by $\{ e_1,e_2 \}$ and by $\{ e_1 \}$, respectively, are both totally isotropic. 
Put $Q_m := \Stab_G(V_m')$, $m = 1,2$. As usual,
$P_m$ denotes the canonical normal sub-group of $Q_m$ defined in \cite[Sect.~4.7]{P}, 
and $W_m$ its unipotent radical
(which equals the unipotent radical of $Q_m$), $m = 1,2$.
Then in our setting, still according to \cite[Ex.~4.25]{P}, 
\[
Q_1 = \biggl\{ \bigl( \begin{array}{cc}
q \cdot A & A \cdot M \\
0 & {}^t A ^{-1}
\end{array} \bigr) \; , \; q \in \Gm \, , A \in GL_{2,\BQ} \, , {}^t M = M \biggr\} \; ,
\] 
\[
P_1 = \biggl\{ \bigl( \begin{array}{cc}
q \cdot I_2 & M \\
0 & I_2
\end{array} \bigr) \; , \; q \in \Gm \, , {}^t M = M \biggr\} \; ,
\] 
\[
W_1 = \biggl\{ \bigl( \begin{array}{cc}
I_2 & M \\
0 & I_2
\end{array} \bigr) \; , \; {}^t M = M \biggr\} \; ,
\] 
while $Q_2$ equals
\[
\Biggl\{ \bigl( \begin{array}{cccc}
a & aq^{-1} (bu+dw) & v & aq^{-1} (cu+ew) \\
0 & b & w & c \\
0 & 0 & a^{-1}q & 0 \\
0 & d & -u & e
\end{array} \bigr) \; , \; a \, , be - cd = q \in \Gm \Biggr\} \; ,
\] 
\[
P_2 = \Biggl\{ \bigl( \begin{array}{cccc}
be - cd & bu+dw & v & cu+ew \\
0 & b & w & c \\
0 & 0 & 1 & 0 \\
0 & d & -u & e
\end{array} \bigr) \; , \; be - cd \in \Gm \Biggr\} 
\] 
and
\[
W_2 = \Biggl\{ \bigl( \begin{array}{cccc}
1 & u & v & w \\
0 & 1 & w & 0 \\
0 & 0 & 1 & 0 \\
0 & 0 & -u & 1
\end{array} \bigr) \Biggr\} 
\] 
(see \cite[pp.~539--540]{W2}, noting that in \loccit, the groups are indexed
by $0$ and $1$ instead of $1$ and $2$, according to the dimension of the
associated strata $M^K (Q_m,\FX)$).
Observe that 
\[
Q_1 \cap Q_2 = \Biggl\{ \bigl( \begin{array}{cccc}
a & aq^{-1} bu & v & a(q^{-1}cu+ b^{-1}w) \\
0 & b & w & c \\
0 & 0 & a^{-1}q & 0 \\
0 & 0 & -u & b^{-1}q
\end{array} \bigr) \; , \; a , b , q \in \Gm \Biggr\} 
\] 
equals the Borel sub-group of $G$ 
stabilizing the flag $V_2' \subset V_1'$. As for the
quotients by $W_1$, we have
\[
\bar{Q}_1 = Q_1/W_1 = P_1/W_1 \times_\BQ GL_2 = \Gm \times_\BQ GL_2 \; ,
\]
the identification given by sending the class of a matrix
\[
\bigl( \begin{array}{cc}
q \cdot A & A \cdot M \\
0 & {}^t A ^{-1}
\end{array} \bigr) 
\] 
to the pair $(q,A)$. Under this identification,
\[
(Q_1 \cap Q_2)/W_1 = \Gm \times_\BQ B \; ,
\]
where $B \subset GL_2$ is the Borel sub-group of upper triangular matrices.

This identification also shows 
that the sub-group $C_1$ of $Q_1$ (Definition~\ref{7A}) is $Q_1$ itself. \\[0.1cm]
(b)~The (unique) boundary component $(P_1',\FX_1')$ of the Shimura data
$(P_2,\FX_2)$ associated to the proper parabolic sub-group
\[
Q_1 \cap P_2 = \Biggl\{ \bigl( \begin{array}{cccc}
q & bu & v & cu+ b^{-1}qw \\
0 & b & w & c \\
0 & 0 & 1 & 0 \\
0 & 0 & -u & b^{-1}q
\end{array} \bigr) \; , \; b , q \in \Gm \Biggr\} 
\]
equals $(P_1,\FX_1)$: indeed, let $Q_1'$ denote the proper parabolic sub-group of $G$ ,
whose associated boundary component is $(P_1',\FX_1')$. 
According to \cite[Lemma~4.19]{P}, we have
$Q_1 \cap P_2 = Q_1' \cap P_2$ and
\[
Q_2 = (Q_1' \cap Q_2)P_2 \; .
 \]
But from the explicit description from (a), we also have
\[
Q_2 = (Q_1 \cap Q_2)P_2 \; .
 \]
It follows that $Q_1 \cap Q_2 = Q_1' \cap Q_2 \,$.
Since the maximal proper sub-groups containing a given parabolic sub-group
are unique, we conclude that $Q_1'$ is equal to $Q_1 \,$.
Consequently, we have $(P_1',\FX_1') = (P_1,\FX_1)$. \\[0.1cm]
(c)~According to (b), the sequence
$\uQ := (Q_1 \prec Q_2)$ belongs to $\CC_{(G,\FX) \tei Q_1} \,$. We have 
\[
\bigcap \uQ = Q_1 \cap Q_2 
\]
and $\uQ \subset \utQ := (Q_1)$.

According to part~(a) of Theorem~\ref{7Q}, 
putting $\Gamma_g := H_C/H_W$, there are canonical and functorial 
isomorphisms between
\[
\CH^{n} \bigl( i^*_{Q_1;P_1(\BA_f)gK} i^* R (i_G)_* \bigr) \circ \mu_K (\BV)
\]
and
\[
\bigoplus_{p+s = n} \mu_{\pi_1(K_1)} \circ 
     H^p \bigl( \Gamma_g, H^s ( W_1, \Res^G_{Q_1}\BV) \bigr) \; ,
\]
and between
\[
\CH^{n} \bigl( i^*_{Q_1;P_1(\BA_f)gK} R(i_{Q_2})_* 
i_{Q_2}^* i^* R(i_G)_* \bigr) \circ \mu_K (\BV)
\]
and 
\[
\bigoplus_{p+s = n} \mu_{\pi_1(K_1)} \circ 
   H^p \bigl( \Gamma_g, {}_{(Q_1 \cap Q_2)} H^s ( W_1, \Res^G_{Q_1}\BV) \bigr) \; ,
\]
for any $n \in \BZ$, and any $V \in D^+ ( \Rep_F G )$.
Note that $\Gamma_g$ is a neat arithmetic sub-group of 
$\bar{Q}_1$ (according to Proposition~\ref{7N}~(a)), \emph{i.e.}, of 
$\Gm \times_\BQ GL_2 \,$, and hence of $SL_2 \,$.
According to part~(b) of Theorem~\ref{7Q},
the adjunction 
\[
\CH^\ast \bigl( i^*_{Q_1;P_1(\BA_f)gK} i^* R(i_G)_* \bigr) 
\longto
\CH^\ast \bigl( i^*_{Q_1;P_1(\BA_f)gK} R(i_{Q_2})_* 
i_{Q_2}^* i^* R(i_G)_* \bigr) 
\]
on $\mu_K (\BV)$ is induced, under the above isomorphisms, 
by the morphism 
\[
H^\ast ( W_1, \Res^G_{Q_1}\BV)
\longto {}_{(Q_1 \cap Q_2)} H^\ast ( W_1, \Res^G_{Q_1}\BV)
\]
from Construction~\ref{7Con}. \\[0.1cm]
(d)~According to Proposition~\ref{7D}~(b) and Remark~\ref{7Db}~(b),
for all integers $p$ and $s$,
\[
H^p \bigl( \Gamma_g \,, {}_{(Q_1 \cap Q_2)} H^s ( W_1, \Res^G_{Q_1}\BV) \bigr)
\]
can be identified with
\[
\bigoplus_{q_1 \in \Omega} 
H^p \bigl( \Gamma_g \cap q_1 B'(\BQ) q_1^{-1}, H^s ( W_1, \Res^G_{Q_1}\BV) \bigr) \; ,
\]
where $B' := SL_2 \cap B$, and $\Omega$ is a set of representatives 
in $SL_2(\BQ)$
of the set of \emph{cusps}
\[
\Gamma_g \backslash \bar{Q}_1(\BQ) / (((Q_1 \cap Q_2)/W_1)(\BQ))
= \Gamma_g \backslash SL_2(\BQ) / B'(\BQ)
\]
of $\Gamma_g \,$. 
By Proposition~\ref{7D}~(c), the composition of
the morphism
\[
H^p \bigl( \Gamma_g \,, H^s ( W_1, \Res^G_{Q_1}\BV) \bigr)
\longto 
H^p \bigl( \Gamma_g \,, {}_{(Q_1 \cap Q_2)} H^s ( W_1, \Res^G_{Q_1}\BV) \bigr)
\]
with this identification equals the direct sum
\[
H^p \bigl( \Gamma_g \,, H^s ( W_1, \Res^G_{Q_1}\BV) \bigr)
\longto 
\bigoplus_{q_1 \in \Omega} 
H^p \bigl( \Gamma_g \cap q_1 B'(\BQ) q_1^{-1}, H^s ( W_1, \Res^G_{Q_1}\BV) \bigr) 
\]
of the restrictions from $\Gamma_g$ to $\Gamma_g \cap q_1 B'(\BQ) q_1^{-1}$,
for $q_1 \in \Omega$. 
\end{Ex}

\begin{Rem} \label{7S}
(a)~Parts~(c) and (d) of Example~\ref{7R} 
yield a positive reply to the question raised in \cite[Rem.~2.10~(b)]{W2}. \\[0.1cm]
(b)~In the general situation considered in the section, and
for elements of $\CC_{(G,\FX) \tei Q_1}$ of the form
$\uQ = (Q_1 \prec Q_2)$, the information
provided by Theorem~\ref{7Q} is formally identical to Example~\ref{7R}~(c): 
putting $\Gamma_g := H_C/H_W$, there are canonical and functorial 
isomorphisms between
\[
\CH^{n} \bigl( i^*_{Q_1;P_1(\BA_f)gK} i^* R(i_G)_* \bigr) \circ \mu_K (\BV)
\]
and
\[
\bigoplus_{p+s = n} \mu_{\pi_1(K_1)} \circ 
     H^p \bigl( \Gamma_g, H^s ( W_1, \Res^G_{Q_1}\BV) \bigr) \; ,
\]
and between
\[
\CH^{n} \bigl( i^*_{Q_1;P_1(\BA_f)gK} R(i_{Q_2})_* 
i_{Q_2}^* i^* R(i_G)_* \bigr) \circ \mu_K (\BV)
\]
and 
\[
\bigoplus_{p+s = n} \mu_{\pi_1(K_1)} \circ 
   H^p \bigl( \Gamma_g, {}_{(Q_1 \cap Q_2)} H^s ( W_1, \Res^G_{Q_1}\BV) \bigr) \; ,
\]
for any $n \in \BZ$, and any $V \in D^+ ( \Rep_F G )$. 
The adjunction 
\[
\CH^\ast \bigl( i^*_{Q_1;P_1(\BA_f)gK} i^* R(i_G)_* \bigr) 
\longto
\CH^\ast \bigl( i^*_{Q_1;P_1(\BA_f)gK} R(i_{Q_2})_* 
i_{Q_2}^* i^* R(i_G)_* \bigr)
\]
on $\mu_K (\BV)$ is induced by the morphism 
\[
H^\ast ( W_1, \Res^G_{Q_1}\BV)
\longto {}_{(Q_1 \cap Q_2)} H^\ast ( W_1, \Res^G_{Q_1}\BV)
\]
from Construction~\ref{7Con}. 

The only formal difference to Example~\ref{7R}~(c) is that
\emph{a priori}, the
group $\Gamma_g$ is neat arithmetic in $C_1/W_1 \,$, and in general not in
$\bar{Q}_1 \,$. \\[0.1cm]
(c)~The situation from (b) occurs in the context of 
\emph{genus $2$ Hilbert--Siegel varieties} (of which Siegel threefolds
are a special case), considered in \cite{C}. More precisely,
the proof of \cite[Prop.~2.5.1.4]{C} (which generalizes \cite[Prop.~2.9]{W2})
necessitates the non-vanishing of the kernel of a certain map denoted $ad$.
In \loccit, this is achieved by showing that the dimension
of the target is strictly smaller than that of the source. 
 
Part~(b) above allows to give an alternative proof of this non-vanishing
in that it relates $ad$ to the residue map on Hilbert
modular forms. Therefore, its kernel contains the cusp forms.
\end{Rem}


\bigskip

%
%

\section{Reformulation in terms of group cohomology. II}
\label{8}



The purpose of this section is to translate Theorem~\ref{6MT2}
(through Corollary~\ref{6D})
into group cohomology (Main Theorem~\ref{8MT} and its Variants~\ref{8Var1} and \ref{8Var2}). 
We keep the setting of Section~\ref{7}: $(P,\FX) = (G,\FX)$ are
pure Shimura data satisfying hypothesis 
$(+)$, and $K$ is a neat open compact sub-group $G(\BA_f)$. 
Furthermore, $Q_1 \ne G$ is
an admissible parabolic sub-group of $G$, $g \in G (\BA_f)$,
\[
H_1= H_1(gK)= Q_1 (\BQ) \cap gKg^{-1} \; ,
\]
and 
\[
H_C= H_C(gK)= C_1 (\BQ) \cap gKg^{-1} \; ,
\]
where
\[ 
C_1 = \{ q \in Q_1 \; , \; \pi_{Q_1}(q) \in \Cent_{\bar{Q}_1}(\pi_{Q_1}(P_1)) \}^0 
\]
(Definition~\ref{7A}). \\

Let us start to set up the cohomological data necessary for the statement of 
Main Theorem~\ref{8MT}. Apart from group cohomology (Definition~\ref{7C}~(a)),
we shall need \emph{cohomology with compact supports} and \emph{boundary 
cohomology} of $H_C \,$. \\

The foundations for these cohomology theories were laid in \cite{W3}.
We consider a
space $X$ underlying a space of type $S - \BQ$ under $C_1$ \cite[Def.~2.3]{BS},
and the cohomology with compact support and boundary cohomology of the quotient $H_C \backslash X$
with coefficients in local systems associated to representations $\BV$ of $H_C$.
Actually, as in the previous section,
we shall need the equivariant version of the theory, with respect to the action of the
quotient $H_1/H_C$, when $\BV$ carries not only an action of $H_C$, but of $H_1$. 
We are therefore led to consider a situation where $X$ occurs as a fibre of a certain family
of spaces. \\

Thus, let $Y$ be a space underlying a space of type $S - \BQ$ under $Q_1$ \cite[Def.~2.3]{BS},
meaning that $Y$ is a homogeneous space under $Q_1(\BR)$, such that the stabilizers in $Q_1(\BR)$
of the points $y \in Y$ are of the form $K_y \cdot S_y(\BR)$, where $S_y$ is a maximal torus of the base change
to $\BR$ of the $\BQ$-split radical $R_d Q_1$ of $Q_1 \,$, and $K_y$ is a maximal compact sub-group of
$Q_1(\BR)$ normalizing (actually \cite[Rem.~2.2]{BS}, centralizing) $S_y \,$.
Define $Z := C_1(\BR) \backslash Y$, and denote by $\pi: Y \onto Z$ the canonical surjection. 

\begin{Prop}[{\cite{W3}}] \label{8a}
(a)~The space $Z$ is contractible.\\[0.1cm]
(b)~Each fibre of $\pi$ underlies a space of type $S - \BQ$ under $C_1$. \\[0.1cm]
(c)~The map $\pi$ can be trivialized $C_1(\BR)$-equivariantly: letting $X$ denote any of the fibres of $\pi$, 
there is a homeomorphism
\[
\Phi: X \times Z \isoto Y 
\]
satisfying
\begin{enumerate}
\item[(i)] for all $(x,z) \in X \times Z$ and $q \in C_1(\BR)$, we have
\[
\Phi (qx,z) = q \Phi (x,z) \; , 
\] 
\item[(ii)] the diagram
\[
\vcenter{\xymatrix@R-10pt{
X \times Z \ar[rr]^-{\Phi} \ar@{->>}[dr]_-{p^X} && 
Y \ar@{->>}[dl]^-{\pi} \\
& Z  &
\\}}
\]
is commutative,
where $p^X$ denotes the projection onto the second factor.
\end{enumerate}
\end{Prop}

The action of $Q_1(\BR)$ on $Y$ induces an action on $Z$, which restricts to the trivial action of $C_1(\BR)$. 
The canonical surjection $\pi: Y \onto Z$ is $Q_1(\BR)$-equivariant.
On $Z$, there is a notion of $Q_1(\BR)/C_1(\BR)$-, and hence of $H_1/H_C$-sheaves of Abelian groups.
The latter notion is also defined on the intermediate quotient $H_C \backslash Y$,
which is related to $Y$ and $Z$ \emph{via} a canonical factorization
\[
Y \longonto H_C \backslash Y \stackrel{\tilde{\pi}}{\longonto} Z
\]
of $\pi$. Proposition~\ref{8a}~(c) implies that $\tilde{\pi}: H_C \backslash Y \onto Z$ is a trivial fibration.
For $\argast \in \{ H_C \backslash Y , Z \}$, let
us denote by $H_1/H_C \, \text{-} \Loc ( \argast )$ the category of $H_1/H_C$-local systems 
on $\argast$, and by $D^+ ( H_1/H_C \, \text{-} \argast )$ the 
derived category of complexes of $H_1/H_C$-sheaves on $\argast$,
that are bounded from below. \\

Since $K$ is supposed neat, 
the arithmetic sub-group $H_1$ of $Q_1 (\BQ)$ acts freely on $Y$ \cite[Sect.~9.5]{BS},
and hence so does $H_C \,$. We thus get a functor
\[
\eta_g: \Rep ( H_1 ) \longto 
H_1/H_C \, \text{-} \Loc \bigl( H_C \backslash Y \bigr) \; ,
\]
by associating to each representation $\BV$ the sheaf of continuous sections of the quotient map
$Y \onto H_C \backslash Y$, equipped with the $H_1/H_C$-action coming from the actions of $H_1$ on $H_C \backslash Y$ and on $\BV$. 
The functor $\eta_g$ is exact, and therefore derives trivially,
to give a functor denoted by the same symbol
\[
\eta_g: D^+ ( \Rep H_1 ) \longto D^+ \bigl( H_1/H_C \, \text{-} (H_C \backslash Y) \bigr) \; .
\]

\begin{Def}[{\cite{W3}}] \label{8A}
Fix a space $Y$ of type $S - \BQ$ under $Q_1 \,$. \\[0.1cm]
(a)~Define
\[
R \Gamma_c (H_C, {\argdot}): D^+ (\Rep H_1) \longto D^+ \bigl( \Rep (H_1/H_C) \bigr)
\]
as being equal to $R \Gamma (Z, R \tilde{\pi}_! (\eta_g ({\argdot})))$. \\[0.1cm]
(b)~Define
\[
\partial R \Gamma (H_C, {\argdot}): D^+ (\Rep H_1) \longto D^+ \bigl( \Rep (H_1/H_C) \bigr)
\]
as being equal to $R \Gamma (Z, \partial R \tilde{\pi}_* (\eta_g ({\argdot})))$.
\end{Def}

Some words of explanation are in order. The space $Z$ being contractible (Proposition~\ref{8a}~(a)),
local systems on $Z$ are necessarily constant. Furthermore, the
global section functor $\Gamma (Z, \argast)$ induces an equivalence of categories between 
$H_1/H_C \, \text{-} \Loc ( Z )$ and $\Rep (H_1/H_C)$.
Similarly,
its right derived functor $R \Gamma (Z, \argast)$ induces an equivalence between
the full sub-category of $D^+ ( H_1/H_C \, \text{-} Z )$ of classes, 
whose cohomology objects are local systems, 
and $D^+ ( \Rep (H_1/H_C) )$. Since $\tilde{\pi}$ is a trivial fibration, the functors 
$R \tilde{\pi}_*$ and $R \tilde{\pi}_!$
map $H_1/H_C$-local systems on $H_C \backslash Y$ to 
objects of $D^+ ( H_1/H_C \, \text{-} Z )$, whose cohomology objects are local systems.
Hence so does $\partial R \tilde{\pi}_*$ (defined as the right derived functor of the left exact functor
$R^0 \tilde{\pi}_* / R^0 \tilde{\pi}_!$ on the category of $H_1/H_C$-sheaves on $H_C \backslash Y$). 
We refer to \cite{W3} for details.

\begin{Cons} \label{8Aa}
Keep the setting of Definition~\ref{8A}. 
Applying the functor $R \Gamma (Z, \argast)$ to $E \eta_g ({\argdot})$, where $E$ is
the exact triangle 
\[
R \tilde{\pi}_! \longto R \tilde{\pi}_* \longto \partial R \tilde{\pi}_* \stackrel{[1]}{\longto}                    
\]
\cite{W3}, we get an exact triangle 
\[
R \Gamma_c (H_C, {\argdot}) \longto R \Gamma (H_C, {\argdot}) 
\longto \partial R \Gamma (H_C, {\argdot}) 
\stackrel{[1]}{\longto}                    
\]
of functors $D^+ (\Rep H_1) \to D^+ ( \Rep (H_1/H_C) )$.
Note that $R \Gamma (H_C, {\argdot})$ is indeed identified with 
\[
R \Gamma \bigl( Z, R \tilde{\pi}_* (\eta_g ({\argdot})) \bigr)
= R \Gamma \bigl( H_C \backslash Y, \eta_g ({\argdot}) \bigr)
\]
as the action of $H_C$ on $Y$ is free, and $Y$ is contractible \cite[Rem.~2.4~(1)]{BS}. 
\end{Cons}

\forget{
\begin{Rem} \label{8Ab}
In the previous setting, the action
of $H_C$ on $X^{BS}$ remains free \cite[Sect.~9.5]{BS}. The immersion of $X$ into $X^{BS}$ being
contractible, the same is therefore true for $\j \,$. By Corollary~\ref{1C}~(c), for any $\BV \in D^+ (\Rep H_C)$, the object 
$R \j_* \eta_g (\BV)$ of $D^+(H_C \backslash X^{BS})$ is represented by a complex of local systems.
More precisely, it equals $\eta_g^{BS} (\BV)$, where 
\[
\eta_g^{BS}: D^+ ( \Rep H_C ) \longto D^+ \bigl( H_C \backslash X^{BS} \bigr) 
\]
is defined in complete analogy to $\eta_g$. Thus, the
exact triangle 
\[
R \Gamma_c (H_C, {\argdot}) \longto R \Gamma (H_C, {\argdot}) 
\longto \partial R \Gamma (H_C, {\argdot}) 
\stackrel{[1]}{\longto}                    
\]
from Construction~\ref{8Aa} equals the one obtained by applying $R \Gamma (H_C \backslash X^{BS}, \argast)$ 
to $T \eta_g^{BS} ({\argdot})$. 
\end{Rem}
}
\begin{Prop}[{\cite{W3}}] \label{8B}
(a)~There is a canonical exact triangle 
\[
R \Gamma_c (H_C, {\argdot}) \longto R \Gamma (H_C, {\argdot}) 
\longto \partial R \Gamma (H_C, {\argdot}) 
\stackrel{[1]}{\longto}                    
\]
of functors from $D^+ (\Rep H_1)$ to $D^+ ( \Rep (H_1/H_C) )$. \\[0.1cm]
(b)~The functors $R \Gamma_c (H_C, {\argdot})$ and $\partial R \Gamma (H_C, {\argdot})$,
and the exact triangle of functors from (a) are all independent of the choice of $Y$.
\end{Prop}

\begin{Proof}
In order to prove that the result of Construction~\ref{8Aa} does not depend on $Y$,
let $Y'$ be another choice of space underlying a space of type $S - \BQ$ under $Q_1 \,$. 
Both $Y$ and $Y'$ being homogeneous under $Q_1(\BR)$, and the collection of stabilizers of points
constituting (only) one conjugation class of sub-groups of $Q_1(\BR)$ \cite[Lemma~2.1~(iii)]{BS},
there is a $Q_1(\BR)$-equivariant homeomorphism $\psi: Y \isoto Y'$. 
Any two
choices of $\psi$ can be interpolated by $Q_1(\BR)$-equivariant homeomorphisms \cite{W3}.    
\end{Proof}

\begin{Rem} \label{8C}
(a)~In general, the isomorphism $\psi$ occurring in the proof of Proposition~\ref{8B}~(b) is not unique. \\[0.1cm]
(b)~The author knows of no purely algebraic definition of $R \Gamma_c (H_C, {\argdot})$ and $\partial R \Gamma (H_C, {\argdot})$. 
\\[0.1cm]
(c)~Proposition~\ref{8a}--Proposition~\ref{8B} generalize from $H_C$ to any neat arithmetic sub-group of 
the rational points of a connected affine algebraic group $C_1$ over $\BQ \,$,
and from $H_1$ to any sub-group of $Q_1(\BR)$ normalizing $H_C \,$, and such that $H_C$ is contained in
$C_1 (\BQ) \cap H_1$ \cite{W3}.  This applies in particular to the quotients $H_C/H_W$ and $H_1/H_W$, where
\[
H_W = H_W(gK) = W_1(\BQ) \cap gKg^{-1} 
\]
(cmp.~Section~\ref{7}). We thus get a canonical exact triangle 
\[
R \Gamma_c (H_C/H_W, {\argdot}) \longto R \Gamma (H_C/H_W, {\argdot}) 
\longto \partial R \Gamma (H_C/H_W, {\argdot}) 
\stackrel{[1]}{\longto}                    
\]
of functors from $D^+ (\Rep H_1/H_W)$ to $D^+ ( \Rep (H_1/H_C) )$. 
The latter exact triangle is related to the one from Proposition~\ref{8B}
by a canonical isomorphism of functors on $D^+ (\Rep H_1)$ between
\[
R \Gamma_c (H_C, {\argdot}) \longto R \Gamma (H_C, {\argdot}) 
\longto \partial R \Gamma (H_C, {\argdot}) 
\stackrel{[1]}{\longto}                    
\]
and
\[
R \Gamma_c (H_C/H_W, \argast) \longto R \Gamma (H_C/H_W, \argast) 
\longto \partial R \Gamma (H_C/H_W, \argast) 
\stackrel{[1]}{\longto} \; ,                   
\]
evaluated at $\argast = R \Gamma ( H_W, {\argdot} )$ (see \cite{W3}). \\[0.1cm]
(d)~According to Proposition~\ref{8a}~(a), (b), application of the forgetful functor
\[
D^+ \bigl( \Rep (H_1/H_C) \bigr) \longto D^+ (\Ab)
\]
transforms the exact triangle
\[
R \Gamma_c (H_C, {\argdot}) \longto R \Gamma (H_C, {\argdot}) 
\longto \partial R \Gamma (H_C, {\argdot}) 
\stackrel{[1]}{\longto}                    
\]
into the analogous exact triangle of functors $D^+ (\Rep H_C) \to D^+ ( \Ab )$,
preceded by the forgetful functor $D^+ (\Rep H_1) \to D^+ (\Rep H_C)$.
\end{Rem}

Recall that $i_G$ denotes the open immersion of $M^K  (G,\FX)$ into $M^K  (G,\FX)^*$, and
$i$ the complementary immersion of $\partial M^K (G,\FX)^*$.
Furthermore, $i_{Q_1}$ denotes the immersion of $M^K ( Q_1,\FX )$ into $M^K  (G,\FX)^*$ ---
by slight abuse of notation, we agreed to use the same symbol for the immersion of $M^K ( Q_1,\FX )$ into 
$\partial M^K  (G,\FX)^*$ as well --- and
\[
h_{Q_1}: \partial M^K  (G,\FX)^* - \overline{M^K \bigl( Q_1,\FX \bigr)} 
\longinto \partial M^K  (G,\FX)^*
\]
the open immersion of the complement of $\overline{M^K \bigl( Q_1,\FX \bigr)}$. The functor
\[
\bar{\mu}: \Rep(H_1/H_C) \longto \Loc \bigl( pr_{gK} \bigl(\coprod \FX_1/W_1 \bigr) \bigr) 
\] 
associates to each representation $\BV$
the sheaf of continuous sections of
\[
(H_1/H_C) \backslash \bigl( \BV \times \coprod \FX_1/W_1  \bigr) 
\longonto (H_1/H_C) \backslash \bigl( \coprod \FX_1/W_1  \bigr)
= pr_{gK} \bigl(\coprod \FX_1/W_1 \bigr) 
\]
(see Main Theorem~\ref{7MT}). It induces an exact functor
\[
\bar{\mu}: D^+ \bigl( \Rep(H_1/H_C) \bigr) \longto 
D^+ \bigl( pr_{gK} \bigl(\coprod \FX_1/W_1 \bigr) \bigr) \; .
\]
For $\BV \in D^+ ( \Rep (G(\BQ) ) )$,
Main Theorem~\ref{7MT}~(a) (see Remark~\ref{7Rem}~(a)) yields a canonical isomorphism between
the restriction to $pr_{gK}(\coprod \FX_1/W_1) \subset M^K ( Q_1,\FX )(\BC)$ of  
$i_{Q_1}^* R(i_G)_* \mu_K(\BV)$ and the image under the functor $\bar{\mu}$ of the
object $R \Gamma (H_C,\Res^{G(\BQ)}_{H_1}(\BV))$ of $D^+ (\Rep (H_1/H_C))$.
Main Theorem~\ref{8MT} will extend this isomorphism to the exact triangle
\[
R i_{Q_1}^! i^*R(i_G)_* \mu_K(\BV) \to  
i_{Q_1}^* R(i_G)_* \mu_K(\BV) \to
i_{Q_1}^* R (h_{Q_1})_* h_{Q_1}^* i^*R(i_G)_* \mu_K(\BV) \stackrel{[1]}{\to} 
\]  

\begin{MThm} \label{8MT}
Let $\BV \in D^+ ( \Rep (G(\BQ) ) )$. \\[0.1cm]
(a)~Assume that $Q_1$ is maximal proper. 
Then the isomorphism
from Theorem~\ref{6MT2}~(a) induces 
a canonical isomorphism of exact triangles in the derived category $D^+ ( pr_{gK} (\coprod \FX_1/W_1 ))$
between the image under $\bar{\mu}$ of the exact triangle  
\[
R \Gamma_c \bigl( H_C, \Res^{G(\BQ)}_{H_1}(\BV) \bigr) \to R \Gamma \bigl( H_C, \Res^{G(\BQ)}_{H_1}(\BV) \bigr) 
\to \partial R \Gamma \bigl( H_C, \Res^{G(\BQ)}_{H_1}(\BV) \bigr) \stackrel{[1]}{\to}           
\] 
(Proposition~\ref{8B}~(a)) and the restriction to $pr_{gK}(\coprod \FX_1/W_1)$ of the exact triangle
\[
R i_{Q_1}^! i^*R(i_G)_* \mu_K(\BV) \to  
i_{Q_1}^* R(i_G)_* \mu_K(\BV) \to
i_{Q_1}^* R (h_{Q_1})_* h_{Q_1}^* i^*R(i_G)_* \mu_K(\BV) \stackrel{[1]}{\to}  
\]  
in $D^+ (M^K ( Q_1,\FX )(\BC))$. \\[0.1cm]
(b)~If $Q_1$ is the intersection of $r$ distinct maximal proper sub-groups of $G$, then
the restriction to $pr_{gK}(\coprod \FX_1/W_1)$
of $R i_{Q_1}^! i^*R(i_G)_* \mu_K(\BV)$ is isomorphic to  
\[
\bar{\mu} \Bigl( R \Gamma_c \bigl( H_C, \Res^{G(\BQ)}_{H_1}(\BV) \bigr) \Bigr) [-(r-1)] \; .
\]
(c)~If $Q_1$ is not maximal proper, then the exact triangle
\[
R i_{Q_1}^! i^*R(i_G)_* \mu_K(\BV) \stackrel{\alpha}{\to}  
i_{Q_1}^* R(i_G)_* \mu_K(\BV) \to
i_{Q_1}^* R (h_{Q_1})_* h_{Q_1}^* i^*R(i_G)_* \mu_K(\BV) \stackrel{[1]}{\to} 
\]  
is split: the morphism $\alpha$ is zero. \\[0.1cm]
(d)~If $Q_1$ is (not maximal proper, and equal to)
the intersection of $r$ distinct maximal proper sub-groups of $G$, 
with $r \ge 2$, then there is an isomorphism from the restriction to $pr_{gK}(\coprod \FX_1/W_1)$
of $i_{Q_1}^* R (h_{Q_1})_* h_{Q_1}^* i^* R(i_G)_* \mu_K(\BV)$ to 
\[ 
\bar{\mu} \Bigl( R \Gamma \bigl( H_C, \Res^{G(\BQ)}_{H_1}(\BV) \bigr) \Bigr)
\oplus \bar{\mu} \Bigl( R \Gamma_c \bigl( H_C, \Res^{G(\BQ)}_{H_1}(\BV) \bigr) \Bigr) [-(r-2)] \; ,
\]
identifying the adjunction
\[
i_{Q_1}^* R(i_G)_* \mu_K(\BV) \longto i_{Q_1}^* R (h_{Q_1})_* h_{Q_1}^* i^*R(i_G)_* \mu_K(\BV)
\]
with the inclusion of the first component (see Remark~\ref{7Rem}~(a)), and the boundary
\[
i_{Q_1}^* R (h_{Q_1})_* h_{Q_1}^* R(i_G)_* \mu_K(\BV) 
\longto R i_{Q_1}^! i^*R(i_G)_* \mu_K(\BV) [1] 
\]
with the projection onto the second component (under the isomorphism of (b)). \\[0.1cm]
(e)~All objects in the exact triangle
\[
R i_{Q_1}^! i^*R(i_G)_* \mu_K(\BV) \to  
i_{Q_1}^* R(i_G)_* \mu_K(\BV) \to
i_{Q_1}^* R (h_{Q_1})_* h_{Q_1}^* i^*R(i_G)_* \mu_K(\BV) \stackrel{[1]}{\to} 
\] 
in $D^+(M^K ( Q_1,\FX )(\BC))$ can be represented by complexes of local systems.
\end{MThm}

\forget{
\begin{Rem}
It would be desirable to have comparison statements concerning objects in $D^+(M^K ( Q_1,\FX )(\BC))$
or $D^+(pr_{gK} ( \coprod \FX_1/W_1 ))$ (as in the previous section), rather than their pull-backs to $\coprod \FX_1/W_1$.   
In order to obtain such statements, it would be necessary to define the action of the quotient $H_1/H_C$ on 
$R \Gamma_c (H_C, \BV)$ and $\partial R \Gamma (H_C, \BV)$, for objects $\BV$ of $D^+ (\Rep H_1)$ (cmp.~Remark~\ref{8C}~(b)).
\end{Rem}

\medskip
}
Main Theorem~\ref{8MT} will be proved together with
the translation of Proposition~\ref{6F} into group cohomology
(Proposition~\ref{8R}), whose statement we prepare now.

\begin{Cons} \label{8P}
Fix a space $(Y,(L_y)_{y \in Y})$ of type $S - \BQ$ under $Q_1 \,$.
By definition \cite[Def.~2.3]{BS}, this means that in addition to the specific type
of homogeneous space $Y$ (see the beginning of this section), 
a family $(L_y)_{y \in Y}$ of Levi sub-groups of   
$Q_{1,\BR}$ is given, such that $L_y(\BR)$ contains the stabilizer of $y$, 
and $L_{qy} = q L_y q^{-1}$, for every $y \in Y$ and $q \in Q_1(\BR)$.

Using this choice $(L_y)_{y \in Y}$, the manifold with corners $Y^{BS}$
is constructed in \cite[Sect.~7.1]{BS} (where the notation $\bar{Y}$ is used instead of $Y^{BS}$).
The set $Y^{BS}$ is the disjoint union 
\[
Y^{BS} = \coprod_Q e(Q) 
\]
of faces $e(Q)$, one for each parabolic sub-group $Q$ of $Q_1 \,$. The space $Y$ is identified with
the (open) face $e(Q_1)$, and the action of the sub-group $Q_1(\BQ)$ of $Q_1(\BR)$ is extended to an
action on $Y^{BS}$. Consider
\[
Y^{BS \, '} := \coprod_{\adm(Q) = Q_1} e(Q) \subset Y^{BS} \; . 
\] 
The disjoint union extends only over those parabolic sub-groups $Q$ satisfying 
$\adm(Q) = Q_1 \,$, \emph{i.e.}, $P_1 \subset Q \subset Q_1 \,$. Given the topology of $Y^{BS}$
\cite[Sect.~7.1 and 5.1]{BS}, the sub-set $Y^{BS \, '}$ of $Y^{BS}$ is open. Also, it is stable
under the action of $Q_1(\BQ)$, thanks to the formula $q e(Q) = e(qQq^{-1})$, for all parabolics $Q$
and $q \in Q_1(\BQ)$ \cite[Sect.~5.6]{BS}. 

We then have \cite{W3}: (a)~the canonical surjection $\pi: Y \onto Z = C_1(\BR) \backslash Y$ 
admits a unique continuous extension
$\pi': Y^{BS \, '} \onto Z$ to $Y^{BS \, '}$ (necessarily $Q_1(\BQ)$-equivariant), 
(b)~for each $z \in Z$, the fibre $\pi^{-1}(z)$ inherits from
$(Y,(L_y)_{y \in Y})$ the structure of a space of type $S - \BQ$ under $C_1 \,$, and the fibre $(\pi')^{-1}(z)$
is canonically identified with the manifold with corners $(\pi^{-1}(z))^{BS}$ from \cite[Sect.~7.1]{BS}, 
(c)~the map $\pi'$ can be trivialized $C_1(\BQ)$-equivariantly. More precisely, a $C_1(\BR)$-equivariant trivialization
of $\pi$ as in Proposition~\ref{8a}~(c) can be chosen, that extends to a trivialization of $\pi'$ (necessarily unique,
and $C_1(\BQ)$-equivariant), (d)~the trivialization $\Phi': (\pi^{-1}(z))^{BS} \times Z \isoto Y^{BS \, '}$ 
of $\pi'$ from (c) can be chosen
to respect the stra\-ti\-fications: for all parabolics $Q$ of $Q_1$ satisfying $\adm(Q) = Q_1 \,$, we have
\[
(\Phi')^{-1} \bigl( e(Q) \bigr) = e(Q \cap C_1) \times Z \; ,
\]
where $e(Q \cap C_1)$ denotes the face of $(\pi^{-1}(z))^{BS}$ associated to the parabolic $Q \cap C_1$
(cmp.~\cite[Cor.~4.10]{W1}).

From (a), we conclude that $\pi'$ factorizes canonically as
\[
Y^{BS \, '} \longonto H_C \backslash Y^{BS \, '} \stackrel{\tilde{\pi}'}{\longonto} Z \; .
\]
Here, the intermediate quotient $H_C \backslash Y^{BS \, '}$ contains $H_C \backslash Y^{BS}$ as an open sub-set,
and $\tilde{\pi}'$ extends the map $\tilde{\pi}$ from Definition~\ref{8A}.
According to~(c) and (d), the map $\tilde{\pi}'$ is a trivial stratified fibration. 
Following~(b), and \cite[Thm.~9.3]{BS}, the map $\tilde{\pi}'$ is proper. 

Therefore, the map $\tilde{\pi}'$ provides a $H_1/H_C$-equivariant
compactification of $\tilde{\pi}$.  
This implies \cite{W3}: (e)~the exact triangle 
\[
(E) \quad\quad R \tilde{\pi}_! \longto R \tilde{\pi}_* \longto \partial R \tilde{\pi}_* \stackrel{[1]}{\longto}                    
\]
from Construction~\ref{8Aa} equals $R \tilde{\pi}_*'$ applied to the exact triangle
\[
(T) \quad\quad \j_!  \longto R \j_* \longto \i_* \i^* R \j_* \stackrel{[1]}{\longto} \; , 
\]
where $\j$ and $\i$ denote the complementary immersions
\[
\j: H_C \backslash Y^{BS} \longinto H_C \backslash Y^{BS \, '}
\]
and
\[
\i: H_C \backslash \partial Y^{BS \, '} \longinto H_C \backslash Y^{BS \, '} \; ,
\]
for 
\[
\partial Y^{BS \, '} := \coprod_{\adm(Q) = Q_1, Q \ne Q_1} e(Q) = Y^{BS \, '} - Y \; . 
\] 
Since $K$ is supposed neat, 
the arithmetic sub-group $H_1$ of $Q_1 (\BQ)$ acts freely on $Y^{BS}$ \cite[Sect.~9.5]{BS},
and hence, on $Y^{BS \, '}$. 
The usual construction provides us with canonical extensions 
\[
\eta_g^{BS}: \Rep ( H_1 ) \longto H_1/H_C \, \text{-} \Loc \bigl( H_C \backslash Y^{BS \, '} \bigr) 
\]
and 
\[
\eta_g^{BS}: D^+ ( \Rep H_1 ) \longto D^+ \Bigl( H_1/H_C \, \text{-} \bigl(H_C \backslash Y^{BS \, '}\bigr) \Bigr) 
\]
of the functors $\eta_g$ from Definition~\ref{8A}, related 
to the latter by the formula $\eta_g = \j^* \circ \eta_g^{BS}$.
The map $\j$ being contractible \cite[Lemma~8.3.1, Sect.~9.5]{BS}, Corollary~\ref{1C}~(b)
implies that $\eta_g^{BS} = R \j_* \circ \eta_g$. Thus, the exact triangle $(T \eta_g)$ equals
\[
(T \eta_g) \quad\quad \j_! \eta_g \longto \eta_g^{BS} \longto \i_* \i^* \eta_g^{BS} \stackrel{[1]}{\longto} \; .
\]
In particular: (f)~we have 
\[
\partial R \tilde{\pi}_* \eta_g = R \tilde{\pi}_*' \i_* \i^* \eta_g^{BS} : D^+ ( \Rep H_1 ) 
\longto D^+ ( H_1/H_C \, \text{-} Z ) \; ,
\] 
hence by Definition~\ref{8A}~(b)
\[
\partial R \Gamma (H_C, {\argdot}) 
= R \Gamma \bigl( H_C \backslash \partial Y^{BS \, '} , \eta_g^{BS} ({\argdot})_{\tei H_C \backslash \partial Y^{BS \, '} } \bigr) \; .
\] 
Now let $Q$ be a parabolic sub-group of $G$ satisfying $\adm(Q) = Q_1 \,$. Consider the sub-set
\[
Y_Q:= \coprod_{q \in Q_1(\BQ) / Q(\BQ)} e \bigl( q Q q^{-1} \bigr) 
\]
of $Y^{BS \, '}$.
Each of the faces $e ( q Q q^{-1} )$ being contractible \cite[Sect.~3.9, Rem.~2.4~(1)]{BS}, 
we get: (g)~there is a canonical isomorphism of functors on $D^+ ( \Rep H_1 )$ between
\[
R \Gamma \bigl( H_C , {}_Q \! {\argdot} \bigr) = 
R \Gamma \bigl( H_C , \App( Q_1(\BQ) / Q(\BQ) , {\argdot} ) \bigr) 
\]
and
\[
R \Gamma \bigl( H_C \backslash Y_Q, \eta_g^{BS} ({\argdot})_{\tei H_C \backslash Y_Q} \bigr) 
\] 
(cmp.~Scholie~\ref{7Scho}~(2)). 
If $Q \ne Q_1$, then $Y_Q \subset \partial Y^{BS \, '}$. 
Altogether: (h)~if $Q \ne Q_1$ is a parabolic of $G$ satisfying $\adm(Q) = Q_1 \,$, then
restriction from $\partial Y^{BS \, '}$ to $Y_Q$ induces a natural transformation
\[
\partial R \Gamma (H_C, {\argdot}) \longto R \Gamma \bigl( H_C , {}_Q \! {\argdot} \bigr)
\]
of functors $D^+ (\Rep H_1) \to D^+ ( \Rep (H_1/H_C) )$. 
Using \cite[Prop.~1.25~(b)]{W1} (cmp.~proof of Proposition~\ref{8B}),
one sees that it does not depend on the choice of $Y$.
\end{Cons}

\forget{ ...to be recycled in \cite{W3}!!!...
\begin{Prop} \label{8Q}
Let $Q \ne Q_1$ be a parabolic sub-group of $G$ satisfying $\adm(Q) = Q_1$.
Then the composition 
\[
R \Gamma (H_C, {\argdot}) \longto R \Gamma \bigl( H_C , {}_Q \! {\argdot} \bigr)                 
\]
of the natural transformations
\[
R \Gamma (H_C, {\argdot}) \longto \partial R \Gamma (H_C, {\argdot})                    
\]
(Construction~\ref{8Aa}) and
\[
\partial R \Gamma (H_C, {\argdot}) \longto R \Gamma \bigl( H_C , {}_Q \! {\argdot} \bigr)
\]
(Construction~\ref{8P}) equals the effect of $R \Gamma(H_C, {\argdot})$ on the natural transformation 
\[
\id_{D^+ ( \Rep H_1 )} = {}_{Q_1} \! {\argdot} \longto {}_Q \! {\argdot}
\]
associated to the inclusion $Q \subset Q_1$.
\end{Prop}

\begin{Proof}
Choose $z \in Z$, and put $X:= (\pi^{-1}(z))^{BS}$. According to point~(b) of Construction~\ref{8P},
\[
X \cap Y_Q = \coprod_{q \in Q_1(\BQ) / Q(\BQ)} e \bigl( q (Q \cap C_1) q^{-1} \bigr) \; .
\]  
Also, under a trivialization $\Phi': X \times Z \isoto Y^{BS \, '}$ of $\pi'$
as in point~(d) of Construction~\ref{8P}, we have 
\[
\Phi \bigl( (X \cap Y_Q) \times Z \bigr) = Y_Q \; .
\]
Apply Scholie~\ref{7Scho}~(3) to $X$, $X_1 :=  X \cap Y_Q \,$, $H':=H_1$ and $H:= H_C$. 
\end{Proof}
}
Fix an element
$\uQ = (Q_1 \prec Q_2 \prec \ldots \prec Q_r)$ of $\CC_{(G,\FX)} \,$, such that
$Q_1$ is maximal proper, and $r \ge 2$. 
As in Section~\ref{6}, 
\[
h_{Q_1}: \partial M^K  (G,\FX)^* - \overline{M^K \bigl( Q_1,\FX \bigr)}
\longinto \partial M^K  (G,\FX)^* 
\]
denotes the open immersion of the complement of the closure $\overline{M^K ( Q_1,\FX )}$
of $M^K ( Q_1,\FX )$. \\

In order to connect Construction~\ref{8P} to the results developed previously in this section,
let $\BV \in D^+(\Rep (G(\BQ)))$.
On the one hand, Main Theorem~\ref{8MT}~(a) 
provides us with a canonical isomorphism between the image under $\bar{\mu}$ of 
\[
\partial R \Gamma \bigl( H_C, \Res^{G(\BQ)}_{H_1}(\BV) \bigr)          
\] 
and the restriction to $pr_{gK}(\coprod \FX_1/W_1)$ of 
\[
i_{Q_1}^* R (h_{Q_1})_* h_{Q_1}^* i^*R(i_G)_* \mu_K(\BV) \; .
\]
On the other hand, the component $\uQ$ of the canonical isomorphism from Main Theorem~\ref{7MT}~(a) identifies
the image under $\bar{\mu}$ of
\[
R \Gamma \bigl( H_C, {}_{\bigcap \uQ} \! \Res^{G(\BQ)}_{H_1}(\BV) \bigr)
\]
and the restriction to $pr_{gK}(\coprod \FX_1/W_1)$ of 
\[
i_{Q_1}^* R(i_{Q_2})_* i_{Q_2}^* \ldots R(i_{Q_r})_* i_{Q_r}^* i^*R(i_G)_* \mu_K(\BV) \; .
\] 

\begin{Prop} \label{8R}
The restriction to $pr_{gK}(\coprod \FX_1/W_1)$ of the adjunction
\[
R (h_{Q_1})_* h_{Q_1}^* \longto R(i_{Q_2})_* i_{Q_2}^* \ldots R(i_{Q_r})_* i_{Q_r}^* \; ,
\]
evaluated at $i^*R(i_G)_* \mu_K(\BV)$,
equals, under the above identifications, 
the image under $\bar{\mu}$ of
the value on $\Res^{G(\BQ)}_{H_1}(\BV)$ of the natural transformation
\[
\partial R \Gamma (H_C, {\argdot}) \longto R \Gamma \bigl( H_C , {}_{\bigcap \uQ} \! {\argdot} \bigr)
\]
from Construction~\ref{8P}.
\end{Prop}

\begin{Proofof}{Main Theorem~\ref{8MT} and Proposition~\ref{8R}}
\forget{Since $K$ is assumed neat, the map
\[
pr_{gK}: \coprod \FX_1/W_1 \longto pr_{gK} \bigl( \coprod \FX_1/W_1 \bigr)
\]
is a covering of the open sub-set $pr_{gK} ( \coprod \FX_1/W_1 )$ of $M^K ( Q_1,\FX )(\BC)$, 
with covering group $H_1/H_C$ \cite[Prop.~7.22~(a)]{W1}.
}
Set $\CV := \mu_K(\BV)$. The proof of Main Theorem~\ref{8MT} will use
the formulae from Corollary~\ref{6D}, which are expressed in terms of the object
\[
\CV^{\partial BS}_{\tei \! e^K ( Q_1,P(\BA_f) )'} 
\]
of $D^+ ( e^K ( Q_1,P(\BA_f) )' )$. In order to control the restrictions of this object to 
$(p^K)^{-1} ( pr_{gK} ( \coprod \FX_1/W_1 ) )$, recall the commutative diagram
\[
(A) \quad \quad \vcenter{\xymatrix@R-10pt{
p^{-1} \bigl( \coprod \FX_1/W_1 \bigr) \ar@{^{>>}}[d]_-{\alpha} 
\ar@{^{ (}->}[r]^-{x \mapsto (x,gK)} & e(Q_1)' \times G(\BA_f)/K \ar@{^{>>}}[d] \\
(p^K)^{-1} \bigl( pr_{gK} \bigl( \coprod \FX_1/W_1 \bigr) \bigr)
\ar@{^{ (}->}[r] & e^K ( Q_1,G(\BA_f) )'
\\}} 
\]
of immersions and coverings, with covering groups $H_1$ and $Q_1(\BQ)$, respectively
(see Auxiliary Construction~\ref{7Aux}~(a), (c)). The inclusion of 
$H_1$ into $Q_1(\BQ)$ is the one of the stabilizer of $p^{-1} ( \coprod \FX_1/W_1)$.
In particular, the usual construction yields functors
\[
\bar{\mu}': \Rep H_1 \longto \Loc \Bigl( (p^K)^{-1} \bigl( pr_{gK} \bigl( \coprod \FX_1/W_1 \bigr) \bigr) \Bigr)
\] 
and
\[
\bar{\mu}': D^+ ( \Rep H_1 )  \longto 
D^+ \Bigl( (p^K)^{-1} \bigl( pr_{gK} \bigl( \coprod \FX_1/W_1 \bigr) \bigr) \Bigr) \; .
\]
The $\uQ = (Q_1)$-component of Proposition~\ref{7G} tells us that for the restriction of 
$\CV^{\partial BS}_{\tei \! e^K ( Q_1,P(\BA_f) )'}$ to $(p^K)^{-1} ( pr_{gK} ( \coprod \FX_1/W_1 ) )$, 
we have
\[
\CV^{\partial BS}_{\tei \! (p^K)^{-1} ( pr_{gK} ( \coprod \FX_1/W_1 ) )} 
= \bar{\mu}' \bigl( \Res^{G(\BQ)}_{H_1}(\BV) \bigr) 
 \in D^+ \Bigl( \bigl( (p^K)^{-1} \bigl( pr_{gK} \bigl( \coprod \FX_1/W_1 \bigr) \bigr) \Bigr) \; .
\]
\forget{Therefore, the inverse image of
\[
\CV^{\partial BS}_{\tei \! (p^K)^{-1} ( [(z_0,gK)] )} \in D^+ \Bigl( \bigl( (p^K)^{-1}( [(z_0,gK)] ) \bigr) \Bigr) 
\]
on the covering $p^{-1} (z_0)$ equals the constant object
\[
\Res^{G(\BQ)}_{H_C}(\BV) \in D^+ ( \Rep H_C ) \; .
\]
Denote (\emph{cf.} Main Theorem~\ref{5MT} and Corollary~\ref{6D})
by $\j$ the open immersion of $e^K ( Q_1,G(\BA_f))$ into 
$e^K ( Q_1,G(\BA_f) )'$, and by 
\[
\i: \partial e^K \bigl( Q_1,G(\BA_f) \bigr)' \longinto e^K \bigl( Q_1,G(\BA_f) \bigr)'
\]
the closed immersion complementary to $\j$.} 
\noindent (a): according to Corollary~\ref{6D}~(a), the isomorphism from Theorem~\ref{6MT2}~(a) identifies 
the exact triangle
\[
R i_{Q_1}^! i^*R(i_G)_* \mu_K(\BV) \to  
i_{Q_1}^* R(i_G)_* \mu_K(\BV) \to
i_{Q_1}^* R (h_{Q_1})_* h_{Q_1}^* i^*R(i_G)_* \mu_K(\BV) \stackrel{[1]}{\to} 
\] 
with 
\[
Rp^K_* \bigl( \j_! \CV^{\partial BS}_{\tei \! e^K ( Q_1,G(\BA_f) )} \bigr) \to
Rp^K_* \bigl( \CV^{\partial BS}_{\tei \! e^K ( Q_1,P(\BA_f) )'} \bigr) \to
Rp^K_* \bigl( \i_* \CV^{\partial BS}_{\tei \! \partial e^K ( Q_1,G(\BA_f) )'} \bigr) \stackrel{[1]}{\to} \; .
\]  
Restriction to $pr_{gK} ( \coprod \FX_1/W_1 )  \subset M^K ( Q_1,\FX )(\BC)$ of that isomorphism induces
an identification between the restriction to $pr_{gK} ( \coprod \FX_1/W_1 )$ of
\[
R i_{Q_1}^! i^*R(i_G)_* \mu_K(\BV) \to  
i_{Q_1}^* R(i_G)_* \mu_K(\BV) \to
i_{Q_1}^* R (h_{Q_1})_* h_{Q_1}^* i^*R(i_G)_* \mu_K(\BV) \stackrel{[1]}{\to} 
\] 
and the image under $Rp^K_*$ of 
$T \CV^{\partial BS}_{\tei \! (p^K)^{-1} ( pr_{gK} ( \coprod \FX_1/W_1 ) )} 
= T \bar{\mu}' ( \Res^{G(\BQ)}_{H_1}(\BV) ) $, where $T$ is the localization triangle 
\[
\j_! \j^* \longto \id_{(p^K)^{-1} ( pr_{gK} ( \coprod \FX_1/W_1 ) )} \longto \i_* \i^* \stackrel{[1]}{\longto} \; .
\]
Here, we keep the notations $\j$ for the open immersion of
\[
(p^K)^{-1} \bigl( pr_{gK} \bigl( \coprod \FX_1/W_1 \bigr) \bigr) \cap e^K \bigl( Q_1,G(\BA_f) \bigr) \; ,
\] 
and $\i$ for the complementary closed immersion of
\[
(p^K)^{-1} \bigl( pr_{gK} \bigl( \coprod \FX_1/W_1 \bigr) \bigr) \cap \partial e^K \bigl( Q_1,G(\BA_f) \bigr)' 
\] 
into $(p^K)^{-1} ( pr_{gK} ( \coprod \FX_1/W_1 ) ) \subset e^K ( Q_1,G(\BA_f) )'$.
By \cite[Prop.~7.22~(c)]{W1}, base change \emph{via} $pr_{gK}$ 
of the restriction to $(p^K)^{-1} ( pr_{gK} ( \coprod \FX_1/W_1 ) )$ of $p^K$ yields the (Cartesian) diagram
\[
(B) \quad \quad \vcenter{\xymatrix@R-10pt{
H_C \backslash p^{-1} \bigl( \coprod \FX_1/W_1 \bigr) \ar[r]^-{\bar{p}} \ar@{->>}[d]_-{\beta} & 
\coprod \FX_1/W_1 \ar@{->>}[d]^-{pr_{gK}}  \\
(p^K)^{-1} \bigl( pr_{gK} \bigl( \coprod \FX_1/W_1 \bigr) \bigr) \ar[r]^-{p^K} & 
pr_{gK} \bigl( \coprod \FX_1/W_1 \bigr) 
\\}} \quad .
\]
The idea is to use this latter diagram 
in order to control the exact triangle $Rp^K_* (T \bar{\mu}' ( \Res^{G(\BQ)}_{H_1}(\BV) ))$.
Note that composition of $p^{-1} ( \coprod \FX_1/W_1 ) \onto H_C \backslash p^{-1} ( \coprod \FX_1/W_1 )$
with the map $\beta$ from diagram~(B) yields the covering $\alpha$ from diagam~(A).
Therefore, $\beta$ is a covering, with covering group $H_1/H_C \,$; actually, this also follows from 
diagram~(B) being Cartesian, and from the corresponding property of $pr_{gK}$ (recall that this latter property
enabled us to define the functor $\bar{\mu}$).
It follows that the inverse image functors $\beta^*$ and $pr_{gK}^*$ induce equivalences of categories
\[
\Loc \Bigl( (p^K)^{-1} \bigl( pr_{gK} \bigl( \coprod \FX_1/W_1 \bigr) \bigr) \Bigr) \isoto
H_1/H_C \, \text{-} Loc \Bigl( H_C \backslash p^{-1} \bigl( \coprod \FX_1/W_1 \bigr) \Bigr)
\]
and
\[
\Loc \Bigl( pr_{gK} \bigl( \coprod \FX_1/W_1 \bigr) \Bigr) \isoto
H_1/H_C \, \text{-} Loc \bigl( \coprod \FX_1/W_1 \bigr) \; ,
\]
respectively, as well as
\[
D^+ \Bigl( (p^K)^{-1} \bigl( pr_{gK} \bigl( \coprod \FX_1/W_1 \bigr) \bigr) \Bigr) \isoto
D^+ \Bigl( H_1/H_C \, \text{-} \bigl( H_C \backslash p^{-1} \bigl( \coprod \FX_1/W_1 \bigr) \bigr) \Bigr)
\]
and
\[
D^+ \Bigl( pr_{gK} \bigl( \coprod \FX_1/W_1 \bigr) \Bigr) \isoto
D^+ \Bigl( H_1/H_C \, \text{-} \bigl( \coprod \FX_1/W_1 \bigr) \Bigr) \; ,
\]
respectively. On the level of derived categories, the inverses of the equivalences are given by 
$R \Gamma (H_1/H_C, {\argdot}) \circ R \beta_*$ and $R \Gamma (H_1/H_C, {\argdot}) \circ R (pr_{gK})_* \,$,
respectively (cmp.~\cite{W3}). In particular, we may, \emph{via} $\beta^*$ and 
$R \Gamma (H_1/H_C, {\argdot}) \circ R \beta_*$, re-interpret the functor $\bar{\mu}'$ as giving
\[
\bar{\mu}': D^+ ( \Rep H_1 )  \longto 
D^+ \Bigl( H_1/H_C \, \text{-} \bigl( H_C \backslash p^{-1} \bigl( \coprod \FX_1/W_1 \bigr) \bigr) \Bigr) \; .
\]
Also, the localization triangle $T$ in $D^+ ( (p^K)^{-1} ( pr_{gK} ( \coprod \FX_1/W_1 ) ) )$ corresponds to the triangle,
still denoted $T$, 
\[
\j_! \j^* \longto \id_{H_C \backslash p^{-1} ( \coprod \FX_1/W_1 )} \longto \i_* \i^* \stackrel{[1]}{\longto} 
\]
in $D^+ ( H_1/H_C \, \text{-} ( H_C \backslash p^{-1} ( \coprod \FX_1/W_1 ) ) )$, where 
we keep the notations $\j$ for the open immersion of
\[
\beta^{-1} \Bigl( (p^K)^{-1} \bigl( pr_{gK} \bigl( \coprod \FX_1/W_1 \bigr) \bigr) \cap e^K \bigl( Q_1,G(\BA_f) \bigr) \Bigr) \; ,
\] 
and $\i$ for the complementary closed immersion of
\[
\beta^{-1} \Bigl( (p^K)^{-1} \bigl( pr_{gK} \bigl( \coprod \FX_1/W_1 \bigr) \bigr) \cap \partial e^K \bigl( Q_1,G(\BA_f) \bigr)' 
                   \Bigr) 
\] 
into $H_C \backslash p^{-1} ( \coprod \FX_1/W_1 )$. 

To summarize the discussion so far, we proved that
the restriction to $pr_{gK} ( \coprod \FX_1/W_1 )$ of
\[
R i_{Q_1}^! i^*R(i_G)_* \mu_K(\BV) \to  
i_{Q_1}^* R(i_G)_* \mu_K(\BV) \to
i_{Q_1}^* R (h_{Q_1})_* h_{Q_1}^* i^*R(i_G)_* \mu_K(\BV) \stackrel{[1]}{\to} 
\] 
--- in other words, our object of interest --- equals 
\[
R \Gamma (H_1/H_C, {\argdot}) \circ R (pr_{gK})_* \circ R \bar{p}_* T \bar{\mu}' \bigl( \Res^{G(\BQ)}_{H_1}(\BV) \bigr) \; .
\]
Given that $\bar{\mu}$ equals the restriction of $R \Gamma (H_1/H_C, {\argdot}) \circ R (pr_{gK})_*$ 
to the sub-category of constant objects, our claim will be established once we show that the exact triangle
$R \bar{p}_* T \bar{\mu}' ( \Res^{G(\BQ)}_{H_1}(\BV) )$ in
$D^+ ( H_1/H_C \, \text{-} ( \coprod \FX_1/W_1 ) )$ is constant, with value equal to
\[
R \Gamma_c \bigl( H_C, \Res^{G(\BQ)}_{H_1}(\BV) \bigr) \to R \Gamma \bigl( H_C, \Res^{G(\BQ)}_{H_1}(\BV) \bigr) 
\to \partial R \Gamma \bigl( H_C, \Res^{G(\BQ)}_{H_1}(\BV) \bigr) \stackrel{[1]}{\to} \; .          
\] 

Recall (Proposition~\ref{7I}) the equality 
\[
p^{-1} \bigl( \coprod \FX_1 / W_1 \bigr) = e(Q_1)' \subset \FX^{BS} \; .
\]
According to \cite[Prop.~3.7]{W1}, the pre-image of 
\[
(p^K)^{-1} \bigl( pr_{gK} \bigl( \coprod \FX_1/W_1 \bigr) \bigr) \cap e^K \bigl( Q_1,G(\BA_f) \bigr) 
\] 
under the covering 
\[
\alpha: e(Q_1)' = p^{-1} \bigl( \coprod \FX_1 / W_1 \bigr) \longonto (p^K)^{-1} \bigl( pr_{gK} \bigl( \coprod \FX_1/W_1 \bigr) \bigr)
\]
equals $e(Q_1)$. It follows that the open immersion in the triangle $T$ is
\[
\j : H_C \backslash e(Q_1) \longinto H_C \backslash e(Q_1)' = H_C \backslash p^{-1} \bigl( \coprod \FX_1 / W_1 \bigr) \; .
\]
\forget{ and 
\[
\i: H_C \backslash \partial e(Q_1)' \longinto H_C \backslash e(Q_1)' \; .
\]
Therefore, the pre-image of the complement
\[
(p^K)^{-1} \bigl( pr_{gK} \bigl( \coprod \FX_1/W_1 \bigr) \bigr) \cap \partial e^K \bigl( Q_1,G(\BA_f) \bigr)' 
\] 
under $\alpha$ equals the complement of $e(Q_1)$ in $e(Q_1)'$, which we shall denote by $\partial e(Q_1)'$.}
Consider $p_{\tei e(Q_1)}: e(Q_1) \to \coprod \FX_1 / W_1$.
Recall that by construction \cite[Def.~3.1~(b)]{W1}, the face $e(Q_1)$ has as many connected components as $\FX \,$.
In particular, the number of connected components is finite. 
According to \cite[Thm.~6.9~(d), Prop.~6.5]{W1}, each fibre of $p_{\tei e(Q_1)}$ underlies a space of type $S - \BQ$ under $C_1$.
Therefore, the fibres of $p_{\tei e(Q_1)}$ are non-empty and connected. It follows that 
$p_{\tei e(Q_1)}$ induces a bijection on the level
of sets of connected components.   
The group $K$ is assumed neat; therefore, so is the 
group $H_1 \,$, which implies that its action on the finite set $\pi_0(\FX) = \pi_0(e(Q_1)) = \pi_0(\coprod \FX_1 / W_1)$ 
is trivial. 

Let $\FX^0$ be a connected component of $\FX \,$, and let $e(Q_1)^0$ and $(\FX_1 / W_1)^0$ be the corresponding connected
components of $e(Q_1)$ and of $\coprod \FX_1 / W_1 \,$, respectively.
According to \cite[Thm.~1.21]{W1}, $\FX^0$ underlies a space of type $S - \BQ$ under $G$. 
Therefore \cite[Sect.~3.9]{BS},
$Y:= e(Q_1)^0$ underlies a space of type $S - \BQ$ under $Q_1 \,$. 
Note that in order to obtain an action of the whole of $G(\BR)$ on 
$\FX^0$, the action of $\Stab_{G(\BR)}(\FX^0)$ coming from our Shimura data needs to be
extended in a precise way \cite[Sect.~1]{W1}. Since $H_1$ is contained in $\Stab_{G(\BR)}(\FX^0)$, the definition
of this extension does not matter to us. 

The fibres of $p_{\tei e(Q_1)}$ underlying spaces of type $S - \BQ$ under $C_1 \,$, 
they are homogeneous under $C_1(\BR)$. This means that
the morphism $\pi:= p_{\tei Y}: Y \onto (\FX_1 / W_1)^0$ factors through the projection to $Z := C_1(\BR) \backslash Y$.
More precisely, there is a (unique) bijective continuous map $\gamma: Z \to (\FX_1 / W_1)^0$, such that the composition
of the projection $Y \onto Z$ with $\gamma$ equals $\pi$. But $\pi$ admits 
(continuous) right inverses \cite[Thm.~7.17~(a), (c), (e)]{W1};
therefore, so does $\gamma$. \emph{Via} the homeomorphism $\gamma$, we identify $Z$ with $(\FX_1 / W_1)^0$, and 
$\pi$ with the projection $Y \onto Z$. Altogether, we find ourselves in the situation of Definition~\ref{8A}.

The quotient $H_C \backslash Y$ is open in 
$H_C \backslash p^{-1} ( \FX_1 / W_1 )^0$. The latter quotient is a connected component of
$H_C \backslash p^{-1} ( \coprod \FX_1 / W_1 )$.
Diagram~(B) being Cartesian, and $p^K$ being proper, the morphism
\[
\bar{p}: H_C \backslash p^{-1} \bigl( \coprod \FX_1/W_1 \bigr) \longonto \coprod \FX_1/W_1  
\]
is proper. Hence so is its base change
\[
\bar{p}_{\tei H_C \backslash p^{-1} (Z)}: H_C \backslash p^{-1} (Z) \longonto Z
\]
\emph{via} the inclusion of $Z$ into $\coprod \FX_1 / W_1$. 
It follows that $\bar{p}_{\tei H_C \backslash p^{-1} (Z)}$ provides a $H_1/H_C$-equivariant
compactification of the map $\tilde{\pi}: H_C \backslash Y \onto Z$.  
But this implies \cite{W3} that the exact triangle 
\[
R \tilde{\pi}_! \longto R \tilde{\pi}_* \longto \partial R \tilde{\pi}_* \stackrel{[1]}{\longto}                    
\]
from Construction~\ref{8Aa} equals $R (\bar{p}_{\tei H_C \backslash p^{-1} (Z)})_*$ applied to the
restriction of the exact triangle
\[
\bigl( T (R \j_*) \bigr) \quad\quad \j_!  \longto R \j_* \longto \i_* \i^* R \j_* \stackrel{[1]}{\longto} 
\]
to $H_C \backslash Y \subset H_C \backslash p^{-1} (Z)$. Now the functor
\[
\eta_g: D^+ ( \Rep H_1 ) \longto D^+ \bigl( H_1/H_C \, \text{-} (H_C \backslash Y) \bigr) 
\]
from Definition~\ref{8A} equals $\j^* \circ \bar{\mu}'$, followed by the restriction to $H_C \backslash Y$.
Our claim follows from contractibility of $\j$ \cite[Prop.~3.6]{W1}, and from Corollary~\ref{1C}~(b).

\noindent (b): according to Corollary~\ref{6D}~(b), we have
\[
R i_{Q_1}^! i^*R(i_G)_* \CV
\cong Rp^K_* \bigl( \j_! \CV^{\partial BS}_{\tei \! e^K ( Q_1,G(\BA_f) )} \bigr)[-(r-1)] \; .
\]  
Repeat the computation of $Rp^K_* ( \j_! \CV^{\partial BS}_{\tei \! e^K ( Q_1,G(\BA_f) )} )$ from (a).

\noindent (c): this claim is identical to Corollary~\ref{6D}~(c).

\noindent (d): use (b) and (c).

\noindent (e): this claim is a formal consequence of (a) (for all $g$).

\medskip

As for the proof of Proposition~\ref{8R}, we place ourselves back in the situation of (a) above
($Q_1$ maximal proper), and assume in addition that $Q_1$ is the first component of an element
$\uQ = (Q_1 \prec Q_2 \prec \ldots \prec Q_r)$ of $\CC_{(G,\FX)} \,$, with $r \ge 2$. 

Recall that the $H_1/H_C$-equivariant morphism
\[
\bar{p}_{\tei H_C \backslash p^{-1} (Z)}: H_C \backslash p^{-1} (Z) \longonto Z \; ,
\]
which is an extension of $\tilde{\pi}: H_C \backslash Y \onto Z$. 
According to Proposition~\ref{7I}, we have $p^{-1} (Z) = Y^{BS \, '}$ in the notation of Construction~\ref{8P}.

Now apply Proposition~\ref{6F} and Remark~\ref{6G} (with $\CV = \mu_K(\BV)$).  
\end{Proofof}

Remark~\ref{8C}~(c) allows to replace in Main Theorem~\ref{8MT} the expressions
$R \Gamma_c ( H_C, \Res^{G(\BQ)}_{H_1}(\BV))$, $R \Gamma ( H_C, \Res^{G(\BQ)}_{H_1}(\BV))$,
$\partial R \Gamma ( H_C, \Res^{G(\BQ)}_{H_1}(\BV))$, 
and the exact triangle
\[
R \Gamma_c \bigl( H_C, \Res^{G(\BQ)}_{H_1}(\BV) \bigr) \to R \Gamma \bigl( H_C, \Res^{G(\BQ)}_{H_1}(\BV) \bigr) 
\to \partial R \Gamma \bigl( H_C, \Res^{G(\BQ)}_{H_1}(\BV) \bigr) \stackrel{[1]}{\to}           
\] 
by $R \Gamma_c ( H_C/H_W, \BX )$, $R \Gamma ( H_C/H_W, \BX )$, $\partial R \Gamma ( H_C/H_W, \BX )$, and
\[
R \Gamma_c ( H_C/H_W, \BX )  \longto R \Gamma ( H_C/H_W, \BX ) 
\longto \partial R \Gamma ( H_C/H_W, \BX )  
\stackrel{[1]}{\longto} \; ,                   
\]
respectively, where $\BX:= R \Gamma ( H_W, \Res^{G(\BQ)}_{H_1}(\BV))$.
As in the preceding section (see Variant~\ref{7Var1}), 
this observation allows to reformulate Main Theorem~\ref{8MT} for those objects
of $D^+ ( \Rep (G(\BQ) ) )$ coming about as restrictions of algebraic representations of $G$. 
The explicit reformulation will be restricted to
part~(a) of Main Theorem~\ref{8MT}; the remaining parts will be left to the reader. 

\begin{Var} \label{8Var1}
Let $F$ be a field of characteristic zero. Assume that $Q_1$ is maximal proper. Let $\BV \in D^+ ( \Rep_F G )$, and
define
\[
\BX:= \Res^{\bar{Q}_1}_{H_1/H_W} R \Gamma \bigl( W_1, \bigl( \Res^G_{Q_1} \BV \bigr) \bigr)
\]
($\Res^G_{Q_1}$ and $\Res^{\bar{Q}_1}_{H_1/H_W} =$ 
the restrictions from $G$ to $Q_1$ and from $\bar{Q}_1 = Q_1/W_1$ to $H_1/H_W$, respectively). 
Then the isomorphism from Main Theorem~\ref{8MT}~(a) restricts to give a canonical isomorphism of exact triangles 
in the derived category $D^+ ( pr_{gK} (\coprod \FX_1/W_1 ))$
between the image under $\bar{\mu}$ of the exact triangle 
\[
R \Gamma_c ( H_C/H_W, \BX )  \longto R \Gamma ( H_C/H_W, \BX ) 
\longto \partial R \Gamma ( H_C/H_W, \BX )  
\stackrel{[1]}{\longto}                  
\]
and the restriction to $pr_{gK}(\coprod \FX_1/W_1)$ of the exact triangle
\[
R i_{Q_1}^! i^*R(i_G)_* \mu_K(\BV) \to  
i_{Q_1}^* R(i_G)_* \mu_K(\BV) \to
i_{Q_1}^* R (h_{Q_1})_* h_{Q_1}^* i^*R(i_G)_* \mu_K(\BV) \stackrel{[1]}{\to}  
\]  
in $D^+ (M^K ( Q_1,\FX )(\BC))$.
\end{Var}

\begin{Proof}
Use the remark preceding our statement, together with the fact that as functors
\[
D^+ \bigl( \Rep_F Q_1 \bigr) \longto D^+ \bigl( \Rep H_1/H_W \bigr) \; ,
\]
the compositions
$\Res^{\bar{Q}_1}_{H_1/H_W} \circ R \Gamma ( W_1, {\argdot} )$ and
$R \Gamma ( H_W, {\argdot} ) \circ \Res^{Q_1}_{H_1}$
($\Res^{Q_1}_{H_1} :=$ the restriction from $Q_1$ to $H_1$)
are canonically isomorphic, as $H_W$ is an arithmetic sub-group of the unipotent
group $W_1$.
\end{Proof}

Recall (Definition~\ref{7K}~(a)) that $M^K (Q_1;P_1(\BA_f)gK,\FX) \subset M^K (Q_1,\FX)$ 
is defined as the image 
of $\coprod \FX_1/W_1 \times P_1(\BA_f)gK/K$ under the projection 
\[
\FX^* \times G (\BA_f) / K 
\longonto G (\BQ) \backslash \bigl( \FX^* \times G (\BA_f) / K \bigr) = M^K (G,\FX)^* (\BC) \; .
\]
The stabilizer in $G(\BQ)$ of $\coprod \FX_1/W_1 \times P_1(\BA_f)gK/K$ equals the group
\[ 
H_1' = Q_1(\BQ) \cap \bigl( P_1(\BA_f) \cdot gKg^{-1} \bigr) 
\]
from Definition~\ref{7K}~(e), and the kernel of the action of $H_1'$
on the space $\coprod \FX_1/W_1 \times P_1(\BA_f)gK/K$
is the group
\[
H_C' = C_1(\BQ) \cap \bigl( W_1(\BA_f) \cdot gKg^{-1} \bigr) 
\]
from Definition~\ref{7M}~(a). The induced action of the quotient $H_1'/H_C'$ is free (Proposition~\ref{7N}~(c)),
giving rise to the functor
\[
\mu_{\pi_1(K_1)}: D^+ \bigl( \Rep(H_1'/H_C') \bigr) \longto 
D^+ \bigl( M^K (Q_1;P_1(\BA_f)gK,\FX)(\BC) \bigr) 
\]
from Variant~\ref{7Var2}.

\begin{Var} \label{8Var2}
Let $F$ be a field of characteristic zero. Assume that $Q_1$ is maximal proper. Let $\BV \in D^+ ( \Rep_F G )$, and
define
\[
\BX':= \Res^{\bar{Q}_1}_{H_1'/(W_1(\BQ))} R \Gamma \bigl( W_1, \bigl( \Res^G_{Q_1} \BV \bigr) \bigr)
\]
($\Res^{\bar{Q}_1}_{H_1'/(W_1(\BQ))} =$ 
the restriction from $\bar{Q}_1$ to $H_1'/(W_1(\BQ))$). 
There is a cano\-nical isomorphism
in the derived category $D^+ ( M^K (Q_1;P_1(\BA_f)gK,\FX)(\BC) )$
bet\-ween the image under $\mu_{\pi_1(K_1)}$ of the exact triangle 
\[
R \Gamma_c ( H_C/H_W, \BX' ) \longto R \Gamma ( H_C/H_W, \BX' ) 
\longto \partial R \Gamma ( H_C/H_W, \BX' )  
\stackrel{[1]}{\longto}                  
\]
and the restriction to $M^K (Q_1;P_1(\BA_f)gK,\FX)(\BC)$ of the exact triangle
\[
R i_{Q_1}^! i^*R(i_G)_* \mu_K(\BV) \to  
i_{Q_1}^* R(i_G)_* \mu_K(\BV) \to
i_{Q_1}^* R (h_{Q_1})_* h_{Q_1}^* i^*R(i_G)_* \mu_K(\BV) \stackrel{[1]}{\to}  
\]  
in $D^+ (M^K ( Q_1,\FX )(\BC))$.
\end{Var}

In order to get the $H_1'/H_C'$-equivariance of the exact triangle
\[
R \Gamma_c ( H_C/H_W, \BX' ) \longto R \Gamma ( H_C/H_W, \BX' ) 
\longto \partial R \Gamma ( H_C/H_W, \BX' )  
\stackrel{[1]}{\longto}               
\]
necessary for $\mu_{\pi_1(K_1)}$ to be applicable, recall first that according to 
Proposition~\ref{7N}~(a), the inclusion of $H_C$ into $H_C'$ induces an identification of 
$H_C/H_W$ with $H_C'/(W_1(\BQ))$. Then apply Remark~\ref{8C}~(c) to the quotients
$H_C/H_W = H_C'/(W_1(\BQ))$ and $H_1'/(W_1(\BQ))$. \\

Recall (Definition~\ref{7K}~(b)) that the immersion of $M^K (Q_1;P_1(\BA_f)gK,\FX)$ into $\partial M^K  (G,\FX)^*$
is denoted by $i_{Q_1;P_1(\BA_f)gK}$.
The terms in the restriction to $M^K (Q_1;P_1(\BA_f)gK,\FX)(\BC)$ of the exact triangle
\[
R i_{Q_1}^! i^*R(i_G)_* \mu_K(\BV) \to  
i_{Q_1}^* R(i_G)_* \mu_K(\BV) \to
i_{Q_1}^* R (h_{Q_1})_* h_{Q_1}^* i^*R(i_G)_* \mu_K(\BV) \stackrel{[1]}{\to}  
\]  
are therefore equal to $Ri_{Q_1;P_1(\BA_f)gK}^! i^*R(i_G)_* \mu_K(\BV)$,
$i_{Q_1;P_1(\BA_f)gK}^* i^*R(i_G)_* \mu_K(\BV)$, and
$i_{Q_1;P_1(\BA_f)gK}^* R (h_{Q_1})_* h_{Q_1}^* i^*R(i_G)_* \mu_K(\BV)$, respectively.

\medskip

\begin{Proofof}{Variant~\ref{8Var2}}
We repeat the relevant information from the proof of Variant~\ref{7Var2}.

The space $M^K (Q_1;P_1(\BA_f)gK,\FX)(\BC)$ is covered by 
open and closed sub-sets of the form 
$pr_{p_1gK}(\coprod \FX_1/W_1)$, for $p_1 \in P_1(\BA_f)$. It is therefore sufficient to
prove the claim after applying the restriction from $M^K (Q_1;P_1(\BA_f)gK,\FX)(\BC)$ to 
$pr_{p_1gK}(\coprod \FX_1/W_1)$, for any $p_1 \in P_1(\BA_f)$. Defining
\[
H_1 (p_1gK) := Q_1(\BQ) \cap (p_1g)K(p_1g)^{-1} \; ,
\]
\[
H_W (p_1gK) := W_1(\BQ) \cap H_1 (p_1gK) \; ,
\]
and noting that $H_1 (p_1gK)$ is the stabilizer of
$\coprod \FX_1/W_1 \times \{ p_1gK \}$ in $H_1' \,$, 
we obtain the statement from Variant~\ref{8Var2}, 
where $g$ is replaced by $p_1g$, except for the use of
$H_C/H_W = H_C'/W_1(\BQ)$ instead of $H_C (p_1gK)/H_W (p_1gK)$, where
\[
H_C (p_1gK) := C_1(\BQ) \cap (p_1g)K(p_1g)^{-1} \; .
\]
Now (1)~$H_C (p_1gK)$ is a sub-group of $H_C'$,
(2)~$H_C'$ does not change when $g$ is replaced by $p_1g$, 
(3)~the inclusion $H_C (p_1gK) \into H_C'$ induces an isomorphism
\[
H_C (p_1gK)/H_W (p_1gK) \isoto H_C'/W_1(\BQ) 
\] 
(Proposition~\ref{7N}~(a), applied to $p_1g$ instead of $g$).
\end{Proofof}

\begin{Rem} \label{8S}
The isomorphism from Variant~\ref{8Var2} restricts to give a 
canonical isomorphism in the derived category of sheaves on 
\[
M^K (Q_1;\FX_1/W_1;P_1(\BA_f)gK,\FX) \subset 
M^K (Q_1;P_1(\BA_f)gK,\FX) 
\]
(see Remark~\ref{7L}~(a)),
for each individual boundary component $(P_1,\FX_1)$ associated to $Q_1$.
As in Remark~\ref{7O},
the correct equivariance statement requires the group $H_1'$ to be replaced by 
\[
\Stab_{Q_1(\BQ)}(\FX_1/W_1) \cap \bigl( P_1(\BA_f) \cdot gKg^{-1} \bigr) \subset H_1' \; .
\]
\end{Rem}

\begin{Cor} \label{8T}
Let $F$ be a field of characteristic zero. Assume that $Q_1$ is maximal proper.
Let $\BV \in D^+ ( \Rep_F G )$. \\[0.1cm]
(a)~There are canonical and
functorial $E_2$-spectral sequences
\[
E^{p,s}_{c,2} \Longrightarrow \CH^{p+s} \bigl( R i^!_{Q_1;P_1(\BA_f)gK} i^* R(i_G)_* \bigr) \circ \mu_K (\BV) \; ,
\]
\[
E^{p,s}_2 \Longrightarrow \CH^{p+s} \bigl( i^*_{Q_1;P_1(\BA_f)gK} i^* R(i_G)_* \bigr) \circ \mu_K (\BV) \; ,
\]
and
\[
E^{p,s}_{\partial,2} 
    \Longrightarrow \CH^{p+s} \bigl( i^*_{Q_1;P_1(\BA_f)gK} R (h_{Q_1})_* h_{Q_1}^* i^* R(i_G)_* \bigr) \circ \mu_K (\BV) \; ,
\]
where
\[
E^{p,s}_{c,2} := \mu_{\pi_1(K_1)} \circ H^p_c \bigl( H_C/H_W, H^s ( W_1, \Res^G_{Q_1}\BV) \bigr) 
\]
($H^{\bullet}_c :=$ the cohomological functor associated to $R \Gamma_c$),
\[
E^{p,s}_2 := \mu_{\pi_1(K_1)} \circ H^p \bigl( H_C/H_W, H^s ( W_1, \Res^G_{Q_1}\BV) \bigr) \; ,
\]
and
\[
E^{p,s}_{\partial,2} := \mu_{\pi_1(K_1)} \circ \partial H^p \bigl( H_C/H_W, H^s ( W_1, \Res^G_{Q_1}\BV) \bigr) 
\]
($\partial H^{\bullet}:=$ the cohomological functor associated to $\partial R \Gamma$). \\[0.1cm]
(b)~The spectral sequence
\[
E^{p,s}_2 \Longrightarrow \CH^{p+s} \bigl( i^*_{Q_1;P_1(\BA_f)gK} i^* R(i_G)_* \bigr) \circ \mu_K (\BV) 
\]
is identical to the spectral sequence $(E^\argast_{\uQ})$ from Corollary~\ref{7P}~(a), for 
$\uQ := (Q_1)$. \\[0.1cm]
(c)~The three types of morphisms 
in the long exact cohomology sequence associated to
the restriction to $M^K (Q_1;P_1(\BA_f)gK,\FX)(\BC)$ of the exact triangle
\[
R i_{Q_1}^! i^*R(i_G)_* \mu_K(\BV) \to  
i_{Q_1}^* R(i_G)_* \mu_K(\BV) \to
i_{Q_1}^* R (h_{Q_1})_* h_{Q_1}^* i^*R(i_G)_* \mu_K(\BV) \stackrel{[1]}{\to}  
\]  
are canonically and functorially extended to morphisms of spectral sequences.
On the $E_2$-terms
\[
E^{p,s}_{c,2} = \mu_{\pi_1(K_1)} \circ H^p_c \bigl( H_C/H_W, H^s ( W_1, \Res^G_{Q_1}\BV) \bigr) \; ,
\]
\[
E^{p,s}_2 = \mu_{\pi_1(K_1)} \circ H^p \bigl( H_C/H_W, H^s ( W_1, \Res^G_{Q_1}\BV) \bigr) \; ,
\]
and
\[
E^{p,s}_{\partial,2} = \mu_{\pi_1(K_1)} \circ \partial H^p \bigl( H_C/H_W, H^s ( W_1, \Res^G_{Q_1}\BV) \bigr) \; ,
\] 
these extensions are the morphisms $E^{p,s}_{c,2} \to E^{p,s}_2$,
$E^{p,s}_2 \to E^{p,s}_{\partial,2}$, and $E^{p,s}_{\partial,2} \to E^{p+1,s}_{c,2}$
induced by the natural transformations of cohomological functors
$H^{\bullet}_c \to H^{\bullet}$, $H^{\bullet} \to \partial H^{\bullet}$, and $\partial H^{\bullet} \to H^{{\bullet}+1}_c$,
respectively.
\end{Cor}

\begin{Thm} \label{8U}
The three spectral sequences of Corollary~\ref{8T} dege\-ne\-rate and split canonically
in a compatible way. More precisely, let $\BV \in D^+ ( \Rep_F G )$. \\[0.1cm]
(a)~For any $n \in \BZ$, there are canonical and
functorial isomorphisms of local systems on the space
$M^K (Q_1;P_1(\BA_f)gK,\FX)$ between
\[
\CH^{n} \bigl( Ri^!_{Q_1;P_1(\BA_f)gK} i^* R(i_G)_* \bigr) \circ \mu_K (\BV)
\]
and 
\[
\bigoplus_{p+s = n} \mu_{\pi_1(K_1)} \circ 
     H^p_c \bigl( H_C/H_W, H^s ( W_1, \Res^G_{Q_1}\BV) \bigr) \; ,
\]
between 
\[
\CH^{n} \bigl( i^*_{Q_1;P_1(\BA_f)gK} i^* R(i_G)_* \bigr) \circ \mu_K (\BV)
\]
and 
\[
\bigoplus_{p+s = n} \mu_{\pi_1(K_1)} \circ 
     H^p \bigl( H_C/H_W, H^s ( W_1, \Res^G_{Q_1}\BV) \bigr) \; ,
\]
and between
\[
\CH^{n} \bigl( i^*_{Q_1;P_1(\BA_f)gK} R (h_{Q_1})_* h_{Q_1}^* i^* R(i_G)_* \bigr) \circ \mu_K (\BV)
\]
and 
\[
\bigoplus_{p+s = n} \mu_{\pi_1(K_1)} \circ 
     \partial H^p \bigl( H_C/H_W, H^s ( W_1, \Res^G_{Q_1}\BV) \bigr) \; .
\]
(b)~The isomorphism between
\[
\CH^{n} \bigl( i^*_{Q_1;P_1(\BA_f)gK} i^* R(i_G)_* \bigr) \circ \mu_K (\BV)
\]
and 
\[
\bigoplus_{p+s = n} \mu_{\pi_1(K_1)} \circ 
     H^p \bigl( H_C/H_W, H^s ( W_1, \Res^G_{Q_1}\BV) \bigr) 
\]
is identical to the isomorphism from Theorem~\ref{7Q}~(a), for $\uQ := (Q_1)$. \\[0.1cm]
(c)~The three types of morphisms 
in the long exact cohomology sequence associated to
the restriction to $M^K (Q_1;P_1(\BA_f)gK,\FX)(\BC)$ of the exact triangle
\[
R i_{Q_1}^! i^*R(i_G)_* \mu_K(\BV) \to  
i_{Q_1}^* R(i_G)_* \mu_K(\BV) \to
i_{Q_1}^* R (h_{Q_1})_* h_{Q_1}^* i^*R(i_G)_* \mu_K(\BV) \stackrel{[1]}{\to}  
\]  
respect the bigradings induced by the isomorphisms from (a)~: their restriction
to the respective summands indexed by $(p,s)$ yields the morphisms identified in Corollary~\ref{8T}~(c), \emph{i.e.},
those induced by the natural transformations of cohomological functors
$H^{\bullet}_c \to H^{\bullet}$, $H^{\bullet} \to \partial H^{\bullet}$, and $\partial H^{\bullet} \to H^{{\bullet}+1}_c$,
respectively.
\end{Thm}

\begin{Proof}
We repeat the proof of Theorem~\ref{7Q}: since $\bar{Q}_1$ is reductive,
there is a canonical and 
functorial isomorphism in $D^+ (\Rep_F \bar{Q}_1)$
\[
R\Gamma(W_1, \BX) \isoto
                  \bigoplus_{s \in \BZ} H^s (W_1, \BX)[-s] 
\]
for any $\BX \in D^+ (\Rep_F Q_1)$.
\end{Proof}

\begin{Rem}
The sub-space $M^K (Q_1;P_1(\BA_f)gK,\FX)$ of $M^K  (G,\FX)^*$ being algebraic, the restriction 
to $M^K (Q_1;P_1(\BA_f)gK,\FX)$ of the exact triangle
\[
R i_{Q_1}^! i^*R(i_G)_* \mu_K(\BV) \to  
i_{Q_1}^* R(i_G)_* \mu_K(\BV) \to
i_{Q_1}^* R (h_{Q_1})_* h_{Q_1}^* i^*R(i_G)_* \mu_K(\BV) \stackrel{[1]}{\to}  
\]
underlies an exact triangle of objects in the derived category of \emph{mixed Hodge modules}. 
It appears reasonable to expect
Theorem~\ref{8U} itself to admit a Hodge theoretic variant, generalizing the main result from \cite{BW}
(which covers the middle terms
\[
\CH^{n} \bigl( i^*_{Q_1;P_1(\BA_f)gK} i^* R(i_G)_* \bigr) \circ \mu_K (\BV)  \; ) .
\]
An analogous remark holds for \emph{$\ell$-adic sheaves}, where one may expect a ge\-neralization 
of the main result from \cite{P2}.
\end{Rem} 

\begin{Rem} \label{8V}
Example~\ref{7R} (concerning Siegel threefolds) can be re-interpreted in the framework of Theorem~\ref{8U}.
Recall that we
fix a four-dimensio\-nal $\BQ$-vector space $V$, together with a $\BQ$-valued
non-degenerate symplectic bilinear form $J$, and set $G$ equal to the group
\[
GSp(V,J) \subset GL (V) 
\]
of symplectic similitudes of $V$. With $Q_1 \prec Q_2$ as in Example~\ref{7R}, we have
\[
\overline{M^K \bigl( Q_1,\FX \bigr)} = M^K \bigl( Q_1,\FX \bigr)
\]
(as $P_1$ is solvable) and
\[
\partial M^K  (G,\FX)^* - M^K \bigl( Q_1,\FX \bigr) = M^K \bigl( Q_2,\FX \bigr)
\]
(as any two totally isotropic sub-spaces of $V$ of the same dimension are conjugate under $G(\BQ)$).

It follows that
\[
h_{Q_1} = i_{Q_2} : M^K \bigl( Q_2,\FX \bigr) \longinto \partial M^K  (G,\FX)^* \; .
\]
\emph{A fortiori},
\[
i^*_{Q_1;P_1(\BA_f)gK} R (i_{Q_2})_* i_{Q_2}^* i^* R(i_G)_*
= i^*_{Q_1;P_1(\BA_f)gK} R (h_{Q_1})_* h_{Q_1}^* i^* R(i_G)_* \; . 
\]
It follows that the adjunction
\[
\CH^\ast \bigl( i^*_{Q_1;P_1(\BA_f)gK} i^* R(i_G)_* \bigr) 
\longto
\CH^\ast \bigl( i^*_{Q_1;P_1(\BA_f)gK} R(i_{Q_2})_* 
i_{Q_2}^* i^* R(i_G)_* \bigr) 
\]
considered in Example~\ref{7R}~(c) is part of a long exact cohomology sequence, whose ``third term'' equals
\[
\CH^\ast \bigl( Ri^!_{Q_1;P_1(\BA_f)gK} i^* R(i_G)_* \bigr) \; .
\]
Theorem~\ref{8U} relates the latter to cohomology with compact supports of $H_C/H_W$. 
In particular, the kernel of adjunction is related to what could be called
\emph{interior cohomology} of $H_C/H_W$, \emph{i.e.}, the image of cohomology with compact supports
in cohomology of $H_C/H_W$. In the context of \cite[Rem.~2.10~(b)]{W2}, this observation yields a direct 
link between the kernel of adjunction and the space of cusp forms (cmp.~Remark~\ref{7S}~(a)).  
\end{Rem}


\bigskip

\bigskip

\bigskip

%
%

\end{document}